\documentstyle[11pt]{article}

\begin{document}
\noindent
\begin{center}
  {\LARGE Virtual neighborhoods  and pseudo-holomorphic curves}
  \end{center}

  \noindent
  \begin{center}
    {\large Yongbin Ruan}\footnote{partially supported by a NSF grant and a
    Sloan
    fellowship}\\[5pt]
      Department of Mathematics, University of Wisconsin-Madison\\
	Madison, WI 53706\\[5pt]
	 
	      \end{center}

	      \def \J{{\cal J}}
	      \def \Map{Map(S^2, V)}
	      \def \M{{\cal M}}
	      \def \A{{\cal A}}
      \def \B{{\cal B}}
	      \def \C{{\bf C}}
	      \def \Z{{\bf Z}}
	      \def \R{{\bf R}}
	      \def \P{{\bf P}}
	      \def \I{{\bf I}}
	      \def \N{{\cal N}}
	      \def \T{{\bf T}}
	      \def \Q{{\bf Q}}
	      \def \D{{\cal D}}
	      \def \H{{\cal H}}
	      \def \S{{\cal S}}
	      \def \e{{\bf E}}
	      \def \C{{\bf C}}
               \def \CP{{\cal CP}}
	      \def \U{{\cal U}}
	      \def \E{{\cal E}}
	      \def \F{{\cal F}}
	      \def \L{{\cal L}}
	      \def \K{{\cal K}}
              \def \G{{\cal G}}
              \def \V{{\cal V}}
	      \section{ Introduction}
   Since Gromov introduced his pseudo-holomorphic curve theory in the 80's, 
 pseudo-holomorphic curve
has soon become an eminent technique in symplectic topology. Many important 
theorems in
this field have been proved by this technique, among them, the squeezing theorem \cite{Gr},
the rigidity theorem \cite{E}, the classification of rational and ruled symplectic
4-manifolds \cite{M2}, the proof of the existence of non-deformation equivalent symplectic structures \cite{R2}. The pseudo-holomorphic curve also plays
a critical role in a number of new subjects such as Floer homology theory,etc. 

 In the meantime of this development,  a great deal of efforts 
has been
made to solidify the foundation of pseudo-holomorphic curve theory, for examples, McDuff's 
transversality theorem for ``cusp curves'' \cite{M1} and the various proofs of Gromov 
compactness theorem.  In the early day of Gromov theory, 
Gromov compactness theorem was enough for its applications to symplectic 
topology. However, 
it was insufficient for its potential applications in algebraic geometry, where a good compactification is often
very important. For example, it is particularly desirable to tie Gromov-compactness
theorem to the Deligne-Mumford compactification of the moduli space of stable curves. Gromov's
original proof is geometric. Afterwards, many works were done to 
prove
Gromov compactness theorem in the line of Uhlenbeck bubbling off. It was succeed by
Parker-Wolfson \cite{PW} and Ye \cite{Ye}. One outcome of their work was a more
delicate compactification of the moduli space of pseudo-holomorphic maps. But 
it didn't attract much attention 
until several years later when Kontsevich and Manin \cite{KM} rediscovered  this new compactification 
in algebraic geometry and initiated an algebro-geometric approach to the same theory.
Now this new compactification becomes  known as the moduli space of
stable maps. The moduli space of stable maps is one of the basic ingredients of this paper.

During last several years, pseudo-holomorphic curve theory entered a period of rapid expansion.
 We has
witnessed its intensive interactions with algebraic geometry, mathematical physics and 
recently with new Seiberg-Witten theory of 4-manifolds \cite{T2}.  One should mention that
those recent activities in pseudo-holomorphic curve theory did not come from the internal
drive of symplectic topology. It was influenced mostly by mathematical physics,
particularly, Witten's theory of topological sigma model. Around 1990, there were
many activities in string theory about ``quantum cohomology'' and mirror symmetry.
The core of quantum cohomology theory is so called ``counting the numbers of rational
curves''. Many incredible predictions were made about those numbers in Calabi-Yau
3-folds, based on results from physics.
But mathematicians were frustrated about the meaning of the so-called ``number of rational
curves''. For example, the finiteness of such number
 is a well-known conjecture
due to H. Clemens which concerns  simplest Calabi-Yau 3-folds-Quintic hypersurface of $\P^4$. It was even worse
that some Calabi-Yau 3-fold never has a finite number of rational curves. 
 One of the basic difficulties at that time was that people usually 
restricted
their attention to Kahler manifolds, where the complex structure is rigid. On the
other hand, the advantage of pseudo-holomorphic curves is that we are allowed
 to choose
almost complex structures, which are much more flexible. Unfortunately,
the most of those exciting developments were little known to symplectic topologists. 
In \cite{R1}, the author brought the machinery of pseudo-holomorphic curves into quantum
cohomology and mirror symmetry. Using  ideas from Donaldson theory, the author provided
a rigorous definition of the symplectic invariants corresponding the "numbers of rational curves" in
string theory. Moreover, the author found many applications of  new symplectic invariants
in symplectic topology \cite{R1}, \cite{R2}, \cite{R3}. These new invariants
are  called ``Gromov-Witten invariants''.

  Gromov-Witten invariants are analogous of  invariants in the enumerative
geometry. However, the actual counting problem (like the numbers of higher degree rational
curves in 
quintic three-fold) did not attract much of attention before the discovery of mirror
symmetry. In general, these numbers are difficult to compute. Moreover,
 computing these number
didn't  seems to help our understanding of Calabi-Yau 3-folds themself. The introduction of
quantum cohomology hence opened a new direction for enumerative geometry. 
According to quantum cohomology theory, these enumerative invariants are not isolated numbers; instead, they are
 encoded in a new cohomological structure of underline manifold.
Note that the quantum cohomology structure is governed by the
associativity law, which corresponds to the famous composition law of topological quantum
field theory.  Therefore, it would be
very important to put quantum cohomology in a rigorous mathematical foundation. 
It was clear that
the enumerative geometry is not a correct framework. (For example, the associativity or composition law of
quantum cohomology
computes certain higher genus invariants, which are always different from enumerative invariants). 
 Based  on \cite{R1}, a correct mathematical framework were layed down by
the author and Tian \cite{RT1}, \cite{RT2} in terms of perturbed holomorphic 
maps. By proving the
crucial associativity law, we put quantum cohomology in a solid mathematical ground. A corollary of the proof of associativity law is a computation of
the number of rational curves in $\P^n$ and many Fano manifolds by recursion 
formulas. Such a formula
for $\P^2$ was first derived by Kontsevich, based on associativity law 
predicted by physicists.
 It should be pointed out that the entire
pseudo-holomorphic curve theory were only established for so-called semi-positive
symplectic manifolds. They includes  most of interesting
examples like Fano and Calabi-Yau manifolds. But, semipositivity is a 
significant obstacle for some
of its important applications like Arnold conjecture and birational geometry.

  Stimulated by the success of symplectic method, the progreses have been made on algebro-geometrical approach.  An important step  is 
Kontsevich-Manin's initiative of using stable (holomorphic) maps. The genus 0  stable
map works nicely for homogeneous space. For example, the moduli spaces of genus
$0$ stable maps always have expected dimension. Many of results in \cite{R1},
\cite{RT1} were redone in this category by \cite{KM},  \cite{LT1}. It 
was soon realized that moduli spaces of stable maps no longer have expected dimension for
non-homogeneous spaces, for example, projective bundles \cite{QR}. To go
beyond homogeneous spaces, one needs new ideas. A breakthrough came with the work
of  Li and Tian \cite{LT2}, where they employ a sophisticated excessive 
intersection theory (normal cone construction) (see another proof in \cite{B}). As a consequence, Li and Tian extended GW-invariant to 
arbitrary algebraic manifolds. In the light of these new developments,  
three obvious problems have emerged: (i) to remove
semi-positivity condition in Gromov-Witten invariants; (ii) to remove semi-positive condition
in Floer homology and solve Arnold conjecture. (iii) to prove that
symplectic GW-invariants are the same as algebro-geometric GW-invariants for algebraic 
manifolds. We will deal with first two problems in this article and leave the
last one to the future research. 

    Recall that, the fundamental difficulty for pseudo-holomorphic curve theory
on non-semi-positive symplectic manifolds is,
 that $\overline{\M}-\M$ may have
arger dimension than that of $\M$, where $\M$ is the moduli space of pseudo-holomorphic
maps and $\overline{\M}$ is a compactification. One view is that this is due to
the reason that the almost complex structure is not
generic at infinity. To deal with this non-generic situation, the  author's idea \cite{R3} (Proposition 5.7)
was to construct
an open smooth manifold (virtual neighborhood ) to contain the moduli space.
Then, we can work on virtual neighborhood, which is much easier to handle 
than the
moduli space itself. In \cite{R4}, the author outlined a scheme to attack the non-generic
problems in Donaldson-type theory using virtual neighborhood technique. Moreover, author 
applied virtual neighborhood technique to monopole 
equation  under a group action. Further application can be found in
\cite{RW}. But the case in \cite{R4} is too
restricted for pseudo-holomorphic case. Recall that in \cite{R4}, we work with
a compact-smooth triple $(\B, \F, S)$ where $\B$ is a smooth Banach manifold
(configuration space), $\F$ is a smooth Banach bundle and $S$ is a section of
$\F$ such that the moduli space 
$\M=S^{-1}(0)$ is compact.
Monopole equation can be interpreted as a smooth-compact triple. However, in
the case of 
pseudo-holomorphic curve, $S^{-1}(0)$ is almost never compact in the 
configuration space. Furthermore, $(\B, \F)$ is often not smooth, but a pair
of $V$-manifold and $V$-bundle. To overcome these difficulty, we need to
generalize the virtual neighborhood technique to handle this  situation. An 
outline of such a generalization were given in \cite{R4}. 

  Another  purpose of this paper is to construct an equivariant quantum 
cohomology theory. For this purpose, we need to study the GW-invariant for
a family of symplectic manifolds. We shall work in
this generality throughout the paper.  Let's outline a definition of 
GW-invariant over a family of symplectic manifolds as follows.
  
  Let $P: Y\rightarrow X$ be a fiber bundle such that both the fiber $V$ and the base
$X$ are smooth compact, oriented manifolds. Furthermore, we assume that $P: Y
\rightarrow X$ is an oriented fibration. Then, $Y$ is also a smooth, compact,
oriented manifold. Let $\omega$ be a closed 2-form on $Y$ such that $\omega$
restricts to a symplectic form over each fiber. A $\omega$-tamed almost complex
structure $J$ is an automorphism of vertical tangent bundle such that $J^2=-Id$
and $\omega(X, JX)>0$ for vertical tangent vector $X\neq 0$. Let $A\in H_2(V,
\Z)\subset H_2(Y, \Z)$. Let $\M_{g,k}$ be the moduli space of 
genus g Riemann surfaces with $k$-marked points such that $2g+k>2$ and 
$\overline{\M}_{g,k}$ be its Deligne-Mumford compactification. Suppose that 
$f: \Sigma\rightarrow Y$ ($\Sigma\in \M_{g,k}$) is a smooth map such that 
$im(f)$ is contained in a fiber and $f$ satisfies Cauchy Riemann equation 
$\partial_J f=0$ with $[f]=A$. Let 
$\M_A(Y, g,k,J)$ be the moduli space of such $f$.  First we need a stable 
compactification of $\M_{A}(Y,g,k, J)$. Roughly speaking, {\em a 
compactification is stable if its local Kuranish model is the quotient of 
vector spaces by a finite group}. In our case, it is provided by the 
moduli space of stable holomorphic maps $\overline{\M}_{A}(Y,g,k,J)$. 

  There are two technical difficulties to use virtual neighborhood technique
to the case of pseudo-holomorphic curve. The first one is that there is a
 finite group action on its
local Kuranish model. An indication is that we should work in the V-manifold 
 and V-bundle category. As a matter of fact, it is easy to
extend virtual neighborhood technique to this category. However, the finite 
dimensional virtual
neighborhood constructed is a V-manifold in this case. It is well-known that
the ordinary transversality theorem fails for V-manifolds. We will overcome 
this problem by using
differential form and integration. We shall give a detail argument in
section 2. The second problem is  the failure of the compactness of
$\M_A(Y,g,k)$. To include $\overline{\M}_A(Y,g,k)$, we have to enlarge our
configuration space to $\overline{\B}_A(Y,g,k)$ of $C^{\infty}$-stable (
holomorphic or not) maps. Then, the obstruction bundle $\F_A(Y,g,k)$
extends to $\overline{\F}_A(Y,g,k)$ over $\overline{\B}_A(Y,g,k)$.
Therefore, we obtained a compact triple $(\overline{\B}_A(Y,g,k), \overline{\F}
_A(Y,g,k), \S)$, where $S$ is Cauchy-Riemann equation. We want to generalize
the virtual neighborhood technique to this enlarge space. Recall that for virtual
neighborhood technique, we construct some stablizations of the equation
 $\S_e=\S+s$, which must satisfy two crucial properties:
 (1) $\{x; Coker \delta_x(\S+s)=0\}$
is open; (2)If $\S+s$ is a transverse section, $U=(\S+s)^{-1}(0)$ is a finite
dimensional smooth V-manifold. By using gluing argument, we can construct a
local
model of $U$ (local Kuranish model). (2) is equivalent to  that the local 
Kuranish model is a quotient of vector spaces by a finite group. By 
definition, it means that our compactification has to be stable. Finally, we 
need
an additional argument to prove that the local models patch together smoothly.
We call a triple satisfying (1), (2) {\em virtual neighborhood technique 
admissible} or {\em VNA}.

   Suppose that $\S$ is already transverse. $\overline{\M}(Y,g,k)$ is naturally a
stratified space whose stratification coincides with that of $\overline{\B}_A
(Y,g,k)$. The attaching map of $\overline{\B}_A(Y,g,k)$ is defined by patching
construction. The gluing theorem shows that if we restrict ourself to
stable holomorphic maps  one can deform this attaching map
slightly such that the image of stable holomorphic maps is again holomorphic.
The deformed attaching map gives a local smooth coordinate
of $\overline{\M}_A(Y,g,k)$. Although it is not necessary in virtual
neighborhood construction, one can also attempt to deform the whole attaching
map by the same implicit funtion theorem argument. Then, it is attempting to
think (as author did) that the deformed attaching map will give a smooth
coordinate of $\overline{\B}_A(Y,g,k)$. It was Tian who pointed out the author
that this is false. However, it is natural to ask if there is any general
property for such an infinite dimensional object. Indeed, some elegant 
properties are formulated by Li and Tian  \cite{LT3} and we refer reader to 
their paper for the detail. 

   Applying virtual neighborhood technique, we construct a finite dimensional
  virtual
neighborhood $(U,F, S)$. More precisely, $U$ is covered by finite many
coordinate charts of the form $U_i/G_i$ ($i=1, \dots, m$) for $U_p \subset
 \R^{ind 
+m}$ and a finite group $G_p$. $F$ is a V-bundle over $U$ and $S: U\rightarrow F$
is a section.
 On the other hand, the evaluation maps over
marked points define a map 
$$\Xi_{g,k}: \overline{\B}_A(Y,g,k) \rightarrow Y^k. \leqno(1.1)$$
We  have another map 
$$\chi: \overline{\B}_A(Y,g,k) \rightarrow \overline{\M}_{g,k}.\leqno(1.2)$$
Recall that $\overline{\M}_{g,k}$ is a V-manifold. To define GW-invariant,
choose  a Thom form $\Theta$ supported in a neighborhood of zero section.
The GW-invariant can be defined as
$$\Psi^Y_{(A,g,k)}(K; \alpha_1, \cdots, \alpha_k)=\int_{U}\chi^*(K)\wedge \Xi^*_{g,k}
\prod_i \alpha_i\wedge S^*\Theta. \leqno(1.4)$$
for $\alpha_i\in H^*(Y, \R)$ and $K\in H^*(\overline{\M}_{g,k}, \R)$ 
represented by differential form. 
Clearly, $\Psi^Y=0$ if $\sum deg(\alpha_i)+deg (K)\neq ind$.

  Recall that $H^*(Y, \R)$ has a modular structure by $P^*\alpha$ for 
$\alpha\in H^*(X, \R)$. In this paper, we prove the following,
\vskip 0.1in
\noindent
{\bf Theorem A (Theorem 4.2): }{\it (i).$\Psi^Y_{(A,g,k)}(K; \alpha_1, \cdots, \alpha_k)$ is
well-defined.
\vskip 0.1in
\noindent
(ii). $\Psi^Y_{(A,g,k)}(K; \alpha_1, \cdots, \alpha_k)$ is independent of  the choice of virtual neighborhoods.
\vskip 0.1in
\noindent
(iii). $\Psi^Y_{(A,g,k)}(K; \alpha_1, \cdots, \alpha_k)$ is independent of $J$ and
 is a symplectic deformation invariant.
\vskip 0.1in
\noindent
(iv). When $Y=V$ is semi-positive, $\Psi^Y_{(A,g,k)}(K; \alpha_1, \cdots, 
\alpha_k)$ agrees with the definition of \cite{RT2}.
\vskip 0.1in
\noindent
(v). $\Psi^Y_{(A,g,k)}(K; \alpha_1, \cdots, \alpha_i\cup P^*\alpha,\cdots, 
\alpha_k)=\Psi^Y_{(A,g,k)}(K; \alpha_1, \cdots, \alpha_j \cup P^*\alpha, 
\cdots, \alpha_k)$}
\vskip 0.1in
Furthermore, we can show that $\Psi$ satisfies the composition law required by
the theory of sigma model coupled with gravity.
Assume $g=g_1+g_2$ and $k=k_1+k_2$ with $2g_i + k_i \ge 3$.
Fix  a decomposition $S=S_1\cup S_2$ of $\{1,\cdots , k\}$ with
$|S_i|= k_i$. Then there is a canonical embedding
$\theta _S: \overline \M_{g_1,k_1+1}\times \overline \M_{g_2,k_2+1}
\mapsto \overline \M_{g,k}$, which assigns to marked
curves $(\Sigma _i; x_1^i,\cdots ,x_{k_1+1}^i)$ ($i=1,2$), their 
union $\Sigma _1\cup \Sigma _2$ with $x^1_{k_1+1}$ identified to
$x^2_{k_2+1}$ and remaining points renumbered by $\{1,\cdots,k\}$ according to $S$.

There is another natural map $\mu : \overline \M_{g-1, k+2}
\mapsto \overline \M_{g,k}$ by gluing together the last two marked 
points.

Choose a homogeneous 
basis $\{\beta _b\}_{1\le b\le L}$ of $H_*(Y,\Z)$ modulo
torsion. Let $(\eta _{ab})$ be its intersection matrix. Note that
$\eta _{ab} = \beta _a \cdot \beta _b =0$ if the dimensions of
$\beta _a$ and $\beta _b$ are not complementary to each other.
Put $(\eta ^{ab})$ to be the inverse of $(\eta _{ab})$.
Now we can state the composition law, which consists of 
two formulas as follows.
\vskip 0.1in
\noindent
{\bf Theorem B. (Theorem 4.7)}  {\it Let $[K_i] \in H_*(\overline \M_{g_i,
k_i+1}, \Q)$ $(i=1,2)$ and $[K_0] \in H_*(\overline \M_{g-1,
k +2}, \Q)$. For any $\alpha _1,\cdots,\alpha _k$ in $H_*(V,\Z)$.
Then we have

$$\begin{array}{rl}
&\Psi ^Y_{(A,g,k)}(\theta _{S*}[K_1\times K_2];\{\alpha _i\})\\
=(-1)^{deg(K_2)\sum_{i=1}^{k_1}deg(\alpha_i)} ~& \sum \limits _{A=A_1+A_2} \sum \limits_{a,b}
\Psi ^Y_{(A_1,g_1,k_1+1)}([K_1];\{\alpha _{i}\}_{i\le k_1}, \beta _a) 
\eta ^{ab}
\Psi ^Y_{(A_2,g_2,k_2+1)}([K_2];\beta _b,
\{\alpha _{j}\}_{j> k_1}) \\
\end{array}
\leqno (1.5) 
$$
$$
\Psi ^Y_{(A,g,k)}(\mu_*[K_0];\alpha _1,\cdots, \alpha _k)
=\sum _{a,b} \Psi ^Y_{(A,g-1,k+2)}([ K_0];\alpha _1,\cdots, \alpha _k,
\beta _a,\beta _b) \eta ^{ab}\leqno (1.6)
$$
}
\vskip 0.1in

There is a natural map $\pi:
\overline{\M}_{g,k}\rightarrow \overline{\M}_{g, k-1}$ as follows: For 
$(\Sigma, x_1, \cdots, x_k)\in \overline{\M}_{g,k}$, if $x_k$ is not in any 
rational component of $\Sigma$ which contains only three special points,
then we define
$$\pi(\Sigma, x_1, \cdots, x_k)=(\Sigma, x_1, \cdots, x_{k-1}),$$
where a distinguished point of $\Sigma$ is either a singular point or a
marked point. If $x_k$ is in one of such rational components, we contract
this component and obtain a stable curve $(\Sigma', x_1, \cdots, x_{k-1})$ in
$\overline{\M}_{g, k-1}$, and define $\pi(\Sigma, x_1, \cdots, x_k)=(\Sigma',
x_1, \cdots, x_{k-1}).$ 

Clearly, $\pi$ is continuous. One should be aware that
there are two exceptional cases $(g,k)=(0,3), (1,1)$ where $\pi$ is not well
defined. Associated with $\pi$,
we have two {\em k-reduction formula} for $\Psi^V_{(A, g, k)}$ as following:
\vskip 0.1in
\noindent
{\bf Proposition C (Proposition 4.4). }{\it Suppose that $(g,k)
\neq (0,3),(1,1)$.
\vskip 0.1in
\noindent
(1) For any $\alpha _1, \cdots , \alpha _{k-1}$ in $H_*(Y, \Z)$, 
we have} 
$$\Psi ^Y_{(A,g,k)}(K; \alpha _1, \cdots,\alpha _{k-1}, [V])~=~
\Psi ^Y_{(A,g,k-1)}([\pi_* (K)]; \alpha _1, \cdots,\alpha _{k-1})
\leqno (1.7)$$
\vskip 0.1in
\noindent
(2)  Let $\alpha _k$ be in $H_{2n-2}(Y, \Z)$, then
$$\Psi ^Y_{(A,g,k)}(\pi^{*}(K); \alpha _1, 
\cdots,\alpha _{k-1}, \alpha _k)~=~\alpha^* _k (A)
\Psi ^Y_{(A,g,k-1)}(K; \alpha _1, \cdots,\alpha _{k-1})
\leqno (1.8)$$
where $ \alpha^* _k$ is the Poincare dual of $\alpha _k$.
\vskip 0.1in
\noindent
When $Y=V$, $\Psi^Y$ is the ordinary GW-invariants. Therefore, we establish a theory of topological sigma model couple with
gravity over any symplectic manifolds.

  It is well-known that GW-invariant can be used to define a quantum 
multiplication. Let's briefly sketch it as follows. First we define a total 3-point function
$$\Psi^V_{\omega}(\alpha_1, \alpha_2, \alpha_3)=\sum_A \Psi^V_{(A,0,3)}(pt;
\alpha_1, \alpha_2, \alpha_3)q^A, \leqno(1.9)$$
where $q^A$ is an element of Novikov ring $\Lambda_{\omega}$ (see \cite{RT1}, 
\cite{MS}). Then, we define a quantum multiplication $\alpha\times_Q \beta$
over  $H^*(V, \Lambda_{\omega})$ by the relation
$$(\alpha\times_Q \beta)\cup \gamma[V]=\Psi^V_{\omega}(\alpha_1, \alpha_2, \alpha_3),\leqno(1.10)$$
where $\cup$ represents the ordinary cup product. As a consequence of Theorem 
B, we have
\vskip 0.1in
\noindent
{\bf Proposition D: }{\it Quantum multiplication is associative over any 
symplectic manifolds. Hence, there is a quantum ring structure over any 
symplectic manifolds.}
\vskip 0.1in
  Given a periodic Hamiltonian function $H: S^1\times V\rightarrow V$, we can
define the Floer homology $HF(V, H)$, whose chain complex is generated by the
periodic  orbits of $H$ and the boundary maps are defined by the moduli spaces 
of 
flow lines. So far, Floer homology $HF(V, H)$ is only defined for semi-positive
symplectic manifolds.  Applying virtual neighborhood technique to Floer 
homology, we show
\vskip 0.1in
\noindent
{\bf Theorem E: }{\it Floer homology $HF(V, H)$ is well-defined for any
symplectic manifolds. Furthermore, $HF(V, H)$ is independent of $H$.}
\vskip 0.1in
Recall that Floer homology was invented to solve the
\vskip 0.1in
\noindent
{\bf Arnold conjecture: }{Let $\phi$ be a non-degenerate Hamiltonian symplectomorphism. Then,
the number of the fixed points of $\phi$ is greater than or equal to the sum of
Betti number of $V$.}
\vskip 0.1in
As a corollary of Theorem E, we prove the Arnold conjecture
\vskip 0.1in
\noindent
{\bf Theorem F: }{\it Arnold conjecture holds for any symplectic manifolds.}
\vskip 0.1in
In this paper, we give another application of our results in higher dimensional
algebraic  geometry. It
was discovered in \cite{R3} that symplectic geometry has a strong connection
with Mori's birational geometry. An important notion in birational geometry
is uniruled variety, generalizing the notion of ruled surfaces in two
 dimension.  An algebraic variety $V$ is uniruled iff $V$ is covered by 
rational curves. Kollar \cite{K1} proved that for 3-folds, uniruledness is a 
symplectic property. Namely, if a 3-fold  $W$ is
symplectic deformation equivalent to an uniruled variety $V$, $W$ is uniruled.
To extend Kollar's result, we need a symplectic GW-invariants defined over any
symplectic manifolds with certain property (Lemma 4.10). We will show that our 
invariant satisfies this properties.  
By combining with Kollar's result, we have
\vskip 0.1in
\noindent
{\bf Proposition G: }{\it If a smooth Kahler manifold $W$ is symplectic 
deformation equivalent to a uniruled variety, $W$ is uniruled.}
\vskip 0.1in
An important topic in quantum cohomology theory is the equivariant quantum 
cohomology group $QH_G(V)$, which generalizes the notion of equivariant cohomology. 
Suppose that a compact Lie group $G$ acts on $V$ as symplectomorphisms. To
define equivariant quantum cohomology, we first have to define equivariant
GW-invariants. There are two approaches. The first approach is to choose a
$G$-invariant tamed almost complex structure $J$ and construct an equivariant
virtual neighborhood. Then, we can use finite dimensional equivariant technique
to define equivariant GW-invariant. This approach indeed works. But a 
technically simpler approach is to consider equivariant GW-invariant as the
limit of GW-invariant over the families of symplectic manifolds. This approach was
advocated by Givental and Kim \cite{GK}. We shall use this approach here.

Let $BG$ be the classifying space of $G$ and $EG\rightarrow BG$ be the 
universal $G$-bundle. Suppose that 
$$BG_1\subset BG_2\cdots\subset BG_m \subset BG \leqno(1.11)$$
such that $BG_i$ is a smooth oriented compact manifold and $BG=\cup_i BG_i$. Let
$$EG_1\subset EG_2\cdots\subset EG_m \subset BG\leqno(1.12)$$
be the corresponding universal bundle. We can also form the approximation of
homotopy quotient $V_G=V\times EG/G$ by $V^i_G=V\times EG_i/G$. Since $\omega$
is invariant under $G$, its pull-back on $V\times EG_i$ descends to $V^i_G$.
So, we have a family of symplectic manifolds $P_i: V^i\rightarrow BG_i$. 
Applying our previous construction, we obtain GW-invariant $\Psi^{P_i}_{(A, 
g,k)}$. We define equivariant GW-invariant 
$$\Psi^G_{(A,g,k)}\in Hom((H^*(V_G, \Z))^{\otimes k}\otimes H^*(\overline{\M}_{g,
k}, \Z), H^*(BG, \Z)) \leqno(1.13)$$
as follow:

  For any $D\in H_*(BG, \Z)$, $D\in H_*(BG_i, \Z)$ for some $i$.
Let $i_{V^i_G}: V^i_G\rightarrow  V_G$. For $\alpha_i\in H^*_G(V)$, we define
$$\Psi^G_{(A,g,k)}(K, \alpha_1, \cdots, \alpha_k)(D)=\Psi^{P_i}_{(A, g,k)}(K,
i^*_{V^i_G}(\alpha_1), \cdots, i^*_{V^i_G}(\alpha_k); P^*_i(D^*_{BG_i})),\leqno(1.14)$$
where $D^*_{BG_i}$ is the Poincare dual of $D$ with respect to $BG_i$.
\vskip 0.1in
\noindent
{\bf Theorem G: }{\it (i). $\Psi^G_{(A, g,k)}$ is independent of the choice of
$BG_i$.
\vskip 0.1in
\noindent
(ii). If $\omega_t$ is a family of $G$-invariant symplectic forms, $\Psi^G_{(A,
g,k)}$ is independent of $\omega_t$.}
\vskip 0.1in
Recall that equivariant cohomology ring $H^*_G(X)$ is defined as $H^*(V_G)$. 
Notes that, for any  equivariant cohomology class $\alpha\in H^*_G(V)$,
$$\alpha [V]\in H^*(BG)\leqno(1.15)$$
instead of being a number in the case of the ordinary cohomology ring. Furthermore, there
is a modulo structure by $H^*_G(pt)=H^*(BG)$, defined by using the projection map
$$V_G\rightarrow BG.\leqno(1.16)$$
The equivariant quantum multiplication is a new multiplication structure
over $H^*_G(V, \Lambda_{\omega})=H^*(V_G, \Lambda_{\omega})$ as follows. We first
define a total 3-point function
$$\Psi^G_{(V,\omega)}(\alpha_1, \alpha_2, \alpha_3)=\sum_A \Psi^G_{(A,0,3)}(
pt; \alpha_1, \alpha_2, \alpha_3)q^A.\leqno(1.17)$$
Then, we define 
an equivariant quantum multiplication by
$$(\alpha\times_{QG}\beta)\cup \gamma [V]=\Psi^G_{(V,\omega)}(\alpha_1, 
\alpha_2, \alpha_3).\leqno(1.18)$$
\vskip 0.1in
\noindent
{\bf Theorem I: }{\it (i) The equivariant quantum multiplication is commutative
with the modulo structure of $H^*(BG)$.
\vskip 0.1in
\noindent
(ii) The equivariant quantum multiplication is skew-symmetry.
\vskip 0.1in
\noindent
(iii) The equivariant quantum multiplication is associative.

Hence, there is a equivariant quantum ring structure for any $G$ and symplectic
manifold $V$}
\vskip 0.1in
 
  Equivariant quantum cohomology has already been defined for monotonic
symplectic manifold by Lu \cite{Lu}.

  The paper is organized as follows: In section 2, we  work out the detail
of the virtual neighborhood technique for Banach V-manifolds. In section 3, we
prove that the virtual neighborhood technique can be applied to pseudo-holomorphic
maps. In the section 4, we 
 prove Theorem A, B, C, D, H and I. We  prove Theorem E, Corollary F in
section 5 and Theorem G in section 6. 

  The results of this paper was announced in a lecture at the IP Irvine conference 
in the end of March, 96.
An outline of this paper was given in \cite{R4}. During the 
preparation of this paper, we received  papers by Fukaya and Ono \cite{FO}, B.
Seibert \cite{S}, Li-Tian \cite{LT3}, Liu-Tian, were informed by Hofer/Salamon 
that they
obtained  some of the results of this paper independently using different methods.
The author would like to thank G. Tian and B. Siebert for pointing out errors in the
first draft and B. Siebert for  suggesting a fix (Appendix) of  an error in Lemma 2.5.
The author would like to thank An-Min Li and Bohui Chen for the valuable discussions.

\section{Virtual neighborhoods for V-manifolds}
  As we mentioned in the introduction, the configuration space 
$\overline{\B}_A(Y,g,k)$ is not a smooth Banach V-manifold in general. But for
the purpose of virtual neighborhood construction, we can treat it as a smooth
Banach V-manifold. To simplify the notation, we will work in the category of
Banach V-manifold in this section and refer to the next section for the proof
that the construction of this section applies to $\overline{\B}_A(Y,g,k,J)$.

 V-manifold is a classical subject dated back at least to \cite{Sa1}.
Let's have a briefly review about the basics of V-manifolds.
\vskip 0.1in
\noindent
{\bf Definition 2.1: }{\it (i).A Hausdorff topological space $M$ is a n-dimensional
V-manifold if for every point $x\in M$, there is an open neighborhood of the form $U_x/G_x$ 
where $U_x$ is a connected open subset of $\R^n$ and $G_x$ is a 
finite group acting on $U_x$ diffeomorphic-ally. Let $p_x: U_x
\rightarrow U_x/G_x$ be the projection. We call $(U_x, G_x, p_x)$ a 
coordinate chart of $x$. If $y\in U_x/G_x$ and $(U_y, G_y, p_y)$ is a coordinate
chart of $y$ such that $U_y/G_y\subset U_x/G_x$, there is an injective
smooth map $U_y\rightarrow U_x$ covering the inclusion $U_y/G_y\rightarrow 
U_x/G_x$. 
\vskip 0.1in
\noindent
(ii). A map between V-manifolds $h: M\rightarrow M'$ is smooth if for every point
$x\in M$, there are local charts $(U_x, G_x,p_x), (U'_{h(x)}, G'_{h(x)}, p'_{
h(x)})$ of $x, h(x)$ such that locally $h$ can be lift to a smooth map 
$$h: U_x\rightarrow U'_{h(x)}.$$
\vskip 0.1in
\noindent
(iii).$P: E\rightarrow M$ is a V-bundle if locally $P^{-1}(U_{\alpha}/
G_{\alpha})$ can be lift to $U_{\alpha}\times \R^k$. Furthermore, the lifting
of a transition map is linear on $\R^k$.
\vskip 0.1in
Furthermore, we can define Banach V-manifold, Banach V-bundle in the same way.}
\vskip 0.1in
An easy observation is that we can always choose a local chart $(U_x,
G_x, p_x)$ of $x$ such that $G_x$ is the stabilizer of $x$  by shrinking the
size of $U_x$. Furthermore, we can assume that $G_x$ acts effectively
and $U_x$ is an open disk neighborhood of the origin $x$ in a linear 
representation $(G_x, \R^n)$. We call such a chart {\em a good chart} and $G_x$ a
{\em local group}. 

Notes that if $S$ is a transverse section of a V-bundle, then $S^{-1}(0)$ is
a smooth V-sub-manifold. But, it is well-known that the ordinary transversality 
theorems fail for V-manifolds.
However, the differential calculus (differential form, orientability, integration,
de Rham theory) extends over V-manifolds. Moreover, the theory of 
characteristic classes and the index theory also extend over $V$ manifolds. We won't
give any detail here. Readers can find a detail expository
in \cite{Sa1}, \cite{Sa2}. In summary, if  we use differential
analysis, we can treat a V-manifold as an ordinary smooth manifold. To simplify
the notation, we will omit the word ``V-manifold'' without confusion when we
work on the differential form and the integration. 

\vskip 0.1in
\noindent
{\bf Definition 2.2: }{\it We call that $M$ to be a fine $V$-manifold if any local 
$V$-bundle is dominated by a global oriented $V$-bundle. Namely, Let $U_{\alpha}\times_{\rho_{
\alpha}} E/G_{\alpha}$ be a local V-bundle, where $\rho_{\alpha}: G_{\alpha}
\rightarrow GL(E)$ is a representation. There is a global oriented V-bundle 
$E\rightarrow M$ such that $U_{\alpha}\times_{\rho_{
\alpha}} E_{\alpha}/G_{\alpha}$ is a subbundle of $E_{U_{\alpha}/G_{\alpha}}$.}
\vskip 0.1in

By a lemma of Siebert (Appendix), $\overline{\B}_A(Y,g,k)$ is fine.

{\em In the rest of the section, we will assume that all the Banach V-manifolds are
fine}

Let $\B$ be a fine Banach V-manifold defined by specifying Sobolov norm of some 
geometric object. Let  $\F\rightarrow \B$ be a Banach V-bundle equipped with a metric
and 
$\S: \B\rightarrow \F$ be a smooth section defined by a nonlinear elliptic 
operator. 
\vskip 0.1in
\noindent
{\bf Definition 2.3: }{\it  $\S$ is a proper section if $\{x; ||\S(x)||\leq C\}$ is
compact for any constant $C$.  We call $\M_S=\S^{-1}(0)$ the moduli space of $F$. We
call  $(\B, \F, \S)$ a compact- V triple if $\B, \F$ is a Banach V-pair and
 $\S$ is proper.}
\vskip 0.1in
When $\S$ is proper, it is clear that $\M_{\S}$ is compact.
\vskip 0.1in
\noindent
{\bf Definition 2.4: }{\it Let $M$ be a compact topological space. We call $(U, E, S)$
 a virtual neighborhood of $M$ if $U$ is a  finite dimensional oriented 
V-manifold (not necessarily compact), $E$ is a finite dimensional V-bundle of 
$U$ and $S$ is a smooth section of $E$ such that $S^{-1}(0)=M$. Suppose 
that $M_{(t)}=\bigcup_t M_t\times \{t\}$ is compact. We call $(U_{(t)}, S_{(t)}, 
E_{(t)})$ a virtual neighborhood cobordism if $U_{(t)}$ is a  finite
dimensional oriented  V-manifold with
boundary and $E_{(t)}$ is a finite dimensional V-bundle and $S_{(t)}$ is a 
smooth section such that $S^{-1}_{(t)}(0)=M_{(t)}$.

  }
\vskip 0.1in
  Let $L_x$ be the linearlization 
$$\delta \S_x: T_x\B\rightarrow \F_x, \leqno(2.12)$$
where the tangent space of a V-manifold at $x$ means the tangent space of 
$U_{\alpha}$ at $x$ where $U_{\alpha}/G_{\alpha}$ is a coordinate chart at $x$.
Then, $L_x$ is an elliptic operator. 
When $Coker L_x=0$ for every $x\in \M$, $\S$ is transverse to the
zero section and $\M_{\S}=\S^{-1}(0)$ is a smooth V-manifold of dimension 
$ind(L_x)$.
The case we are interested in is the case that $Coker L_x\neq 0$ and it may 
even jump the dimension. The original version of following Lemma is erroneous.
The new version is corrected by B. Siebert (appendix).
\vskip 0.1in
\noindent
{\bf Lemma 2.5: }{\it Suppose that $(\B,\F,\S)$ is a compact-V triple. 
There exists an open set $\U$ such that $\M_{\S}\subset 
\U \subset \B$ and a finite dimensional oriented  V-bundle $\E$  over $\U$ with a 
V-bundle map $s: \E\rightarrow \F_{\U}$ such that
$$L_x+s(x, v): T_x \U\oplus \E\rightarrow \F\leqno(2.13)$$ 
is surjective  for any  $x\in \U$. Furthermore, the linearlization of $s$ is a 
compact operator.}
\vskip 0.1in 
{\bf Proof: } For each  $x\in \M_S$, there is a good chart $(\tilde{U}_x, G_x,
p_x)$. Suppose that $\tilde{U}_x$ is open disk of radius 1 in $H$ for some Banach space $H$.
Let  $(\F_{\tilde{U}}, G_x, \pi_x)$ be the corresponding chart of $\F$.
Let $H_x=Coker L_x$. Then, $G_x$ acts on $H_x$. Since $\M_{\S}$ is compact,
there is a finite cover $\{(\frac{1}{2}\tilde{U}_{x_i}, G_{x_i}, p_{x_i})\}^m_1
$. Each $\frac{1}{2}\tilde{U}_{x_i}\times H_{x_i}/G_{x_i}$ is a local V-bundle. Since
$\B$ is fine, there exists an oriented global finite dimensional V-bundle $\E_i$ over $\U=\bigcup_i
\frac{1}{2} U_{x_i}$ such that 
$\frac{1}{2}\tilde{U}_{x_i}\times H_{x_i}/G_{x_i}$ is a subbundle of $(\E_i)|_{\frac{1}{2}
\tilde{U}_{x_i}/G_{x_i}}$. Let
$$\E=\oplus_i \E_i.\leqno(2.14)$$

Next, we define $s$.  Each element $w$ of $H_{x_i}$ can be extended to a local
section of $\F_{\tilde{U}_{x_i}}$. Then one can multiply it
by  a cut-off function $\phi$ such that $\phi=0$ outside of the disk of radius
$\frac{3}{4}$ and $\phi=1$ on $\frac{1}{2}\tilde{U}_{x_i}$. Then, we obtain
a section supported over $\tilde{U}_{x_i}$ (still denoted it by $s$). Define 
$$\bar{s}_i(x, w)=s(x).\leqno(2.15)$$
Then, 
$$s_i(x,w)=\frac{1}{|G_{x_i}|}\sum_{g_i\in G_{x_i}} 
(g_i)^{-1}\bar{s}(g_i(x), g_i(w)). \leqno(2.16)$$
By the construction, $s_i$ descends to a map  $U_{x_i}\times H_{x_i}/G_{x_i}\rightarrow \F_{
U_{x_i}}$. Clearly, $s_i$ can be viewed as a bundle map from $\E_i$ to $\F$
since it  is supported in $U_{x_i}$. Moreover, 
$$s(x_i, w): (\E_i)_{x_i}\rightarrow H_{x_i}\subset \F_{x_i} \leqno(2.17)$$ 
is projection. Then, we define 
$$s=\sum s_i.$$
By (2.17), $L_x+s_i$ is surjective at $x_i$ and hence it is surjective at a
neighborhood of $x_i$. By shrinking $U_{x_i}$, we can assume that $L_x+s_i$
is surjective over $\frac{1}{2}U_{x_i}$. Hence, $L_x+s$ is surjective over
$\U$. We have finished the proof. $\Box$

 Next we define the extended equation 
$$\S_e: \E \rightarrow \F\leqno(2.18)$$
by 
$$\S_e(x, w)=\S(x)+s(x, w) \leqno(2.19)$$
for $w\in E_x$. We call that $s$ {\em a stabilization term} and $\S_e$ {\em a stabilization of $\S$}. 
$\S_e$ can be identified with a section of $\pi^*\F$ where $\pi: \E\rightarrow \U$ is the
projection. We shall use the same $\S_e$ to denote the corresponding section.
Notes that  $\M_{\S}\subset \S^{-1}_v(0)$, where we view $\U$ as the zero section of $\E$. 
Moreover,
its linearlization 
$$(\delta \S_e)_{(x, 0)}(\alpha,  u)=L_{x}(\alpha)+s(x,u). \leqno(2.20)$$
By lemma 2.5, it is surjective. Hence, $\S_e$ is a transverse section
over a neighborhood of $\M_{\S}$. Since we only want to construct a
neighborhood of $\M_{\S}$, without the loss of generality,
we can assume that $\S_e$ is transverse to the zero section of $\pi^*\F$. 
Therefore,
$$U=(\S+ s)^{-1}(0)\subset \E \leqno(2.21)$$
is a smooth V-manifold of dimension $ind(L_x)+dim \E$. Clearly, 
$$\M_{\S}\subset U. \leqno(2.22)$$
\vskip 0.1in
\noindent
{\bf Lemma 2.5: }{\it If $det(L_A)$ has a nowhere vanishing section, it
defines an orientation of $U$.}
\vskip 0.1in
{\bf Proof: } $T_{(x, w)}U=Ker (\delta \S_v)$ and
$Coker (\delta \S_v)=0$ by the construction.
Hence, an orientation of $U$ is equivalent to a nowhere vanishing section of 
$det(ind (\delta \S_v))$. 
$$(\delta \S_v)_{(x,  w)}(\alpha,  u)=L_{x}(\alpha)+s(x,u)+\delta s_{(x,w)}(
\alpha). \leqno(2.13)$$
Let
$$(\delta^t \S_v)_{(x,  w)}(\alpha,  u)=L_{x}(\alpha)+ts(x,u)+t\delta s_{(x,w)}
(\alpha).\leqno(2.14)$$
Then, 
$$det(ind( \delta \S_v))=det(ind (\delta^t \S_v))=det(ind (\delta^0 \S_v))=det(
ind(L_x))\otimes det(\E).$$
Therefore, a nowhere vanishing section of $det(ind(L_A))$ decides an orientation of
 $U$. $\Box$

 Furthermore, we have
the inclusion map
$$S: U \rightarrow \E, \leqno(2.25)$$
which can be viewed as a section of $E=\pi^*\E$. $S$ is proper since $\S$ is proper.
Moreover, 
$$S^{-1}(0)=\M_{\S}. \leqno(2.26)$$
Here, we construct a virtual neighborhood $(U, E, S)$ of $\M_{\S}$. To simplify the
notation, we will often use the same notation to denote the bundle (form) and its pull-back,

Notes that for any cohomology class $\alpha\in H^*(\B, \Z)$, we can pull back
$\alpha$ over $U$. Suppose that it is represented
by a closed differential form on $U$ (still denoted it by $\alpha$)
\vskip 0.1in
\noindent
{\bf Definition 2.8: }{\it Suppose that $det(ind (L_A))$ has a nowhere 
vanishing section so that $U$ is oriented.
\vskip 0.1in
\noindent
(1). If $deg(\alpha)\neq ind(L_A)$, we define  virtual neighborhood invariant 
$\mu_{\S}$ to be zero. 
\vskip 0.1in
\noindent
(2).When $deg(\alpha)=ind(L_A)$, choose a Thom form $\Theta$
supported in a neighborhood of zero section of $E$. We define 
$$\mu_{\S}(\alpha)=\int_U\alpha\wedge S^*\Theta.$$}
\vskip 0.1in
\noindent
{\bf Remark: }{\it In priori, $\mu_S$ is a real number. However, it was pointed to the author
by S. Cappell that when $\alpha$ is a rational cohomology class, $\mu_S(\alpha)$ is a rational
number. This is because both $U, E$ have  fundamental classes in compacted supported
rational homology. Then, $\mu_S(\alpha)$ can be interpreted as paring with the
fundamental class in rational cohomology.}
\vskip 0.1in
\noindent
{\bf Proposition 2.9: }{\it 
\noindent
(1). $\mu_{\S}$ is independent of $\Theta, \alpha$.
\vskip 0.1in
\noindent
(2). $\mu_{\S}$ is independent of the choice of $s$ and $\E$.}
\vskip 0.1in
{\bf Proof: }(1). If $\Theta'$ is another Thom-form  supported in a neighborhood of zero
section, there is a  $(k-1)$-form $\theta$ supported a neighborhood of zero section
such that
$$\Theta-\Theta'=d\theta.\leqno(2.28)$$
Then,
$$\int_U \alpha\wedge S^*\Theta-\int_U \alpha\wedge S^*\Theta'=\int_U\alpha
\wedge d(S^*\theta)=\int_U d(\alpha\wedge S^*\theta)=0. \leqno(2.28)$$
If $\alpha'$ is another closed form representing the same cohomology class,
it is the same proof to show
$$\int_U \alpha\wedge S^*\Theta=\int_U \alpha'\wedge S^*\Theta. \leqno(2.29)$$

To prove (2), suppose that $(\E', s')$ is another choice and $(U',E',
S')$ is the virtual neighborhood constructed by $(\E', s')$. Let $\Theta'$ be the 
Thom form of $E'$ supported in a neighborhood of zero section.
Consider 
$$\S^{(t)}_e=\S+(1-t)s+ts': \E\oplus \E'\times [0,1]\rightarrow
\F.\leqno(2.30)$$
Let $(U_{(t)},  \E\oplus \E', S_{(t)})$ be the virtual neighborhood
cobordism constructed by $\S^{(t)}_e$. By Stokes theorem,
$$\int _{U_0}\alpha\wedge S^*_0(\Theta\wedge \Theta')-\int_{U_1}\alpha \wedge
S^*_1(\Theta\wedge \Theta')=\int_{U_{(t)}}d (\alpha\wedge S^*_{(t)}(\Theta\wedge
\Theta'))=0,\leqno(2.31)$$
since both $\alpha$ and  $\Theta\wedge\Theta'$ are closed. 
It is easy to check that $U_0=\pi^* E'$ where $\pi: E\rightarrow U$ is the projection, 
$S_0=S\times Id$. Therefore,
$$\int _{U_0}\pi^*\alpha\wedge S^*_0(\Theta\wedge \Theta')=\int_U\alpha\wedge S^*(
\Theta)=\int_U\alpha\wedge S^*(\Theta).\leqno(2.32)$$
In the same way, one can show that
$$\int_{U_1}\alpha\wedge S^*_1(\Theta\wedge \Theta')=\int_{U'}\alpha\wedge
(S')^*(\Theta').$$
We have finished the proof. $\Box$
\vskip 0.1in
\noindent
{\bf Proposition 2.9: }{\it Suppose that $\S_t$ is a family of elliptic 
operators over $\F_t\rightarrow \B_t$ such that $\B_{(t)}=\bigcup_t\B_t\times
\{t\}$ is a smooth Banach V-cobordism and $\F_{(t)}=\bigcup_t \F_t 
\times \{t\}$ is a smooth V-bundle over $\B_{(t)}$. Furthermore, we assume that
$\M_{\S_{(t)}}=\bigcup_t\M_{\S_t}\times \{t\}$ is compact. We call $(\B_{(t)}, 
\F_{(t)}, \S_{(t)})$ a compact-V cobordism triple. Then
$\mu_{\S_0}=\mu_{\S_1}$. }
\vskip 0.1in
{\bf Proof: }  Choose 
$(\E_{(t)}, s)$ of $\F_{(t)}\rightarrow \U_{(t)}$ such that
$$\delta (\S^t+s) \leqno(2.33)$$
is surjective to $\F_{\U_{(t)}}$ where $\M_{\S_{(t)}}\subset \U_{(t)}
\subset \B_{(t)}$. Repeating previous argument, we construct a virtual 
neighborhood cobordism $(U_{(t)}, E_{(t)}, S_{(t)})$. Then, it is easy to check
that $(U_0, E_0, S_0)$ is a virtual neighborhood of $\S_0$ defined by
 $(\E_0, s(0))$  and
$(U_1,  E_1, S_1)$ is a virtual neighborhood of $\S_1$ defined by 
 $(\E_1, s(1))$. Applying the Stokes theorem as before, we have 
 $\mu_{\S_0}=\mu_{\S_1}$. $\Box$

Recall that by \cite{Sa2} one can define connections and curvatures on a
V-bundle. Then, characteristic classes can be defined by Chern-Weil formula in
 the category of V-bundle.
Next, we prove a proposition which  is very useful to calculate $\mu_{\S}$.
\vskip 0.1in
\noindent
{\bf Proposition 2.10: }{\it (1) If $F$ is a transverse section, 
$\mu_{\S}(\alpha)=\int_{\M_{\S}}\alpha$.
\vskip 0.1in
\noindent
(2) If $Coker L_A$ is constant and $\M_{\S}$ is a smooth V-manifold such that
$dim (\M_{\S})=ind(L_A)+\dim Coker L_A$, $Coker L_A$ forms an 
obstruction V-bundle $\H$ over $\M_{\S}$. In this case,
$$\mu_{\S}(\alpha)=\int_{\M_{\S}}e(\H)\wedge \alpha.\leqno(2.34)$$}
\vskip 0.1in
Before we prove the proposition, we need following lemma
\vskip 0.1in
\noindent
{\bf Lemma 2.11: }{\it Let $E\rightarrow M$ be a V-bundle over a V-manifold. Suppose
that $s$ is a transverse section of $E$. Then the Euler class $e(E)$ is dual
to $s^{-1}(0)$ in the following sense: for any compact supported form $\alpha$ 
with $deg(\alpha)=\dim M-\dim E$,
$$\int_M e(E)\wedge \alpha=\int_{s^{-1}(0)} \alpha.\leqno(2.35)$$}
\vskip 0.1in
{\bf Proof: } When $\dim \H=\dim \M_{\S}$, it is essentially Chern's proof
of Gauss-Bonnett theorem. By \cite{Sa2}, Chern's proof in smooth case  holds for
V-bundle. For general case, it is an easy generalization of Chern's proof using
normal bundle. We omit it. $\Box$

{\bf Proof of Proposition 2.10: } (1) follows from the definition where we take $k=0$.

To prove (2), let $F_{b}$ be the eigenspace of Laplacian $L_AL^*_A$ of an
eigenvalue $b$. Since $rank(Coker L_A)$ is constant, there is a $a\not\in 
Spec(L_A)$ for $A\in \M_{\S}$ such that the eigenspaces 
$$F_{\leq a}=\oplus_{b\leq a} F_b=Coker L_A \leqno(2.36)$$
has dimension $dim Coker(L_A)$ over $\M_{\S}$. Then, the same is true for an open
 neighborhood of
$\M_{\S}$. Without the loss of generality, we can assume that the open
neighborhood is $\U$. Therefore $F_{\leq a}$ form a V-bundle (still denoted
by $F_{\leq a}$) over $\U$ whose
restriction over $\M_{\S}$ is $\H$. In this case, we can choose $s$ such
that $s\in F_{\leq a}$ and $s$ satisfy Lemma 2.4.
Let $(U, E, S)$ be the virtual neighborhood constructed from $s$. 
 Recall that
$$U=(\S_e)^{-1}(0), \leqno(2.37)$$
where 
$$\S_e=\S+s.\leqno(2.38)$$
Let 
$$p_{\leq a}: \F\rightarrow F_{\leq a}\leqno(2.39)$$
be the projection. Then,
$$\S_e=p_{\leq a}(\S+s)+(1-p_{\leq a})(\S+s)=p_{\leq a}(\S+s)+(1-p_{\leq a})
(\S). \leqno(2.40)$$
The last equation follows from the fact that $s\in F_{\leq a}$.
So, $\S_e=0$ iff 
$$p_{\leq a}(\S+ s)=0 \mbox{ and } (1-p_{\leq a})(\S)=0.\leqno(2.41)$$
By our assumption, $(1-p_{\leq a})(\S)$ is transverse to the zero section over 
$\M_{\S}$ since $Coker(L^A)= F_{\leq a}$. Therefore, we can assume
that $(1-p_{\leq a})(\S)$ is transverse to the zero section over $\U$. Hence,
$((1-p_{\leq a})(\S))^{-1}(0)$ is a smooth V-manifold of dimension $ind(L_A)+\dim
F_{\leq a}=ind(L_A)+\dim Coker(L_A).$ But 
$$\M_{\S}\subset ((1-p_{\leq a})(\S))^{-1}(0)\leqno(2.42)$$
is a compact submanifold of the same dimension. Then, $\M_{\S}$ consists
of the components of $((1-p_{\leq a})(\S))^{-1}(0)$. In particular, other components
are disjoint from $\M_{\S}$. Therefore, we can choose smaller $\U$ to
exclude those components. Without the loss of generality, we can assume that
$$((1-p_{\leq a})(\S))^{-1}(0)=\M_{\S}.\leqno(2.43)$$
Since $\S=0$ over $\M_{\S}$, the first equation of (2.31)becomes 
$$p_{\leq a}(F+ s)= s=0.\leqno(2.44)$$
Therefore, $U\subset E_{\M_s}$ and
$$U=s^{-1}(0).\leqno(2.45)$$
However, $s$ is a transverse section by the construction. By Lemma 2.11,
$$\int_U\alpha \wedge S^*(\Theta)=\int_{E_{\M_{\S}}}\pi^*(e(\H)\wedge 
\alpha) \wedge \Theta=\int_{\M_{\S}}e(\H)\wedge \alpha, \leqno(2.46)$$
since $S: E_{\M_{\S}}\rightarrow E_{\M_{\S}}$ is identity. Then,
 we proved (2). $\Box$

\section{Virtual neighborhoods of Cauchy-Riemann equation}
   This is a technical section about the local structure of $\overline{\B}_A(
Y,g,k)$ and Cauchy-Riemann equation. Roughly speaking, we will show that for all
 the applications of this
article $\overline{\B}_A(Y,g,k)$ behaves like a Banach V-manifold. Namely,
$\overline{\B}_A(Y,g,k)$ is  VNA.  If 
readers only want to get a sense of big picture, one can skip over this section.

   There are roughly two steps in the virtual neighborhood construction. First step
is to define an extended equation $\S_e$ by the stabilization. Then, we need to prove that
(i) The set  $\U_{\S_e}=\{x, Coker D_x \S_v=\emptyset\}$ is open; (ii) $\U_{\S_e}\cap
\S^{-1}_e(0)$ is a smooth, oriented V-manifold. Ideally, we would like to set up
some Banach manifold structure on our configuration space and treat $\U_{\S_e}
\cap\S^{-1}_e(0)$ as a smooth submanifold. However, there are some basic
analytic difficulty against such an approach, which we will explain now.
For $\B_A(Y,g,k)$, we allow the domain of the map to vary to accomendate the variation
of complex structures of Riemann surfaces. Let's look at a simpler model.
 Suppose that $\pi: M\rightarrow N$ be a fiber bundle with
fiber $F$. We want to
put a Banach manifold structure on $\bigcup_{x\in N}C^k(\pi^{-1}(x))$. A natural way is 
to choose a local trivialization $\pi^{-1}(U)\cong U\times F$. It induces a 
trivialization $\bigcup_{x\in U}C^k(\pi^{-1}(x))\rightarrow U\times C^k(F)$. Then,
we can use the natural Banach manifold structure on $C^k(F)$ to induce
a Banach manifold structure on $\bigcup_{x\in U}C^k(\pi^{-1}(x))$. However, if we have
a different local trivialization, the transition function is a map $g:U\rightarrow
Diff(F)$. The problem is that $Diff(F)$ only acts on $C^k(F)$ continuously. For 
example, suppose that $\phi_t$ is a one-parameter family of diffeomorphisms generated by a
vector field $v$. Then, the derivative of the path $f\circ g_t$ is $v(f)$, which
decreases the differentiability of $f$ by one. So we do not have a natural Banach
manifold structure on  $\bigcup_{x\in N}C^k(\pi^{-1}(x))$ in general. It is obvious
that we have a natural Frechet manifold structure on $\bigcup_{x\in N}C^{\infty}(
\pi^{-1}(x))$. However, we only care about the zero set $\M$ of some elliptic operator
$\S_e$ defined over Frechet manifold  $\bigcup_{x\in N}C^{\infty}(\pi^{-1}(x))$. The
crucial observation is that locally we can choose any local trivialization and use
 Banach manifold structure induced from the local trivialization to show that
$\M_U=\M\cap U\times C^k(F)$ is smooth. The elliptic regularity implies that $\M_U \subset U\times 
C^{\infty}(F)$.  Although the transition map is not smooth for
$C^k(F)$, but it is smooth on $\M_U$. Therefore, $\M_U$ patches
together to form a smooth manifold. Our strategy is to define the extended
equation $\S_e$ over the space of $C^{\infty}$-stable map. In each coordinate chart,
we enlarge our space with Sobolev maps. Then, we can use usual analysis to show that
the moduli space can be given a local coordinate chart of a smooth manifold. Elliptic 
regularity guarantees that every element
of the moduli space is indeed smooth. Then, we show that the moduli space in each coordinate
chart patches up to form a $C^1$-V-manifold.

 Suppose that $(Y, \omega)$ is a family of
symplectic manifold and $J$ is a tamed almost complex structure. Choose a metric
tamed with $J$.
\vskip 0.1in
\noindent
{ \bf Definition 3.1 ([PW], [Ye], [KM]). }{\it Let $(\Sigma, \{x_i\})$ be a stable Riemann surface.
A stable holomorphic map (associated with $(\Sigma, \{x_i\})$) is an equivalence class 
of continuous maps $f$ from $\Sigma'$ to
$Y$ such that  $f$  has the 
image in a fiber of $Y\rightarrow X$ and is  smooth at smooth points of $\Sigma'$, 
where the domain $\Sigma'$ is obtained by joining chains of $\P^1$'s
at some double points of $\Sigma$ to separate the two components, and then 
attaching some
trees of $\P^1$'s. We call components of $\Sigma$ {\em principal components} 
and others {\em bubble components}. Furthermore,
\begin{description}
\item[(1)] If we attach a tree of $\P^1$ at a marked point $x_i$, then $x_i$ will
be replaced by a point different from intersection points on some component of the tree. 
Otherwise, the marked points do not change.
\item[(2)] The singularities of $\Sigma'$ are normal crossing and there are at 
most two components intersecting at one point.
\item[(3)]  If the restriction of $f$ on a bubble component is constant, then it has
at least three special points (intersection points or marked points). We call
this component  {\em a ghost bubble} \cite{PW}.
\item[(4)]  The restriction of $f$ to each component is  $J$-holomorphic. 
\end{description}
Two such maps are equivalent if one is the composition of the other with
an automorphism of the domain of $f$.

   If we drop the condition (4), we simply call $f$ a stable map. Let
$\overline{\M}_{A}(Y,g,k,J)$ be the moduli space of stable holomorphic maps and
$\overline{\B}_A(Y,g,k)$ be the space of stable maps.}
\vskip 0.1in
\noindent
{\bf Remark 3.2: }{\it There are two types of automorphism  here. Let
$Aut_f$ be the group of automorphisms $\phi$ of the domain of $f$ such that
$f\circ \phi$ is also holomorphic. This is the group we need to module out when
we define $\overline{\M}_{A}(Y,g,k,J)$ and $\overline{\B}_A(Y,g,k)$. It
consists two kinds of elements. (1) When some
bubble component is not stable with only one or two marked points,
there is a continuous subgroup of $PSL_2\C$ preserving the marked points. 
(2) Another type of element comes from the automorphisms of domain interchanging
 different components, which form a finite group. Let $stb_f$ be the 
subgroup of $Aut_f$ preserving $f$. It is easy to see that $stb_f$ is always a 
finite group.  Type (1) elements of $stb_f$ appear with multiple
covered maps.} 
\vskip 0.1in
\noindent
{\bf Proposition 3.3: }{\it $\overline{\B}_A(Y,g,k)$ (whose topology is defined 
later) is a stratified Hausdorff Frechet V-manifold of finite many strata.}
\vskip 0.1in
The proof consists of several lemmas.
\vskip 0.1in
\noindent
{\bf Lemma 3.4: }{\it $\B_A(Y,g,k)$ is a Hausdorff Frechet V-manifold for any
$2g+k\geq 3$ or $g=0, k\leq 2, A\neq 0$.}
\vskip 0.1in
{\bf Proof: } Recall
$$\B_A(Y,g,k)=\{(f, \Sigma); \Sigma\in \M_{g,k}, f: \Sigma \stackrel{F}{
\rightarrow} Y \},\leqno(3.1)$$
where $\stackrel{F}{\rightarrow}$ means that the image  is in a fiber.
When $2g+k\geq 3$, $\Sigma$ is stable and $\M_{g,k}$ is a V-manifold. 
Hence, the automorphism group $Aut_{\Sigma}$ is finite. Furthermore, there is a
$Aut_{\Sigma}$-equivariant holomorphic fiber bundle
$$\pi_{\Sigma}: U_{\Sigma}\rightarrow O_{\Sigma}$$
such that $O_{\Sigma}/Aut_{\Sigma}$ is a neighborhood of $\Sigma$ in $\M_{g,k}$ and fiber
$\pi^{-1}_{\Sigma}(b)=b$. Consider
$$\U_{\Sigma, f}=\{(b, h); h: b\stackrel{F}{\rightarrow} Y, h\in C^{\infty}.\}
\leqno(3.2)$$
As we discussed in the beginning of this section, $\U_{\Sigma, f}$ has a natural
Frechet manifold structure. Let $stb_f\subset Aut_{\Sigma}$ be the subgroup preserving $f$. 
One can observe that  $\U_{\Sigma, f}/stb_f$ is a neighborhood of $(\Sigma, f)$
in $\B_A(Y,g,k)$. Hence, $\B_A(Y,g,k)$ is a  Frechet V-manifold. Since
only a finite group is involved, $\B_A(Y,g,k)$ is obviously Hausdorff.

For the case $g=0, k\leq 2, A\neq 0$, $\Sigma$ is no longer stable and the
automorphism group $Aut_{\Sigma}$ is infinite. Here, we fix our marked points
at $0$ or $0,1$. First of all, $stb_f$ is finite for any $f\in Map^F_A(Y,0,k)$
with $A\neq 0$.
$$\B_A(Y,g, k)=Map^F_A(Y,0,k)/Aut_{\Sigma}.$$
We first show that $B_A(Y,
g,k)$ is Hausdorff. It requires showing that the graph
$$\Delta=\{(f, f\tau); f\in Map^F_A(Y,0,k), \kappa\in Aut_{\Sigma}\}\leqno(3.3)$$
is closed. Suppose that $(f_n, f_n\tau_n)$ converges to $(f,h)$ uniformly for all its
derivatives. We claim that
$\{\tau_n\}$ has a convergent subsequence. Suppose that $\infty$ is one of marked point
which $\tau_n$ fixes. They, $\tau_n$ can be written as $a_nz+b_n$ for $a_n\neq 0$.

   Suppose that $\tau_n$ is degenerated. Then, (i)$ b_n \rightarrow \infty$, (ii)
$a_n\rightarrow 0$ or (iii) $a_n\rightarrow \infty$. In each case, we observe that
$\tau_n$ converges pointwisely to $\tau$ which is either a constant map taking value 
$\infty$ or a map
taking two different values. Since $f_n$ converges uniformly, $f_n\tau_n$ converges to
$f\tau$ pointwisely. Hence, $h=f\tau$ which is either a constant map or discontinous.
We obtain a contradiction. Suppose that $\tau_n$ converges to $\tau$. Then, $f_n\tau_n$
converges to $f\tau$. Therefore, $\Delta$ is closed. 

Notes that 
$$||df||_{L^2}\geq \omega(A).\leqno(3.4)$$
Choose the standard metric on $\P^1$ with volume 1. Then, for a holomorphic map,
there are points $p$ (hence
a open set of them) 
such that $df(p)$ is of maximal rank and $|df(p)|\geq \frac{1}{2}\omega(A)$. 
Since we only want to construct a neighborhood and the condition above is an open
condition. Without the loss of generality, we assume that it is true for any $f$.

We marked 
extra points $e_i$ such
that $df(e_i)$ is of maximal rank, $|df(e_i)|\geq \frac{1}{2}\omega(A)$ and 
$(\Sigma, e_i)$ has three marked points. 

Next we   construct a slice $W_f$ of the 
action $Aut_{\Sigma}$. Note that $Map^F_A(Y,0,k)$ is only a Frechet manifold. We can not
use implicit function theorem. Since $stb_f$!is finite, we can construct a $stb_f$ invariant
metric on $f^*TY$ by averaging the existing metric. Using $stb_f$ invariant metric, the set 
$$\{w\in \Omega^0(f^*T_FY); ||w||_{L^p_1}< \epsilon\}\leqno(3.5)$$
is $stb_f$-invariant and open in $C^{\infty}$-topology. Now, we fix the $stb_f$-invariant metric. For each extra 
marked point $e_i$ constructed in previous paragraph, 
$df(e_i)$ is a $2$-dimensional vector space. Clearly, 
$$f_{e_i}=\oplus_{\tau\in stb_f} df(\tau(e_i))\subset (T_{f(e_i)}
Y)^{|stb_f|}$$
 is
$stb_f$-invariant. Now we want to construct a 2-dimensional subspace $E_{e_i}\subset f_{e_i}$ which is the
orbit of action $Aut_{\Sigma}$. For simplicity, we assume that we only need one extra marked point
$e_1$
 to stabilize $\Sigma$. The proof of the case with two extra marked points is the same.

 In this case, a neighborhood of $id$ in $Aut_{\Sigma}$ can be identified with a neighborhood of
$e_1$ by the relation $\tau_x(e_1)=x$ for $x\in D^2(e_1)$.  $\frac{d}{dx}\tau_x(f)(y)|_{x=e_1}
=df(y)(v(y))$, 
where $v=\frac{d}{dx}\tau_x|_{x=e_1}$ is a holomorphic vector field. By our identification, 
$v$ is
decided by its value $v(e_1)\in T_{e_1}S^2$. Given any $v\in T_{e_1}S^2$, we use $v_{e_1}\in T_{id}
Aut_{\Sigma}$ to denote its extension. Therefore, $v$ decides $v_{e_1}(\tau(e_1))$.
To get a precise relation, we can differentiate $\tau_x(\tau(e_1))=\tau (\tau^{-1} \tau_x \tau)(
e_1)$ to
obtain 
$$v_{e_1}(\tau(e_1))=D\tau Ad_\tau(v),\leqno(3.6)$$
 where $Ad_\tau$ is the adjoint action. 
$$E_{e_i}=\{ \oplus_{\tau\in stb_f}  df(\tau(e_1))(v_{e_1}(\tau(e_1))); v\in T_{e_1}S^2\}\leqno(3.7)$$
It is easy to check that $E_{e_1}$ is indeed $stb_f$-invariant. We can identify $E_{e_1}$ with
$T_{e_1}S^2$ by 
$$v\rightarrow \oplus_{\tau\in stb_f}  df(\tau(e_1))(D\tau Ad_\tau(v)),\leqno(3.8)$$
Hence,  $E_{e_i}$ is 2-dimensional. Given any $w\in \Omega^0(f^*T_FY)$, we say that
 {\em $w\perp E_{e_i}$ if $\oplus_{\tau\in stb_f} w(\tau(e_i))$
is orthogonal to  $E_{e_i}$}. The slice $W_f$ can be 
constructed as
$$W_f=exp_f \{  w\in \Omega(f^*T_FY); ||w||_{L^p_1}< \epsilon, ||w||_{C^1(D_{\delta_0}(g(e_i)))}<
\epsilon \mbox{ for } g\in stb_f, w\perp E_{e_i}\},\leqno(3.9)$$
where $T_FY$ is the direct sum of vertical tangent bundle and $P^*TX$ and $\delta_0$ is
a small fixed constant. We need to show that
\begin{description}
\item[(1)] $W_f$ is invariant under $stb_f$.
\item[(2)] If $h\tau\in W_f$ for $h\in W_f$, then $\tau\in Stb_f$.
\item[(3)] There is a neighborhood $U$ of $id\in Aut$ such that the 
multiplication $F: U\times W_f\rightarrow Map^F_A(Y,0,k)$ is a homeomorphism
onto a neighborhood of $f$.
\end{description}

  (1) follows from the definition. For (2), we claim that the set of $\tau$ satisfying $(2)$ is
close to an element of $stb_f$ for small $\epsilon$. If not, there is a 
neighborhood $U_0$ of $stb_f$  and a sequence of $(h_n ,\tau_n)$ such that $\tau_n\not
\in U_0$,  $h_n$
converges to $f$ and $h_n\tau_n$ converges to $f$. By the previous argument, $\tau_n$
has a convergent subsequence.  Without the loss of generality, we can assume
that $\tau_n$ converges to $\tau\not\in U_0$. Then, $h_n\tau_n$ converges to $f\tau=f$.
This is a contradiction. By (1), we can assume that $\tau$ is close to identity.
Then, (2) follows from (3). 

  Next we prove (3). 
Consider the local model around $f(\tau(e_1))$. Since $df(\tau(e_1))$
is injective, we can choose a local coordinate system of $V$ such that
$Im(f)$ is a ball of $\C_{\tau}\subset \C_{\tau}\times \C^{n-1}_{\tau}$ in which
the origin corresponds to $f(\tau(e_1))$.. 
Furthermore, we may assume that the metric
is standard. For any $w$, let 
$$P(w): \Omega^0(f^*T_FY)\rightarrow E_{e_1}.$$
be the projection
Then, $w\in W_f$ iff $P(w)=0$.
Suppose that $w$ is bounded. 
$$\tau_x(w)(\tau(e_1))=w(\tau_x(\tau(e_1)))+f(\tau_x(\tau(e_1)))-f(\tau(e_1))+O(r^2)=
w(\tau_x(\tau(e_1)))+f(\tau_x(\tau(e_1)))+O(r^2),\leqno(3.10)$$
where $r=|\tau_x(\tau(e_1))|$.
Then, 
$$P(\tau_x(w))=P(w\circ \tau_x)+P(f\circ \tau_x)+O(r^2).\leqno(3.10.1)$$
Hence $P(\tau_x(w))=0$ iff $-P(w\circ \tau_x)=P(f\circ \tau_x)+O(r^2)$, where
$$P(w\circ \tau_x), P(f\circ \tau_x): D^2\rightarrow E_{e_1}.\leqno(3.10.2)$$
Notes that $P(f\circ \tau_0)=0$. 
$$dP(f\circ \tau_x)(v)|_{x=0}=P(df(v_{e_1})).\leqno(3.10.3)$$ 
 Under the identification (3.8), $dP(f\circ \tau_x)_0$ is the identity. 
Let $\bar{f}=P(f\circ \tau_x)$. Then, $\bar{f}^{-1}$ exists and $d\bar{f}^{-1}$ is bounded on
a small disc. Consider $\bar{w}(x)=\bar{f}^{-1}P(w\circ \tau_x+O(r^2))$. Then, $P(\tau_x(w))=0$ iff
$x$ is a fixed point of $\bar{w}$. 
 Suppose that $\epsilon <<1$. Since $||w||_{C^1(D_{\delta_0}(g(e_i)))}<\epsilon$, 
 $|\bar{w}(0)|<C\epsilon$.
Furthermore, $|dw|<\epsilon$. $\bar{w}: D^2_{\delta_0}
\rightarrow D^2_{\delta_0}$ for fixed $\delta_0$. The small bound on the derivative also implies that
$\bar{w}$ is a contraction mapping. Therefore, there is a unique fixed point $x(w)$  in 
$D_{\delta_0}$ and hence $\tau_w=\tau_{x}$. Moreover, $x(w)$ depends 
smoothly on $w$. Therefore, $\tau_w$ depends smoothly on $w$.
We define $H(w)=(\tau^{-1}_w, f_w\tau_w)$. By our construction, $H$ is continuous and
an inverse of $F$. $\Box$

$\overline{\M}_A(Y,g,k,J)$ has an obvious stratification indexed
by the combinatorial type of the domain. The later can be viewed as the topological 
type of the domain as an abstract 2-manifold with marked points such that each 
component
is associated with a nonzero  integral 2-dimensional class $A_i$ unless this 
component is genus zero with at least three marked points. Furthermore, each
component is represented by a $J$-holomorphic map with fundamental class $A_i$
and total energy is equal to
$\omega(A)$. Suppose that $\D^{J,A}_{g,k}$ is the set of indices.
\vskip 0.1in
\noindent
{\bf Lemma 3.5: }{\it $\D^{J,A}_{g,k}$ is a finite set.}
\vskip 0.1in
{\bf Proof:}  Let
 $(A_1, \cdots, A_k)$ be the integral 2-dimensional
nonzero classes associated with the components. The last condition implies that 
$$\omega(A_i)>0, \sum A_i=A. \leqno(3.11)$$
\vskip 0.1in
In \cite{RT1}(Lemma 4.5), it was shown that the set of tuple (3.11) is 
finite. Therefore, the number of non-ghost components is bounded. We claim
that the number of ghost bubbles  is bounded by the number of non-ghost
bubbles. Then, the finiteness of $\D^{J,A}_{g,k}$ follows automatically.

We prove our claim by the induction on the number of non-ghost bubbles. 
It is easy to observe that any ghost bubble must lie in some bubble tree $T$.
 By the construction, this ghost bubble can not lie on the
tip of any branch. Otherwise, it has at most two marked points.
Choose $B$ to be the ghost bubble closest to the tip. We remove the subtree $T_B$ with base $B$. Then, we obtain an abstract 2-manifold with
marked points. If it is  the domain of another stable map, we denote it by
$T'$. If not,  $B$ is based on another ghost bubble $B'$ with only three marked
points. Then, we remove $T_B$ and contract $B'$ to obtain $T'$ the domain of
another stable map.  Let $gh(T')$ be the number of ghost 
bubbles and $ngh(T')$ be the number of non-ghost bubbles. By the induction,
$$gh(T')\leq ngh(T').\leqno(3.12)$$
However,
$$gh(T)\leq gh(T')+2, ngh(T_B)\geq 2.$$
Therefore, 
$$gh(T)\leq ngh(T')+2\leq ngh(T)+ngh(T_B)=ngh(T).\leqno(3.13)$$
We finish the proof. $\Box$
\vskip 0.1in
For any $D\in \D^{J,A}_{g,k}$, let $\B_D(Y,g,k)
\subset \overline{\B}_A(Y,g,k)$ be the set of stable maps whose domain
and the corresponding fundamental class of each component have type $D$. Then,
$\B_D(Y,g,k)$ is a strata of $\B_A(Y,g,k)$.
\vskip 0.1in
\noindent
{\bf Lemma 3.6: }{\it $\B_D(Y,g,k)$ is a  Hausdorff Frechet  
V-manifold.}
\vskip 0.1in
{\bf Proof: } 
  $\B_D(Y,g,k)$ is a subset of $\prod_i \B_{A_i}(
Y, \Sigma_i)$ such that the components intersect each other according to
the intersection pattern specified by $D$. Therefore, it is Hausdorff.
For the simplicity, let's
consider the case that $D$ has only two components. The general case is the
same.

 Let $D=\Sigma_1\wedge \Sigma_2$ joining at $p\in \Sigma_1, q\in \Sigma_2$.
Assume that $A_i$ is associated with $\Sigma_i$.
Then,
$$\B_D(Y,g,k)=\{(f_1, f_2)\in \B_{A_1}(Y, g_1, k_1+1)\times 
\B_{A_2}(Y, g_2, k_2+1); f_1(p)=f_2(p)\}.\leqno(3.14)$$
It is straightforward to show that $\B_D(Y,g,k)$ is Frechet V-manifold with
the tangent space
$$T_{(f_1, f_2)}\B_D(Y,g,k)=\{(w_1, w_2)\in \Omega^0(f^*_1T_FV)\times \Omega^0
(f^*_2T_FV); w_1(p)=w_2(q)\} \leqno(3.15)$$
We leave it to readers. $\Box$

Next, we discuss how  different strata fit together. It amounts to show
how a stable map deforms when it changes domain. A natural starting point is the
deformation theory of the domain of stable maps as abstract nodal Riemann
surfaces. However, it is well-known that unstable components cause a problem in
the deformation theory. For example, the moduli space will not be Hausdorff.
To have a good deformation theory, we have to consider a map with its domain together
for unstable components. 

 Let $\overline{\M}_{g,k}$ be the space of  stable Riemann surfaces.
The important properties of $\overline{\M}_{g,k}$ are that (i)
$\overline{\M}_{g,k}$ is a V-manifold; (ii) there is a local universal 
V-family in following sense: for each $\Sigma\in \overline{\M}_{g,k}$, let 
$stb_{\Sigma}$ be its automorphism group.
 There is a $stb_{\Sigma}$-equivariant (holomorphic) fibration 
$$\pi_{\Sigma}: U_{\Sigma} \rightarrow O_{\Sigma}\leqno(3.16)$$
such that $O_{\Sigma}/Aut_{\Sigma}$ is a neighborhood of $\Sigma$ in 
$\overline{\M}_{g,k}$ and the fiber $\pi^{-1}_{\Sigma}(b)=b$.

 Suppose that the components of $f$ are 
$(\Sigma_1, f_1), \cdots, (\Sigma_m, f_m)$, where $\Sigma_i\in \M_{g_i, k_i}$ is a
marked
Riemann surface. If $\Sigma_i$ is stable, locally $\M_{g_i. k_i}$ is a
V-manifold  and have a local universal V-family. Suppose that they are
$$\pi: U_i\rightarrow O_i\leqno(3.17)$$
divided by the automorphism group $Aut_i$ of $\Sigma_i$ preserving the marked
points. Stability means that $Aut_{\Sigma_i}$ is finite. However, the relevant
group for our purpose is $stb_i=stb_{f_i}\subset Aut_i$. Suppose that $x_{i1}, 
\cdots, x_{ik_i}$ are the
marked points. We choose a disc $D_{ij}$ around each marked point $x_{ij}$
invariant under $stab_{\Sigma_i}$. 
For each $\tilde{\Sigma}_i\in O_i$, $x_{ij}$ may vary. We can find a 
diffeomorphism $\phi_{\Sigma}: \Sigma \rightarrow \tilde{\Sigma}_i$ to carry $x_{ij}$ together
with $D_{ij}$  to the corresponding marked point and its neighborhood on 
$\tilde{\Sigma}_i$. Pulling back the complex
structures by $\phi_{\tilde{\Sigma}_i}$, we can view
$O_i$ as the set complex structure on $\Sigma_i$ which have the same marked
points and moreover are the same on $D_{ij}$. $\phi_{\tilde{\Sigma}_i}$
gives a local smooth trivialization
$$\phi_{\Sigma}: U_i\rightarrow O_i\times \Sigma.\leqno(3.18)$$
When $\Sigma_i$ is unstable, $\Sigma_i$ is a sphere with one or two marked 
points and we have to divide it by the subgroup $Aut_i$ of $\P^1$ preserving the marked
points. But to glue the Riemann surfaces, we have to choose a parameterization.
Recall that $\B_{A_i}(\Sigma_i)=Map^F_{A_i}(\Sigma_i, Y)/Aut_{i}$.
For any $f_i\in Map^F_{A_i}(\Sigma_i, Y)$, one  constructs a slice $W_{f_i}$
(Lemma 3.4) at $f_i$ such that $W_{f_i}/stb_{f_i}$ is diffeomorphic to a neighborhood of 
$[f_i]$ in the quotient. Moreover, we only want to construct a neighborhood of $f$.
To abuse notation, we identify $\B_{A_i}(\Sigma_i)$ with the slice
$W_{f_i}/stb_{f_i}$. Then, we can proceed as before. Fix a standard $\P^1$.
We choose a disc $D_{ij}$ ($j\leq 2$) around each marked point invariant under
$stb_{f_i}$. Then, $O_i=pt, U_i=\P^1$. 

   Let $\N$ be the set of the
nodal points of $\Sigma$. For each $x\in \N$, we associate a copy of $\C$ (gluing
parameter) and denote it by $\C_x$. Let $\C_{f}=\prod_{x\in \N} \C_x$, which
is a finite dimensional space. For each $v\in C_{f}$ with $|v|$ small 
and $\tilde{\Sigma}_i\in O_i$ 
, we can 
construct  a Riemann surface $\tilde{\Sigma}_v.$ 
Suppose that $x$ is the intersection point of $\Sigma_i, \Sigma_j$ and 
$\Sigma_i, \Sigma_j$ intersect at 
$p\in 
\Sigma_i, q\in \Sigma_j$.  For any small complex number
$v_x=re^{iu}$. We construct $\Sigma_i \#_{v_x} \Sigma_j$ by cutting discs
with radius $\frac{2r^2}{\rho}$-$D_p(\frac{2r^2}{\rho}), D_q(\frac{2r^2}{\rho})$,
where $\rho$ is a small constant to be fixed later. Then, we identify two annulus $N_p(\frac{\rho r^2}{2}
,\frac{2r^2}{\rho}), N_q(\frac{\rho r^2}{2},\frac{2r^2}{\rho})$ by holomorphic map
$$(e^{i\theta}, t)\cong (e^{i\theta}e^{iu}, \frac{r^4}{t}).\leqno(3.19)$$
Notes that (3.19) sends inner circle to outer circle and vis versus. Moreover, we
identify the circle of radius $r^2$. Roughly speaking,
we cut off the discs of radius $r^2$ and glue them together by rotating $e^{i\theta}$.
When $v_x=0$, we define $\Sigma_i\#_0 \Sigma_j=\Sigma_i\wedge \Sigma_j$-the
one point union at $p=q$. Given any metric $\lambda=(\lambda_1, \lambda_2)$ on $\Sigma$, we can
patch it up on the gluing region as follows. Choose coordinate system of $N_p(\frac{\rho r^2}
{2}, \frac{2r^2}{\rho})$. The metric of $\Sigma_1$ is $t(ds^2+dt^2)$ and the metric from $\Sigma_2$
is $\frac{r^4}{t}(ds^2+dt^2)$. Suppose that $\beta$ is a cut off function vanishing
for $t<\frac{\rho r^2}{2}$ and equal to one for $t>\frac{2r^2}{\rho}$. We define a metric $\lambda_v$ which is
equal to $\lambda$ outside the gluing region and 
$$\lambda_v=(\beta t+(1-\beta)\frac{r^4}{t})(ds^2+dr^2)\leqno(3.20)$$
over the gluing region. 
We observe that on the annulus $N_p(\frac{\rho r^2}{2}, \frac{2r^2}{\rho})$ the metric $g_v$ has the same
order as standard metric.  For any complex structure on $\Sigma_i$ which is 
fixed on the gluing region, it induces a complex structure on $\Sigma_i\#_{v_x}
\Sigma_j$. If we start from the complex structure of $\tilde{\Sigma}$, by
repeating  above process for each nodal point we construct a marked  Riemann
surface  $\tilde{\Sigma}_v$. Clearly, $\tilde{\Sigma}_0=\tilde{\Sigma}$.
\vskip 0.1in
\noindent
{\bf Remark 3.7: }{\it The reader may wonder why we glue in a disc of radius $r^2$ 
instead of $r$. The reason is a technical one. If we use $r$, the gluing map is only 
continuous at $r=0$. Using $r^2$, we can show that the gluing map is $C^1$ at $r=0$.}
\vskip 0.1in

  Let 
$$\tilde{O}_{f}=\prod_i O_i\times \C_{f}.\leqno(3.21)$$
The previous construction yields a universal family
$$\tilde{U}_{f}=\cup \{ \tilde{\Sigma}_v; \tilde{\Sigma}\in \prod_i O_{
f}, v\in \C_{f} \mbox{ small }\}.\leqno(3.22)$$
The projection 
$$\pi_{f}: \tilde{U}_{f}\rightarrow  \tilde{O}_{f}\leqno(3.23)$$ 
maps $\tilde{\Sigma}_v$ to $(\tilde{\Sigma},v)$. We still need to show that (3.23)
is  $stb_{f}$-equivariant. $\prod_i stb_i$ induces an obvious action on (3.23). 
There are other types of automorphisms of 
$\Sigma$ by switching the different components and $stb_{f}$ is a finite
extension of $\prod_i stb_i$ by such automorphisms. The gluing construction with perhaps
different gluing parameter is
clearly commutative with such automorphisms. Hence,  $stb_{f}$ acts
on (3.23). $(\tilde{U}_{f}, \tilde{O}_{f})/stb_{f}$ is the local 
deformation of domain we need. After we stabilize the unstable component,
$\tilde{\Sigma}_v$ should be viewed as an element of $\overline{\M}_{g,k+l}$, where
$l$ is the number of extra marked points. Hence, $\tilde{O}_{f}\subset
\overline{\M}_{g,k+l}$ and $\tilde{U}_{f}$ is just the local universal
family of $\overline{\M}_{g,k+l}$. Forgetting the extra marked points, we map
$\tilde{O}_{f}$ to $\overline{\M}_{g,k}$ by the map
$$\pi_{k+l}: \overline{\M}_{g,k+l}\rightarrow \overline{\M}_{g,k}\leqno(3.24)$$
 Suppose that the extra marked points are $e^v_1,\cdots,
e^v_l$. Sometimes, we also use notation $e^f_1, \cdots e^f_l$.

  To describe a neighborhood of $f$, without the loss of generality, we can
assume that  $dom(f)=\Sigma_1
\wedge \Sigma_2$ and $f=(f_1, f_2)$, where $\Sigma_1, \Sigma_2$ are marked
Riemann surfaces of genus $g_i$ and $k_i+1$ many marked points such that
$g=g_1+g_2, k=k_1+k_2$. Furthermore, suppose that $\Sigma_1, \Sigma_2$ intersects at 
the last marked points  $p, q$ of $\Sigma_1, \Sigma_2$ respectively.
 The general case is identical and we just repeat our
construction for each nodal point. In this case, the gluing parameter $v$ is a complex number. We
choose $v$ small enough such that marked points  other than $p,q$ are away from
the gluing region described above. Let $f_1(p)=f_2(q)=y_0\in V
\subset Y$. Let $U_{P(y_0)}$ be a small neighborhood of $P(y_0)\in X$. We can
assume that $P^{-1}(U_{y_0})=V\times U_{P(y_0)}$ and $y_0=(x_0, x_1)$. 
Suppose that the fiber exponential map $exp: 
T_{x_0} V\rightarrow V \times \{x\}$  is a diffeomorphism from $B_{\epsilon}(
x_0, T_{x_0}V)$ onto its image for
any $x\in U_{P(y_0)}$, where $B_{\epsilon}$ is a ball of radius $\epsilon$. Furthermore, we define
$$f^w=exp_f w.\leqno(3.25)$$

Next, we  construct attaching maps which define the topology of $\overline{\B}_A
(Y,g,k)$.
 First we construct a neighborhood  $\U_{f,D}/stb_{f}$ of $f\in \B_D(Y,g,k)$. 
Recall that
if $dom(f)=\Sigma$ is an irreducible stable marked Riemann surface, then a neighborhood
of $f$ can be described as 
$$O_{f}\times \{f^w; w\in \Omega^0(f^*T_FY), ||w||_{L^p_1}< \epsilon\}
\leqno(3.26)$$
divided by $stb_{f}$.

  If $\Sigma$ is unstable, we needs to find a slice $W_{f}$. By lemma 3.4, we mark 
additional points $e^f_i$ on 
$\Sigma$ such that $\Sigma$ has three marked points. We call the resulting
Riemann surface $\bar{\Sigma}$. Furthermore, we choose $e^f_i$ such that 
$df_{e^f_i}$ is of maximal rank. 
Then,
$$W_f=\{f^w; w\in \Omega^0(f^*T_FY); ||w||_{L^p_1}<\epsilon, ||w||_{C^{1}(D_{\delta_0}(g(e_i)))}< \epsilon, g\in
sbt_f,  w\perp E_{e^f_i}\}.\leqno(3.28)$$
If $dom(f)=\Sigma_1\wedge \Sigma_2$ joining at $p\in\Sigma_1, q\in \Sigma_2$ and 
$f=f_1\wedge f_2$,
We define 
$$\Omega^0(f^*T_FY)=\{(w_1, w_2)\in \Omega^0(f^*_1T_FY)\times 
\Omega^0(f^*_2T_FY); w_1(p)=w_2(q), w\perp E_{e^f_i}\}.\leqno(3.29)$$
A neighborhood of $f$ in $\B_D(Y,g,k)$ is 
$$\prod_i O_i\times \{f^w; w\in \Omega^0(f^*T_FY), ||w||_{L^p_1}<\epsilon, ||w||_{C^{1}(D_{\delta_0}(g(e_i)))}< \epsilon, g\in
sbt_f, w\perp E_{e^f_i}\}/stb_f.
\leqno(3.30)$$
If $dom(f)$ is an arbitrary configuration, we repeat above construction over
each nodal point to define $\Omega^0(f^*T_FY)$. A neighborhood of $f$ in $\B_D(Y,
g,k)$ is 
$$\U_{f,D}=\prod_i O_i\times \{f^w; w\in \Omega^0(f^*T_FV), ||w||_{L^p_1}<\epsilon, ||w||_{C^{1}(D_{\delta_0}(g(e_i)))}< \epsilon, g\in
sbt_f, w\perp E_{e^f_i}\}/stb_f.\leqno(3.31)$$

 We want 
to construct an attaching map
$$\bar{f}^{w,v}: \U_{f,D}\times \C^{\epsilon}_f\rightarrow \overline{\B}_A(Y,g,k)$$
invariant under $stb_f$, 
where $\C^{\epsilon}_f$ is a small $\epsilon$-ball around the origin of $\C_f$. We simply denote
$$\bar{f}^v=\bar{f}^{0,v}.\leqno(3.32)$$
Again, let's focus on the case
that $D=\Sigma_1\wedge \Sigma_2$ and the general case is similar. Recall the
previous set-up.   $f_1(p)=f_2(q)=y_0=(x_0, x_1)\in V
\subset Y$. Let $U_{P(y_0)}$ be a small neighborhood of $P(y_0)\in X$. We can
assume that $P^{-1}(U_{y_0})=V\times U_{P(y_0)}$ and $y_0=(x_0, x_1)$. 
Suppose that the fiber exponential map $exp: 
T_{x_0,x} V \rightarrow V \times \{x\}$  is a diffeomorphism from 
$B_{\epsilon}(x_0, T_{x_0}V)$ to its image for any $x\in U_{P(y_0)}$. 
 In the 
construction of $dom(f)_v$, we can choose $r$ small enough  such that 
$$f^w_1(D_p(\frac{2r^2}{\rho})), f^w_2(D_q(\frac{2r^2}{\rho}))
\subset  B_{\epsilon}(x_0, T_{x_0} V)\times P(y'_1),$$
for any $w\in \Omega^0(f^*T_FY)$ and $||w||_{C^{1}}<\epsilon$.
Following \cite{MS}, we choose a special cut-off function as follows. Define $\beta_{\rho}$ to
be the involution of the function 
$$1-\frac{log(t)}{log \rho}.\leqno(3.33)$$
for $t\in [\rho, 1]$ and equal to $0, 1$ for $t<\rho, t>1$
respectively. This function has the property that
$$\int |\bigtriangledown \beta|^2<\frac{C}{-log \rho}.$$
Such a cut-off function was first introduced by Donaldson and Kroheimer \cite{DK} 
in 4-dimension case. We refer to \cite{DK}, \cite{MS} for the discussion of the importance of such
a cut-off function. Then, we define
$$\bar{\beta}_r(t)=\beta(\frac{2t}{r^2}), \leqno(3.34)$$
which is a cut-off function for the annulus $N_p(\frac{\rho r^2}{2}, \frac{r^2}{2}).$
Clearly, $\bar{\beta}_r$ is the convolution of the function
$$1-\frac{log(\frac{2t}{r^2})}{log \rho}.\leqno(3.35)$$

Let $\Sigma^w=dom(f^w)$, where we have already marked the extra marked
$e^v_1, \cdots, e^v_l$ to stabilize the unstable components. Then, we define 
$$f^{v,w}: \Sigma^w_v\rightarrow Y$$
as
$$f^{v,w}=\left\{\begin{array}{ll}
f^w_1(x);  x\in \Sigma_1-D_p(\frac{2r^2}{\rho})\\
\bar{\beta}_r(t)(f^w_1(s, t)-y_w)+f_2(\theta+s, \frac{r^4}{t}); 
   x=re^{i\theta}\in N_p(\frac{\rho r^2}{2}, \frac{r^2}{2})\cong N_q(2r^2, \frac{2r^2}{\rho})\\
 f^w_1(s,t)+f^w_2(\theta+s, \frac{r^4}{t})-y_w; x=re^{i\theta}\in N_p(\frac{r^2}{2},
                              2r^2)\cong N_q(\frac{r^2}{2}, 2r^2)\\
\bar{\beta}_r(t)(f^w_2(s, t)-y_w)+f_1(\theta+s, \frac{r^4}{t}); 
   x=re^{i\theta}\in N_q(\frac{\rho r^2}{2}, r^2)\cong N_p(r^2, \frac{2r^2}{\rho})\\
f^w_2(x); x\in \Sigma_2-D_q(\frac{2r^2}{\rho})
\end{array}\right. \leqno(3.36)$$
where $y_w=f^w_1(p)=f^w_2(q).$ 
To get an element of $\overline{\B}_A(Y,g,k)$,
we have to view $f^{w,v}$ as a function  $\pi_{k+l}(\tilde{\Sigma}_v)$ by
forgetting the extra marked points. We denote it by $\bar{f}^{w,v}$.

  There is a right inverse of the map $f^{v,w}$ defined as follows. Suppose that 
$$f: \Sigma^w_v\rightarrow Y.\leqno(3.37)$$
Let $\tilde{beta}_r(t)$ be a cut-off function on the interval $(\frac{r^2}{2},
2r^2)$, which is symmetry with respect to $t=r^2$. Namely, 
$$\tilde{\beta}_r(t)=1-\tilde{\beta}_r(-2t+3r^2), \mbox{ for $t<r^2$ }.$$
 We define 
$$f_v=(f^1_v, f^2_v): \Sigma^w_1\wedge \Sigma^w_2\rightarrow Y.\leqno(3.38)$$
by
$$f^1_v=\left\{\begin{array}{lll}
        f(x); x\in \Sigma_1-D_p(2r^2)\\
        \tilde{\beta}_r(f(x)-\frac{1}{2\pi r^2}\int_{S^1} f(s, r^2))+
           \frac{1}{2\pi r^2}\int_{S^1} f(s, r^2); x\in D_p(2r^2)
      \end{array} \right. \leqno(3.39)$$
$$f^2_v=\left\{\begin{array}{lll}
        f(x); x\in \Sigma_2-D_q(2r^2)\\
        \tilde{\beta}_r(f(x)-\frac{1}{2\pi r^2}\int_{S^1} f(s, r^2))+
           \frac{1}{2\pi r^2}\int_{S^1} f(s, r^2); x\in D_q(2r^2)
      \end{array} \right. \leqno(3.40)$$
Roughly speaking, we cut the $f$ over the annulus with $\frac{r^2}{2}<t<2r^2.$

   By the construction, the attaching map is really the composition of two maps.
The intermediate object is
$$\U_f=\bigcup_{\tilde{\Sigma}_v\in \tilde{O}_{f}}\{exp_{f^v}\{w\in
\Omega^0((f^v)^*T^*_FY); w\perp E_{e^f_i}, ||w||_{L^p_1}<\epsilon, ||w||_{C^{1}(D_{\delta_0}(g(e_i)))}< \epsilon, g\in
sbt_f\}\}.\leqno(3.41)$$
$\U_f$ is clearly a stratified Frechet V-manifold.
Then,
$$f^{.,.}: \U_{f,D}\times \C^{\epsilon}_f\rightarrow \U_f\leqno(3.42)$$
and
$$\bar{\{.\}}:  \U_f\rightarrow \overline{\B}_A(Y,g,k).\leqno(3.43)$$
Let $\tilde{\U}_f=Im(\U_f)$ under $\bar{\{.\}}$.

 The different gluing 
parameters give rise to different $\tilde{\Sigma}_v\in \overline{\M}_{g, k+l}$.
However, we want to study the injectivity of attaching map, where we have to
consider $\bar{f}^{w,v}$. It would be more convenient to construct $\pi_{k+l}(
\tilde{\Sigma}_v)$ directly. We shall give such an equivalent description of
gluing process. 

Recall that the domain of a stable map can be constructed by
first adding a chain of $\P^1$'s to separate double point and then add trees of
$\P^1$'s. Now we distinguish principal components and bubble components in our
construction. We first glue the principal components. In this case, the different 
gluing parameters give rise to the different marked Riemann surfaces. Then,
we glue the maps according to formula (3.36). When we glue a bubble component, we
gives an equivalent description. Suppose that $\Sigma_i$ is a stable Riemann
surface and $\Sigma_j$ is a bubble component. Moreover, $\Sigma_i, \Sigma_j$
intersects at $p\in \Sigma_i, q\in \Sigma_j$. Suppose that the gluing parameter
is $v=re^{i\theta}$. We can view the previous construction as follow. We cut
off the balls $D^p_i\subset \Sigma_i, D^q_j\subset \Sigma_j$ of radius $\frac{2r^2}{\rho}$ centered at 
marked points we want to glue. The complement $\Sigma_j-D_j$ is conformal
equivalent to a ball of radius $\frac{2r^2}{\rho}$. Then, we glue back the disc along the 
annulus by rotating 
angel $\theta$. Clearly, this is
just a different parameterization of $\Sigma_i$. But we do obtain a holomorphic
map from $\Sigma_i\#_v \Sigma_j$ to $\Sigma_i$. Furthermore, we 
obtain a local universal family 
$$\bar{\tilde{U}}_{f}\rightarrow \bar{\tilde{O}}_{f}\leqno(3.44)$$
of $\Sigma=dom(f)$ as an element of $\overline{\M}_{g,k}$. Although 
 $\Sigma_i\#_v \Sigma_j$ is just $\Sigma_i$ in our alternative gluing construction,
 the different 
gluing parameters may give different maps. Let $\tau_v$ be the composition of
 rescaling and rotation conformal transformations described above.
Let $e_{i}$ be the marked
points of $\Sigma_j$ other than $q$. We observe that $\tau_v$ rescaled $|df(e_i)|$
at the order $\frac{1}{r^2}$. Then, we repeat above construction for each
bubble component.
\vskip 0.1in
\noindent
{\bf Lemma 3.8: }{\it Suppose that $\bar{f}^{v,w}=\bar{f}^{v',w'}$. Then, 
$$v=v'. \mbox{ mod }(stb_f) \leqno(3.45)$$}
\vskip 0.1in

As we
mentioned above, $\Sigma^w_v\neq \Sigma^{w'}_{v'}$ if
$v\neq v'$. If $\pi_{k+l}(\Sigma^w_v)\neq \pi_{k+l}(\Sigma^{w'}_{v'})$, 
$$\bar{f}^{v,w}\neq \bar{f}^{ v', w'}\leqno(3.46)$$
by the definition. If $\pi_{k+l}(\Sigma^w_v)=\pi_{k+l}(\Sigma^{w'}_{v'})$,
there are two possibilities. Since $\pi_{k+1}(\Sigma^w_v)$ is the quotient of $\bar{\Sigma}^w_v$
by $stb_{\bar{\Sigma}^w_v}$, either  $\bar{\Sigma}^w_v=\bar{\Sigma}^{w'}_{v'}$ or 
they are different by an element  of
$stb_{\bar{\Sigma}^w_v}\subset stb_f$. Since the attaching map
is invariant under $stb_f$, we can apply this element to $(w', v')$. Therefore,
we can just simply assume that $\bar{\Sigma}^w_v=\bar{\Sigma}^{w'}_{v'}$. 
On the other hand, $\Sigma^w_v$ is just $\bar{\Sigma}^w_v$ with 
additional marked points $e^v_1, \cdots, e^v_l$. Then, it is enough
to show that 
$$e^v_i=e^{v'}_i. \mbox{ mod }(stb_f) \leqno(3.47)$$

Suppose that $\Sigma_j$ contains extra marked point $e_s$.
We choose small $r$ such that 
$$\frac{1}{r^2}>>\frac{max \{|df^w_1| |df^{w'}_1|\}}{min \{|df^{w}_2(e_s)|,
|df^{w'}_2(e_{s})|\}}.$$
When $\epsilon$ is small, $|df^{w}_2(e_s)|, |df^{w'}_2(e_s)|>0$. Therefore, we can
assume that
$$|d(\tau_vf)^{w}_2(e_s)|, |d(\tau_v)f^{w'}_2(e_s)|>max \{|df^w_1|, |df^{w'}_1|\}.
\leqno(3.48)$$
Hence, $\tau(e^v_{s}), \tau(e^{v'}_{s})\in D^p_i\cap D^{p'}_j, \tau\in stb_{f_i}$. Furthermore,
$$\tau_v f^w=\tau_{v'}f^{w'}.\leqno(3.49)$$
on a smaller open subset $D_0$ of $D^p_i\cap D^{p'}_j$ containing 
$\tau(e^v_s), \tau(e^{v'}_s)$. 
Hence, 
$$f^{w'}=\tau^{-1}_{v'}\tau_v f^w.\leqno(3.50)$$
on an open set containing $e_s$. However, both $f^w, f^{w'}$ are in the slice 
$W_f$. Hence, (3.50) is valid for $f^w, f^{w'}$ over a component of $\Sigma_f$ containing
$e^v_s, e^{v'}_s$. Hence
$$\tau^{-1}_{v'}\tau_v\in stb_f.\leqno(3.51)$$
Therefore,
$$e^v_s=e^{v'}_s. \mbox{ mod }(stb_f) \leqno(3.52)$$
Furthermore, we also observe that
$$f^{w}=f^{w'} \mbox{ on } \Sigma-\bigcup D_{ij}. \leqno(3.53)$$
$\Box$

$\bar{.}$ is obviously  invariant under $stb_f$. Moreover,
\vskip 0.1in
\noindent
{\bf Lemma 3.9: }{\it The induced map of $\bar{\{.\}}$ from $ \U_f/stb_f$ to 
$\tilde{\U}_f\subset \overline{\B}_A(Y,g,k)$ is one-to-one. Furthermore, the 
intersection of $\tilde{\U}_f$ with each strata is open and homeomorphic to the corresponding
strata of $\U_f$.}
\vskip 0.1in
{\bf Proof:} Let
$$\V_f=\U_{f,D}\times \C_f.$$
By (3.39), (3.40), $f^{v,w}$ is onto. Suppose that $\bar{f}^{v,w}=\bar{f}^{v',w'}$. By
the Lemma 3.8, $v=v'$ mod($stb_f$). Therefore, we can assume that $v=v'$. 
Moreover, we can assume that $\Sigma^w_v=\Sigma^{w'}_{v'}.$ However, it is
obvious that
$$\bar{.}: Map^F_A(\Sigma^w_v)\rightarrow \overline{\B}_A(Y,g,k)$$
is injective. So we show that 
$$f^{w,v}=f^{w',v'}.\leqno(3.54)$$
To prove the second statement, let $w_0\in \Omega^0(f^*T_FY)$ with $w_0\perp E_{e^f_i}$.
For any map close to $\bar{f}^{v, w_0}$, it is of the form $f^{v, w_0+w}$ with
 $||w||_{L^p_1}< \epsilon, ||w||_{L^p_1}<\epsilon, ||w||_{C^{1}(D_{\delta_0}(g(e_i)))}< \epsilon$. We want to show that we can perturb $e^f_i$ such that
$$w_0+w\perp E_{e^f_i}.\leqno(3.55)$$  The argument of Lemma 3.4 applies. 

Now we define the topology of $\overline{\B}_A(Y,g,k)$ by specifying the 
converging sequence.
\vskip 0.1in
\noindent
{\bf Definition 3.10: }{\it A sequence of stable maps $f_n$ converges to $f$ if for any
$\tilde{\U}_f$, there is $N>0$ such that if $n>N$ $f_n\in \tilde{\U}_f$. Furthermore, $f_n$ 
converges to $f$ in $C^{\infty}$-topology in any compact domain away from the gluing region.}
\vskip 0.1in
\noindent
{\bf Proposition 3.11: }{\it If a sequence of stable holomorphic
maps weakly converge to $f$ in the sense of \cite{RT1}, they converge to $f$ in the topology 
defined in the Definition 3.8.
}
\vskip 0.1in
  The proof is delayed after Lemma 3.18.
\vskip 0.1in
Define
$$\chi: \overline{\B}_A(Y,g,k)\rightarrow \overline{\M}_{g,k} \leqno(3.56)$$
by $\chi(f)=\pi_{k+l}(dom(f)).$
\vskip 0.1in
\noindent
{\bf Corollary 3.12: }{\it  $\chi$ is continuous.}
\vskip 0.1in
  The proof follows from the definition of the topology of $\overline{\B}_A(Y,g,k)$.
\vskip 0.1in
\noindent
{\bf Theorem 3.13: }{\it $\overline{\B}_A(Y,g,k)$ is Hausdorff.}
\vskip 0.1in
{\bf Proof: } Suppose that $f\neq f'$. By the corollary 3.12, we can assume that 
$\pi_{k+l}(dom(f))=\pi_{k+l}(dom(f')$.
We want to show that $\tilde{\U}_f\cap \tilde{\U}_{f'}=\emptyset$ for some $\epsilon$. 
Suppose that it is false. We claim that $dom(f), dom(f')$ have the same topological type. 
Namely, $f,f'$ are in the same strata. We start from the underline stable Riemann surfaces
$\pi_{k+l}(dom(f))=\pi_{k+l}(dom(f'))$ which are the same by the assumption. We want to show that
they always have the same way to attach bubbles to obtain $dom(f), dom(f')$. Suppose that
we attach a bubble to $\pi_{k+l}(dom(f))$ at $p$. Recall that the energy concentrates at $D_p(\frac{2r^2}{\rho})$,i.e., $\int_{D_p(\frac{2r^2}{\rho})}|df|^2\geq \epsilon_0$.
 The same is
true for $f^{w,v}$ when $||w||_{L^p_1}<\epsilon$. On the other hand, we have the same property for
$(f')^{w',v'}$ for some $||w'||_{L^p_1}, |v'|<\epsilon$. If $\bar{f}^{w,v}=\bar{f'}^{w',v'}$, 
$f'$ must have  a bubbling point in $D_p(\frac{2r^2}{\rho})$. In fact, the bubbling point must be $p$. Otherwise, 
we can construct a small ball $D_p(\frac{2r^2}{\rho})$ 
containing no bubbling points of $f'$. 
Then, we proceed inductively on the next bubble. Now the energy concentrates at a ball of
radius $r^2r^2_1$, where $r_1=|v_1|$ is the next gluing parameter. By the induction, we can
show that $dom(f), dom(f')$ have the same topological type. In fact, we proved that $dom(f), dom(f')$
have the same bubbling points and hence the same holomorphic type.

   Suppose that $f, f'\in \B_D(Y,g,k)$. Then, some component of $f,f'$ are different. Suppose
that the component $f_i\neq f'_i$, where $f_i, f'_i\in \B_{A_i}(Y,g,k)$. Note that $f^v_i$ is equal
to $f$ outside the gluing region. $f^v\neq (f')^v$ for small $v$. By Lemma 3.4,
$\B_{A_i}(Y,g,k)$ is Hausdorff and the neighborhoods of $f_i, f'_i$ are described by slice
$W_{f_i}, W_{f'_i}$ for a small constant $\epsilon$. Add extra
marked points to stabilize unstable components.  $||f_i-f'_i||_{L^p_1}\geq 2\epsilon$ for
small $\epsilon$.  Then, it is obvious that
$$W_{f_i}\cap W_{f'_i}=\emptyset.\leqno(3.57)$$
Note that $f^{w,v}(e_0)=f^w(e_0), f^{w',v'}(e_0)=f^{w'}(e_0)$. It is straightforward to check that
$$\tilde{\U}_f \cap \tilde{\U}_{f'}=\emptyset\leqno(3.58)$$
for the same $\epsilon$. This is a contradiction. $\Box$
\vskip 0.1in
\noindent
{\bf Corollary 3.14: }{\it $\overline{\M}_A(Y,g,k)$ is Hausdorff.}
\vskip 0.1in
To construct the obstruction bundle $\overline{\F}_A(Y,g,k)$, we start from the top
strata $\B_A(Y,g,k)$. Let $\V(Y)$ be vertical tangent bundle. With an almost
complex structure $J$, we can view $\V(Y)$ as a complex vector bundle. 
Therefore, for each $f\in \B_A(Y,g,k)$ we can decompose
$$\Omega^1(f^*\V(Y))=\Omega^{1,0}(f^*\V(Y))\oplus \Omega^{0,1}(f^*\V(Y)).
\leqno(3.59)$$
Both bundles patch together to form Frechet V-bundles over $\B_A(Y,g,k)$. We 
denote them by $\Omega^{1,0}(\V(Y)), \Omega^{0,1}(\V(Y))$. Then,
$$\F_A(Y,g,k)=\Omega^{0,1}(\V(Y)).\leqno(3.60)$$
For lower strata $\B_D(Y,g,k)$, $\B_D(Y,g,k)\subset \prod_i \B_{A_i}(Y,g_i,
k_i)$, where $\B_{A_i}(Y,g_i,k_i)$ are components. When a component is stable,
we already have an obstruction bundle $\F_{A_i}(Y,g,k)$. When the i-th 
component is unstable, we first form the obstruction bundle over $Map^F_{A_i}(
Y,,0,k_i)$ in the same way and divide it by $Aut_i$. In the quotient, we obtain
a V-bundle denoted by $\Omega^{0,1}(\V(Y))$. Let
$$i: \B_D(Y,g,k)\rightarrow \prod_i \B_{A_i}(Y,g_i,k_i)\leqno(3.61)$$
be inclusion. We define
$$\F_D(Y,g,k)=i^*\prod_i \F_{A_i}(Y,g_i,k_i).\leqno(3.62)$$
Finally, we define
$$\overline{\F}_A(Y,g,k)|_{\B_D(Y,g,k)}=\F_D(Y,g,k).
\leqno(3.63)$$
 For any $f\in \B_D(Y,g,k)$, consider a chart
$(\U_f, V_f, stb_f)$. Suppose that $D=\Sigma_1
\wedge \Sigma_2$. For $\eta^w\in \Omega^{0,1}((f^w)^*\V(Y))$, define
$$\eta^{w,v}\in \Omega^{0,1}((f^{w,v})^*\V(Y))$$
by
$$\eta^{w, v}=\left\{\begin{array}{ll}
\eta_1(x);  x\in \Sigma_1-D_p(\frac{2r^2}{\rho})\\
\bar{\beta}_r(t)\eta_1(s, t)+\eta_2(\theta+s, \frac{r^4}{t}
); 
   x=te^{is}\in N_p(\frac{\rho r^2}{2},\frac{r^2}{2})\cong N_q(2r^2, \frac{2r^2}{\rho})\\
   \eta_1(s,t)+\eta_2(\theta+s, \frac{r^4}{t}); x=te^{is}\in N_p(\frac{r^2}{2},2r^2)
\cong N_q(\frac{r^2}{2},2r^2)\\
\bar{\beta}_r(t)\eta_2(s, t)+\eta_1(\theta+s, \frac{r^4}{t}
); 
   x=te^{is}\in N_q(\frac{\rho r^2}{2}, \frac{r^2}{2})\cong N_p(2r^2, \frac{2r^2}{\rho})\\
\eta_2; x\in \Sigma_2-D_q(\frac{2r^2}{\rho})
\end{array}\right. \leqno(3.64)$$
$\bar{\partial}_J$ is clearly a continuous section of $\bar{\F}_A(Y,g,k,J)$. Let
$\bar{\partial}_{J, D}$ be the restriction of $\bar{\partial}_J$ over $\B_D$.

  Next, we define the local sections by repeating the constructions in section
2. Let $f\in \B_D(Y,g,k)$.  $Coker D_f \bar{\partial}_{J,D}$ is a finite dimensional
 subspace of
$\Omega^{0,1}(f^*\V(Y))$ invariant under $stb_f$. We first choose a  
$stb_f$-invariant cut-off function vanishing in a small neighborhood of the
intersection points. Then we multiple it to the element of $Coker D_f \bar{
\partial}_{J,D}$ and denote the resulting finite dimensional space as $F_f$. By the
construction, $F_f$ is $stb_f$-invariant. When the support of the cut-off function
is small, $F_f$ will have the same dimension as $Coker D_f \bar{\partial}_J$ and
$$D_f\bar{\partial}_{J,D}+Id: \Omega^0(f^*T_FY)\oplus F_f\rightarrow \Omega^{0,1}(
f^*\V(Y))$$
is surjective.  We first extend each element $s$ of $F_f$ to a smooth 
 section  $s^w\in \Omega^{0,1}((f^w)^*\V(Y))$ of $\F_D(Y,g,k,J)$ supported in $U_{f,D}$ such
that it's value vanishes in a neighborhood of the intersection points. Hence, $s^w$ can be
naturally viewed as an element of $\Omega^{0,1}((f^{w,v})^*\V(Y))$ supported 
away from the gluing region.
Let $\beta_f$ be a smooth cut-off function on a polydisc $\C_f$ vanishing 
outside of a polydisc of radius $2\delta_1$ and equal to 1 in the polydisc of
radius $\delta_1$. One can construct $\beta_f$ by first constructing such $\beta$
over each copy of gluing parameter $\C_{x}$ and then multiple them together. We
now extend $s^w$ over $\U_f$ by the map
$$s^{v}_c(f^{w,v})=\beta_f(v)s^{w}.\leqno(3.65)$$
Then, we use the method of the section 2 (2.5) to extend the identity map of $F_f$ to
a map
$$s_f: F_f \rightarrow \overline{\F}_A(Y,g,k)|_{\U_f}.\leqno(3.65.1)$$
invariant under $stb_f$ and supported in $\U_f$. Then, it descends to a map 
over $\overline{\B}_A(Y,g,k)$. We will use $s_f$ to denote the
induced map on $\overline{\B}_A(Y,g,k)$ as well. We call such $s_f$ {\em 
admissible}. 
Our new equation 
will be of the form 
$$\S_e=\bar{\partial}_f+\sum_i  s_{f_i}: \E \rightarrow \overline{\F}_A(Y,g,k,J),
\leqno(3.66)$$
where $s_{f_i}$ is  admissible. We observe that the restriction $\S_D$ of
 $\S$
over each strata is smooth. Let $U_{\S_e}=(\S_e)^{-1}(0)$  and
$$S: U_{\S_e}\rightarrow E.$$
\vskip 0.1in
{\bf Lemma 3.15: }{\it $S$ is a proper map.}
\vskip 0.1in
{\bf Proof: } Since the value of $s_{f_i}$ is supported away from the 
gluing region, the proof of lemma is completely same as the case to show
that the moduli space of stable holomorphic maps is compact. We omit it. $\Box$

For $f\in \B_D(Y,g,k)$, we define the tangent space
$$T_f\overline{\B}_A(Y,g,k)=T_f \B_D(Y,g,k)\times \C_f$$
and  the derivative
$$D_{f,t}\S_e=D_{f,t} \S_e|_{\B_D(Y,g,k)}: T_f\overline{\B}_A(Y,g,k)\rightarrow \Omega^{0,1}(f^*
\V(Y)).\leqno(3.67)$$
\vskip 0.1in
\noindent
{\bf Lemma 3.16: }{\it 
$$Ind D_{f,t}\S=2C_1(V)(A)+2(3-n)(g-1)+2k+\dim X+dim E.\leqno(3.68)$$
}
\vskip 0.1in
{\bf Proof:} 
$$D_{f,t}\S_D(W,u)=D_f\bar{\partial}_J(W)+\sum_i  D_{f,t}s_{f_i}(W,u).
\leqno(3.69)$$
$$Ind D_{f,t}\S_D=Ind D_f\bar{\partial}_J + dim E.$$
If $\Sigma_f=dom(f)$ is irreducible, the lemma follows from
Riemann-Roch theorem. Suppose that $\Sigma_f=\Sigma_1\wedge \Sigma_2$ and
$f=(f_1, f_2)$ with $f_1(p)=f_2(q)$.
$$\begin{array}{lll}
Ind D_f\bar{\partial}_J&=&Ind D_{f_1}\bar{\partial}_J+Ind D_{f_2}\bar{
\partial}_J-dim Y\\
&=&2C_1(V)([f_1])+2(3-n)(g_1-1)+2(k_1+1)+\dim X+2C_1(V)([f_2])+2(3-n)(g_2-1)\\
&&+2(k_2+1)+\dim X-\dim Y\\
&=&2C_1(V)(A)+2(3-n)(g-1)+2k+\dim X-6+2n+2-2n\\
&=&2C_1(V)(A)+2(3-n)(g-1)+2k+\dim X-2
\end{array}$$
Adding the dimension of gluing parameter, we derive Lemma 3.16. The general case
can be proved inductively on the number of the components of $\Sigma_f$. We omit 
it.

This is the end of the construction of the extended equation. Next, we shall prove
that 
$$(\overline{\B}_A(Y,g,k), \overline{\F}_A(Y,g,k), \bar{\partial}_J)\leqno(3.70)$$ 
is
VNA. The openness of
 $\U_S=\{(x,t); Coker D_{f,t}\S_e=\emptyset\}$ is a local property. To prove
the second property, we first construct a local coordinate chart for each point of virtual
neighborhood. Then, we prove that the local chart patches together to form a
$C^1$-V-manifold. The construction of a local coordinate chart is  basically a 
gluing theorem. The first gluing theorem for 
pseudo-holomorphic curve was given by \cite{RT1}. There were two new proofs by
\cite{Liu}, \cite{MS} which are more suitable to the set-up we have here. Here we follow that
of \cite{MS}. For reader's convenience, we outline the proof here.

We need to enlarge our space to include Sobolev maps. Suppose that $f\in
\M_D(Y,g,k), t_0\in \R^m$ such that $\S_e(f,t_0)=0$ and $Coker \D_{f,t_0}\S_e=0$.
Choose metric $\lambda$ on 
$\Sigma_1\wedge \Sigma_2$. Using the trivialization of (3.18), we can define Sobolev norm on $\U_{f,D}$. Let
$$L^p_1(\U_{f,D})=U_{\Sigma}\times\{f^{w}; w\in \Omega^0(f^*T_FY), ||w||_{L^p_1}<
\epsilon, ||w||_{C^1(D_{\delta_0}(g(e_i)))}< \epsilon, w\perp E_{e^f_i}\}.\leqno(3.71)$$
By choosing small $\delta_0$, we can assume that $D_{\delta_0}(e_i)$ is away from gluing
region. For the rest of this section, we assume that $2<p<4$. Then, $L^p_1(\U_{f,D})$ is 
a Banach manifold. To simplify the notation, we shall
assume that $dom(f)=\Sigma_1\wedge \Sigma_2$ for the argument below. However, it
is obvious that the same argument works for the general case. 
Let $\lambda_v$ be the metric on $\Sigma_v$ defined in (3.20).
We use $L^p_v, L^p_{1,v}$ to denote the Sobolev norms on $\Sigma_v$, where $v$ is used
to indicate the dependence on $v$. By \cite{MS} (Lemma A.3.1), the Sobolev constants of the 
metric $\lambda_v$ are
independent of $v$.
Let
$$L^p_1(\U_f)=\bigcup_{\tilde{\Sigma}_v}\{f^{v,w}; w\in \Omega^0((f^v)^*T_FY),
  w\perp E_{e^f_i}, ||w||_{L^p_{1,v}}<\epsilon, ||w||_{C^1(D_{\delta_0}(g(e_i)))}< \epsilon\}.\leqno(3.72)$$

  First of all, the map
$$f^{w,v}: \U_{f,D}\times \C_f \rightarrow \U_f$$
induces a natural map
$$\phi_f: \Omega^0((f^w)^*T_FY)\rightarrow \Omega^0((f^{w,v})^*T_FY)$$
by the formula
$$u^{w,v}=\phi_f(u)=\left\{\begin{array}{ll}
u_1(x);  x\in \Sigma_1-D_p(\frac{2r^2}{\rho})\\
\bar{\beta}_r(t)(u_1(s, t))-u_1(0))+u_2(\theta+s, \frac{
         r^4}{t}); 
     x=re^{is}\in N_p(\frac{\rho r^2}{2}, \frac{r^2}{2})\cong N_q(2r^2, \frac{2r^2}{\rho})\\
        u_1(s,t)+u_2(\theta+s, \frac{r^4}{t})-u(0); x=re^{is}\in N_p(\frac{r^2}{2}, 
        2r^2)\cong N_q(\frac{r^2}{2}, 2r^2)\\
\bar{\beta}_r(t)(u_2(s, t))-u_2(0))+u_1(\theta+s, \frac{
         r^2}{t}); 
     x=re^{is}\in N_q(\frac{r}{2}, r)\cong N_p(r, 2r)\\
u_2(x); x\in \Sigma_2-D_q(2r)
\end{array}\right. \leqno(3.73)$$
where $u=(u_1, u_2)\in \Omega^0((f^w)^*T_FY)$. Notes that $u_1(0)=u_2(0)$.

One can construct an inverse of $\psi_f$. For any $u\in
\Omega^0((f^{w,v})^*T_FY)$,
we define 
$$u_v=(u^1_v, u^2_v)$$
by
$$u^1_v=\left\{\begin{array}{lll}
        u(x); x\in \Sigma_1-D_p(2r^2)\\
        \tilde{\beta}_r(u(x)-\frac{1}{2\pi r^2}\int_{S^1} u(s, r^2))+
           \frac{1}{2\pi r^2}\int_{S^1} u(s, r^2); x\in D_p(2r^2)
      \end{array} \right. \leqno(3.74)$$
$$u^2_v=\left\{\begin{array}{lll}
        u(x); x\in \Sigma_2-D_q(2r^2)\\
        \tilde{\beta}_r(u(x)-\frac{1}{2\pi r^2}\int_{S^1} u(s, r^2))+
           \frac{1}{2\pi r^2}\int_{S^1} u(s, r^2); x\in D_q(2r^2)
      \end{array} \right. \leqno(3.75)$$

For any $\eta\in \Omega^{0,1}((f^{w,v})^*\V(Y))$, we cut $\eta$ along the circle
of radius $r^2$ and extend as zero inside the $D_p(r^2), D_q(r^2)$. We denote 
resulting 1-form as $\eta^f_1\in \Omega^{0,1}((f^w)^*
\V(Y)), \eta^f_2\in \Omega^{0,1}(f^w)^*\V(Y))$. Clearly, $(\eta^f_1, \eta^f_2)$ is
an right inverse of $\eta^{w,v}$. 
\vskip 0.1in
\noindent
{\bf Lemma 3.17: }{\it Let $u$ be a 1-form over a disc of radius $\frac{2r^2}{\rho}<1$. Then,
$$||\bigtriangledown \bar{\beta}_r (u-u(0))||_{L^p}\leq c|log \rho|^{1-\frac{4}{p}}||u||_{L^p_1}.
\leqno(3.76)$$
}
\vskip 0.1in
The inequality is just the lemma A.1.2 of \cite{MS}, where we use $r^2$ instead of $r$. 
\vskip 0.1in
\noindent
{\bf Lemma 3.18: }{\it $||\phi_f(u^w)||_{L^p_{1,v}}\leq C 
||u^w||_{L^p_1}, ||u^i_v||_{L^p_1}\leq C||u||_{L^p_{1,v}}$.}
\vskip 0.1in
{\bf Proof: } We only have to consider $u^w$ over $N_p(\frac{\rho r^2}{
2}, \frac{r^2}{2})$, where 
$$\phi_f(u_w)=\bar{\beta}_r(t)(u^w_1(s,t)-u^w_1(0))+u^w_2(s+\theta,
\frac{r^4}{t}).\leqno(3.77)$$
$$\begin{array}{lll}
||\phi_f(u^w)||_{L^p(N_p(\frac{\rho r^2}{2}, \frac{r^2}{2}))}&\leq &C(||u^w_1||_{L^p_1(
N_p(\frac{\rho r^2}{2}, \frac{r^2}{2}))}+||u^w_2||_{L^p_1(N_q(2r^2,\frac{2r^2}{\rho}))}+|u^w_1(0)|)\\
 &\leq& C(||u^w_1||_{L^p_1(N_p(\frac{\rho r^2}{2},\frac{r^2}{2}))}+||u^w_2||_{L^p_1(N_q(2r^2,\frac{2r^2}{\rho}))}).
\end{array}\leqno(3.78)$$
$$\begin{array}{lll}
&&||\bigtriangledown \phi_f(u^w)||_{L^p(N_p(\frac{\rho r^2}{2},\frac{r^2}{2}))}\\
&\leq &C(||\bigtriangledown 
u^w_1||_{L^p_
1(N_p(\frac{\rho r^2}{2}, \frac{r^2}{2}))}+||\bigtriangledown u^w_2||_{L^q_1(N_q(2r^2,\frac{2r^2}{\rho}))}+
||\bigtriangledown\bar{
\beta}_r (u^w_1-u^w_1(0))||_{L^p(N_p(\frac{\rho r^2}{2}, \frac{r^2}{2}))})\\
&\leq& C||u^w||_{L^p_1(N_p(\frac{r^2}
{4},\frac{r^2}{2}))},
\end{array}\leqno(3.79)$$
 where the last inequality follows from Lemma 3.17. The proof of the second inequality is
the same and we omit it. $\Box$

{\bf Proof of Proposition 3.11: } Suppose that $f_n\rightarrow f$ as a weakly 
convergent sequence of holomorphic
stable maps in the  sense of \cite{RT1}. Then, $f_n$ converges to $f$ in $C^{\infty}$-norm
in any compact domain outside the gluing region, in particular on $D_{\delta_0}(g(e_i))$.
Now, we want to show that $f_n$ is in the
open set $\U_{f, D}$  for $n>N$. Note that formula (3.74,3.75) is a left
inverse of formula (3.73). By Lemma 3.18, the formula (3.73) preserves $L^p_1$
norm. Hence, it is enough
to show that $f_n$ is close to $f^v$ when $n$ is large. Namely, we want to
estimate $||f_n-f^v||_{L^p_{1,v}}$. Outside of gluing region, $f_n$ converges
to $f^v$ in the $C^{\infty}$ norm. So $||(1-\beta)(f_n-f^v)||_{L^p_{1,v}}$ 
converges to zero, where $\beta$ is a cut-off function vanishing outside gluing
region.  Over the gluing region, it is enough
to show that $||\beta(f_n-pt)||_{L^p_1}$ is small where $pt$ is the intersection
point of two components of $f$. Here we assume that $f$ has only two 
components to simplify the notation. The argument for general case is the same. By the decay
 estimate in \cite{RT1}(Lemma 6.10), $||f_n-pt||_{C^0}$ converges to zero
over the gluing region with  cylindric metric. However, $C^0$-norm is independent of the metric of domain.
Hence, we have a $C^0$ estimate for the metric in this paper. Furthermore, 
$f_n$ is holomopophic. By elliptic estimate,
$$||\beta(f_n-pt)||\leq c(||\bar{\partial}_J(\beta(f_n-pt)||_{L^p_v}+||\beta(
f_n-pt)||_{C^0}\leq c(||\bigtriangledown \beta(f_n-pt)||_{L^p_v}+||f_n-pt||_{
C^0}\leq c||f_n-pt||_{C^0}.$$
We will finish the argument
by showing that the constant in elliptic estimate is independent of the
gluing parameter $v$. The later is easy since our metric is essentially 
equivalent to the metric on the anulus N(1,r) in $R^2$, where $r=|v|$ and
$\beta(f_n-pt)$ is compact supported. $\Box$

Suppose that $D_{f,t_0}\S_e$ is surjective. Since $\S_e$ is smooth over $\B_D(Y,g,k)$,
$D_{f^w,t}\S_e$ is surjective for $||w||_{L^p_1}<\delta, |t-t_0|<\delta$ with some small
$\delta$. We choose a family of right inverse $Q_{f^w,t}$. Then,
$$||Q_{f^w,t}||\leq C.\leqno(3.79.1)$$
We want to construct right inverse of $D_{f^{w,v},t}\S_e$. 
\vskip 0.1in
\noindent
{\bf Definition 3.19: }{\it Define $AQ_{f^{w,v},t}(\eta)=\phi_fQ_{f^{w},t}(\eta^f_1, 
\eta^f_2)$.}
\vskip 0.1in
Then, it was shown in \cite{MS} that 
\vskip 0.1in
\noindent
{\bf Lemma 3.20:}{\it $$||AQ_{f^{w,v},t}||\leq C, ||D_{f^{w,v},t}AQ_{f^{w,v},t}-
Id||<\frac{1}{2} \mbox{ for small } r, \rho. \leqno(3.80)$$}
\vskip 0.1in
   Now, we fix a $\rho$ such that Lemma 3.20 holds. 

The right inverse of $D_{f^{w,v},t}$ is given by
$$Q_{f^{w,v},t}=AQ_{f^{w,v},t}(D_{f^{w,v},t}AQ_{f^{w,v},t})^{-1}. \leqno(3.81)$$
Furthermore,
$$||Q_{f^{w,v},t}||\leq C.\leqno(3.82)$$
Therefore, we show that
\vskip 0.1in
\noindent
{\bf Corollary 3.21:}{\it $$\U_{\S_e}=\{(x,t); Coker D_{f,t}\S_e=\emptyset\}$$
is open.}
\vskip 0.1in
Next, we have an estimate of error term.
\vskip 0.1in
\noindent
{\bf Lemma 3.23: }{\it Suppose that $\S_e(f^w)=0$. Then,
$$||\S_e(f^{v,w})||_{L^p_v}\leq C r^{\frac{4}{p}}.\leqno(3.83)$$}
\vskip 0.1in
{\bf Proof: } It is clear that $\S_e(f^{v,w})=0$ away from the gluing region.
Notes that the value of $s_{f_i}$ is supported away from the gluing region. Hence, 
$\S_e=\bar{\partial}_J$ over the
gluing region. Then, the lemma follows from \cite{MS} 
(Lemma A.4.3). $\Box$

Next we construct the coordinate charts of $\M_{\S_e} \cap \U_{\S_e}$. Suppose that $(f,t_0)\in 
\M_{\S_e}\cap \U_{\S_e}$. By the previous argument,
we can assume that some neighborhood $\U_f\times B_{\delta}(t_0)\subset \U_{\S_e}$. 
To simplify the notation, we drop $t$-component. It is understood that $s_{f_i}$ will
not affect the argument since it's value is supported away from the gluing region.
Since $L^p_1(\U_{f,D})$
is a Banach manifold and the restriction to $\S_e$ is a Fredholm map, $\M_{\S_e}\cap \B_D(Y,g,k)$
is a smooth V-manifold by ordinary transversality theorem. Let 
$$f\in E^D_f\subset \M_{\S_e} \cap \B_D(Y,g,k)\leqno(3.84)$$ 
be a small $stb_f$-invariant neighborhood. 
\vskip 0.1in
\noindent
{\bf Theorem 3.24: }{\it There is a one-to-one continuous map
$$\alpha_f: E^D_f \times B_{\delta_f}(\C_f) \rightarrow \U_f \leqno(3.85)$$
such that $im(\alpha_f)$ is an open neighborhood of $f\in \M_{\S_e}$, where $\delta_f$ is a small
constant.}
\vskip 0.1in
{\bf Proof:}
  For any $w\in E^D_f$ and small $v$, we would
like to find an element $\xi(w,v)\in \Omega^0((f^v)^*T_FY)$ with $\xi\perp E_{e_i}$
and $\xi(w,v)\in Im Q_{f^{w,v}}$ such that 
$$\S_e((f^{v,w})^{\xi(w,v)})=0.\leqno(3.86)$$
Consider the Taylor expansion
$$\S_e((f^{v,w})^{\xi})=\S(f^{w,v})+D_{f^{w,v}}(\xi)+N_{f^{w,v}}(\xi),$$
for $w\in E^D_f, \xi\in \Omega^0((f^v)^*T_FY)$ with $\xi(e^v_i)\perp df(e^v_i), 
||w||_{L^p_{1,v}}, ||\xi||_{L^p_{1,v}}< \epsilon$. 
Then,
$$\xi(w,v)=-Q_{f^{w,v}}(S(f^{w,v})+N_{f^{w,v}}(\xi(w,v)).\leqno(3.87)$$
Hence, $\xi(w,v)$ is a fixed point of the map
$$H(w,v; \xi)=-Q_{f^{w,v}}(S(f^{w,v})+N_{f^{w,v}}(\xi)).\leqno(3.88)$$
Conversely, if $\xi(w,v)$ is a fixed point, 
$$\S_e((f^{v,w})^{\xi(w,v)})=0.\leqno(3.89)$$
$N_{f^{w,v}}$ satisfies the condition
$$||N_{f^{w,v}}(\eta_1)-N_{f^{w,v}}(\eta_2)||_{L^p_v}\leq C(||\eta_1||_{L^p_{1,v}}+||
\eta_2||_{L^p_{1,v}})||\eta_1-\eta_2||_{L^p_{1,v}}.\leqno(3.90)$$
Next, we show that $H$ is a contraction map on a ball of radius $\delta/4$ for some 
$\delta$. 
$$||H(w,v; \xi)||_{L^p_{1,v}}\leq C(||\S_e(f^{w,v})||_{L^p_v}+
 ||N_{f^{w,v}}(\xi)||_{L^p_v})$$
$$\leq C(r^{\frac{4}{p}}+||\xi||^2_{L^p_{1,v}})\leq \frac{\delta}{4},\leqno(3.91)$$
for $ ||\xi||_{L^p_{1,v}}\leq \frac{\delta}{4}$ and $2C\delta<1, r<(\frac{\delta^2}{4})^{
-\frac{4}{p}}$.
$$||H(w,v; \xi)-H(w,v; \eta)||_{L^p_{1,v}}\leq C ||N_{f^{w,v}}(\xi)-N_{f^{w,v}}(\eta)||_{L^p_v}$$
 $$\leq C (||\xi||_{L^p_{1,v}}+||\eta||_{L^p_{1,v}})||\xi-\eta||_{L^p_{1,v}}<2\delta C
 ||\xi-\eta||_{L^p_{1,v}}.\leqno(3.92)$$
Therefore, $H$ is a contraction map on the ball of radius $\frac{\delta}{4}$. Then,
there is a unique fixed point $\xi(w,v)$. Furthermore, $\xi(w,v)$ depends smoothly
on $w$. Recall that $\xi(w,v)$ is obtained by iterating $H$. One can check that
$$||\xi(w,v)||_{L^p_{1,v}}\leq C r^{\frac{4}{p}}.\leqno(3.93)$$
Our coordinate chart at $f$ is $(E^D_f\times B_{\delta_f}(\C_f), \alpha_f(v,w))$ 
where $\delta_f=(\frac{\delta^2}{4})^{\frac{4}{p}}$. and 
$$\alpha_f(v,w)=(f^{v,w})^{\xi(w,v)}.\leqno(3.94)$$
Notes that all the construction is $stb_f$-invariant. Hence $\alpha_f$ is 
$stb_f$-invariant. It is clear that $\alpha_f$ is one-to-one by contraction mapping principal.
Notes that $\S_e=\bar{\partial}_J$ over the gluing region. It follows from Proposition 3.11
and uniqueness of contraction mapping principal that $\alpha_f$ is surjective onto a 
neighborhood of $f$ in $\M_{\S_e}$.  $\Box$

Furthermore, 
$E^D_f\times \C_f$ has a natural orientation induced by the orientation of
$J$, $\R^m$ and $\C_f$. 

   Next, we show that the transition map is a $C^1$-orientation preserving map. In the previous
argument, we expand $\S_e$ up to the second order, which is given in \cite{F}, \cite{MS}. To
prove the transition map is $C^1$, we need to expand $\S_e$ up to third order.  Let $z=s+it$ be
the complex coordinate of $\Sigma_v$. Let $\bigtriangledown^v \xi=\bigtriangledown_t\xi+
\bigtriangledown_s\xi$ 
to indicate the  dependence on $v$. Let
$$f^{w,v}=exp_{f^v} w^v.\leqno(3.95)$$
Let $ \xi\in L^p_{1,v}(\Omega^0((f^{v})^*T_FY))$ with 
$||w^v||_{L^p_{1,v}}, ||\xi||_{L^p_{1,v}}\leq
\delta$ for small $\delta$.  A similar calculation of 
\cite{MS} (Theorem 3.3.4) implies
$$\bar{\partial}_J(f^{v})^{w^v+\xi})=\bar{\partial}_J(f^{v,w})+D_{f^{v,w}}(\xi)+D^2_{f^{v,w}}(\xi^2)
 +\tilde{N}_{f^{v,w}}(\xi),\leqno(3.96)$$
where
$$D_{f^{w,v}}(\xi)=\bigtriangledown^v_s \xi+J \bigtriangledown^v_t\xi+(C_1 \bigtriangledown w^v+
C_2\bigtriangledown^v f^v+C_2 \bigtriangledown^v w)\xi,\leqno(3.97)$$
$$D^2_{f^{w,v}}\xi=(C_1 \bigtriangledown^v f^v+C_2 \bigtriangledown^v w^v)\xi^2+C_3 \xi
\bigtriangledown^v \xi,\leqno(3.98)$$
$$\tilde{N}_{f^{w,v}}(\xi)=(C_1\bigtriangledown^v f^v+C_2 \bigtriangledown^v w^v)\xi^3+ C_3
(\bigtriangledown^v \xi)\xi^2,\leqno(3.99)$$
where  $C_1, C_2, C_3$ are smooth bounded functions for each of the identities. 
Furthermore, we have
$$D_{(f^v)^{w^v+\tilde{w}}}(\xi)=D_{f^{w,v}}\xi+(2C_1\bigtriangledown^v f^v+2C_2 
\bigtriangledown^v w^v)
 \tilde{w}\xi+C_3 \tilde{w}\bigtriangledown^v \xi+C_4 \tilde{w}\bigtriangledown \xi+ 
O(\tilde{w}^2), \leqno(3.100)$$
where the coefficients of higher order terms are independent from $\tilde{w}$ by (3.99).
\vskip 0.1in
\noindent
{\bf Lemma 3.25: }{\it The derivative with respect to $w$ 
$$||\frac{\partial}{\partial w}D_{f^{v,w}}(\tilde{w})(\xi)||_{L^p_v}\leq C (||f^v||_{L^p_{1,v}}+||w^v||_{
L^p_{1,v}})||\tilde{w}||_{L^p_{1,v}}
||\xi||_{L^p_{1,v}}.\leqno(3.101)$$
$$||\frac{\partial}{\partial w}N_{f^{v,w}}(\tilde{w})(\xi)||_{L^p_v}\leq C(||f^v||_{L^p_{1,v}}
+||w^v||_{L^p_{1,v}})||\tilde{w}||_{L^p_{1,v}}
||\xi||^2_{L^p_{1,v}}.\leqno(3.102)$$}
\vskip 0.1in
{\bf Proof: } The first inequality follows from 3.100. To prove the second inequality, recall that
$$N_{f^{v,w}}(\xi)=\S_e((f^{v})^{w^v+\xi})-\S_e(f^{v,w})-D_{f^{v,w}}(\xi).\leqno(3.103)$$
Hence
$$\begin{array}{lll}
&&N_{(f^{v})^{w^v+\tilde{w}}}(\xi)-N_{f^{v,w}}(\xi)\\
&=&\S_e(f^{v,w^v+\tilde{w}+\xi})-\S_e((f^{v})^{w^v+\xi})-(\S_e((f^{v})^{w^v+\tilde{w}})-
\S_e(f^{v,w}))-(D_{(f^{v})^{w^v+\tilde{w}}}(\xi)-D_{f^{v,w}}(\xi))\\
&=&D_{(f^{v})^{w^v+\xi}}(\tilde{w})-D_{f^{v,w}}(\tilde{w})-\frac{\partial}{\partial w}D_{f^{v,w}}(\tilde{w})(\xi)+O(\tilde{w}^2)\\
&=&\frac{\partial}{\partial w} D_{f^{v,w}}(\xi)(\tilde{w})-\frac{\partial}{\partial w}D_{f^{v,w}}(\tilde{w})(\xi)+O(\tilde{w}^2)
\end{array}. \leqno(3.104)$$
Therefore, the second inequality follows from the first one.

Next, we consider the derivative of $D, N$ with respect to the $v$. First of all,
\vskip 0.1in
\noindent
{\bf Lemma 3.26: }{\it Let $|v-v_0|<\delta$ for small $\delta$ and $j_v$ be the complex structure
on $\Sigma_v$, there is a smooth family of
diffeomorphism $\Phi_v: \Sigma_{v_0}\rightarrow \Sigma_v$ such that $\Phi_v=id$ outside gluing
region and 
$$|\frac{\partial}{\partial v}|_{v=v_0}(\Phi_v j_v (\frac{\partial}{\partial t}))|\leq
 \frac{C}{r_0}.
\leqno(3.105)$$
$$|\frac{\partial}{\partial v}|_{v=v_0}(\Phi_v j_v (\frac{\partial}{\partial s}))|\leq 
\frac{C}{r_0}.
\leqno(3.106)$$}
\vskip 0.1in
{\bf Proof: } The complex structure outside the gluing region does not change. Over the gluing
region, it is conformal equivalent to a cylinder. Constructing
$\Phi_v$ in the cylindric model, we will obtain the estimate of Lemma 3.26.$\Box$

Suppose that we want to estimate the derivative at $v_0$. We fix $u=f^{v_0}$ and the trivialization
given by $\Phi_v$. To abuse the notation,
let
$f^{v,w}=exp_{f^{v_0}} w^v$. We still have the same Taylor expansion (3.96)-(3.100). Furthermore, we
can estimate $\frac{\partial}{\partial v}|_{v=v_0}\bigtriangledown^v \xi$ by the norms of 
$\bigtriangledown^{v_0} \xi$ and the derivative of $\Phi_v$. Hence,
\vskip 0.1in
\noindent
{\bf Corollary 3.27: }{\it Under the same condition of Lemma 3.26,
$$||\frac{\partial}{\partial v}|_{v=v_0} D_{f^{v,w}}(\xi)||_{L^p_v}\leq \frac{C}{|v_0|}
(||f^{v_0}||_{L^p_{1,v}}+||w^{v_0}||_{L^p_{1,v}})||\frac{\partial}{\partial v}|_{v=v_0} w^v||_{
L^p_{1,v}} ||\xi||_{L^p_{1,v}}.\leqno(3.107)$$
$$||\frac{\partial}{\partial v}|_{v=v_0}N_{f^{v,w}}(\xi)||_{L^p_v}\leq \frac{C}{|v_0|}
(||f^{v_0}||_{L^p_{1,v}}+||w^{v_0}||_{L^p_{1,v}})||\frac{\partial}{\partial v}|_{v=v_0} w^v||_{
L^p_{1,v}} ||\xi||^2_{L^p_{1,v}}.\leqno(3.108)$$}
\vskip 0.1in
Next, we compute the derivative of $Q_{f^{v,w}}$. Recall that
$$Q_{f^{v,w}}=AQ_{f^{v,w}}(D_{f^{v,w}} Q_{f^{v,w}})^{-1}.\leqno(3.109)$$
Therefore, it is enough to compute $AQ_{f^{v,w}}=\phi_f Q_{f^w}$ and $((D_{f^{v,w}} Q_{f^{v,w}})^{-1})'$.
Clearly,
$$\frac{\partial}{\partial w}AQ_{f^{w,v}}=\phi_f(\frac{\partial}{\partial w}Q_{f^w}).\leqno(3.110)$$
$$\frac{\partial}{\partial v}AQ_{f^{w,v}}=\frac{\partial }{\partial v}(\phi_f)Q_{f^w}.\leqno(3.111)$$
Recall that in the gluing construction, only the cut-off function has variable $v$. Hence, we 
need to compute the derivative of the cut-off function with respect to $v$.
\vskip 0.1in
\noindent
{\bf Lemma 3.28:}{\it $$|\frac{\partial}{\partial r} \bar{\beta}_r|<\frac{C}{r}.$$}
\vskip 0.1in
{\bf Proof: }
$$\bar{\beta}_r(t)=\int^{\frac{r^2}{2}}_{\frac{r^2\rho}{2}}(1-\frac{log(u)-log(r^2)}
{-log \rho})T(t-u)du+\int^{\infty}_{\frac{r^2}{2}} T(t-u)du.$$
$$\frac{\partial}{\partial r}\bar{\beta}_r(t)=\int \frac{\partial}{\partial r}
\frac{log(u)-log(r^2)}{-log \rho}T(t-u)du+(1-\frac{log(\frac{r^2}{2})-log(r^2)}
{-log \rho})T(t-\frac{r^2}{2})r$$
$$+(1-\frac{log(\frac{r^2\rho}{2})-log(r^2)}
{-log \rho})T(t-\frac{r^2\rho}{2})r\rho+ T(t-\frac{r^2}{2})r,\leqno(3.112)$$
where $T(t-u)$ is a positive smooth function with compact supported and integral 1. Then,
$$|\frac{\partial}{\partial r} \bar{\beta}_r|<\frac{C}{r}.\leqno(3.113)$$
$\Box$

Furthermore, we can choose $Q_{f^w}$ such that $\frac{\partial}{\partial w} Q_{f^w}$ is bounded.
Therefore,
$$||\frac{\partial}{\partial w}AQ_{f^{w,v}}||< C.\leqno(3.114)$$
$$||\frac{\partial}{\partial v}|_{v=v_0}AQ_{f^{w,v}}||< \frac{C}{|v_0|}.\leqno(3.115)$$
Notes that
$$D_{f^{w,v}}AQ_{f^{v,w}}(D_{f^{w,v}}AQ_{f^{v,w}})^{-1}=Id.\leqno(3.116)$$
Hence, 
$$((D_{f^{w,v}}AQ_{f^{v,w}})^{-1})'=-(D_{f^{w,v}}AQ_{f^{v,w}})^{-1}(D_{f^{w,v}}
AQ_{f^{v,w}})'(D_{f^{w,v}}AQ_{f^{v,w}})^{-1}.\leqno(3.117)$$
Combined (3.114)-(3.117), we obtain
\vskip 0.1in
\noindent
{\bf Lemma 3.29: }{\it $$||\frac{\partial}{\partial w} Q_{f^{v,w}}||\leq C.\leqno(3.118)$$
$$||\frac{\partial}{\partial v}|_{v=v_0} Q_{f^{v,w}}||\leq \frac{C}{|v_0|}.\leqno(3.119)$$}
Next, let's compute the derivative of $\S_e(f^{v,w})$. Let $w_{\mu}\in E^D_f$ be a smooth path such
that $w_0=w$ and $\frac{d}{d\mu}|_{\mu=0} w_{\mu}=\tilde{w}.$
\vskip 0.1in
\noindent
{\bf Lemma 3.30: }{\it For $w\in E^D_f$, we view $\S_e(f^{v,w})$ as a map from $E^D_f\times 
B_{\delta_f}(\C_f)$
to $\U_f$ where we use local trivialization given by $\Phi_v$ in Lemma 3.26. Then,
$$||\frac{d}{d\mu}|_{\mu=0} \S_e(f^{v,w_{\mu}})||_{L^p_v}\leq C r^{\frac{4}{p}},\leqno(3.120)$$
$$||\frac{\partial}{\partial v}|_{v=v_0} \S_e(f^{v,w})||_{L^p_{v_0}}\leq Cr_0^{\frac{4}{p}-1}.
\leqno(3.121)$$}
\vskip 0.1in
{\bf Proof: } $\S_e(f^{v,w_{\mu}})=0$ outside the gluing region and over $N_p(\frac{r^2}{2}, 2r^2)$. 
Therefore, the derivative is zero
outside the gluing region and over $N_p(\frac{r^2}{2}, 2r^2)$. Here, we work over a slightly larger domain 
$N_p(\frac{\rho r^2_0}{2.1}, \frac{(2.1)r^2_0}{\rho})$ so that we can vary $r$ in a fixed domain.

 It is enough to work
over $N_p(\frac{\rho r^2_0}{2.1}, \frac{r^2_0}{2})$, where  
$$f^{v,w_{\mu}}=\bar{\beta}_r(t)(f^{w_{\mu}}_1(s, t))-f^{w_{\mu}}_1(s,0))+
f^{w_{\mu}}_2(s+\theta,\frac{r^4}{t}).\leqno(3.124)$$
$$\S_e(f^{v,w_{\mu}})=\bigtriangledown \bar{\beta}_r(t)(f^{w_{\mu}}_1(s, t))-
f^{w_{\mu}}_1(s,0)).\leqno(3.125)$$
Therefore,
$$\frac{d}{d\mu}|_{\mu=0}\S_e(f^{v,w_{\mu}})=\bigtriangledown \bar{\beta}_r(t)(
\tilde{w}_1(s,t)-\tilde{w}_1(s,0)).\leqno(3.126)$$
$$||\frac{d}{d\mu}|_{\mu=0}\S_e(f^{v,w_{\mu}})||_{L^p_v}\leq C r^{\frac{4}{p}} 
 ||\tilde{w}||_{C^1}.\leqno(3.127)$$
Since $\tilde{w}$ varies in a finite dimensional space and $\tilde{w}$ is smooth, we can
replace $C^1$ norm by $L^p_1$-norm. Hence,
$$||\frac{d}{d\mu}|_{\mu=0}\S_e(f^{v,w_{\mu}})||_{L^p_v}\leq C r^{\frac{2}{p}} 
 ||\tilde{w}||_{L^p_1}.\leqno(3.128)$$ 
When we pull it back to the $\Sigma_{v_0}$ by $\Phi_v=(\Phi^1_v, \Phi^2_v)$,
$$\S_e(f^{v,w})=\bigtriangledown \bar{\beta}_r(\Phi^2_v(s,t))
(f^{w}(\Phi_v(t,s))-f^{w}_1(\Phi^1_v(t,s),0)).\leqno(3.129)$$
Using  Lemma 3.26 and Lemma 3.28, it is easy to estimate that
$$|\frac{\partial}{\partial v}|_{v=v_0}\S_e(f^{v,w})|\leq  C\frac{1}{r_0}||f^{w}||_{C^1}.\leqno(3.130)$$
Hence,
$$||\frac{\partial}{\partial v}|_{v=v_0}\S_e(f^{v,w})||_{L^p_{v_0}}\leq C\frac{1}{r_0}
vol(N_p(\frac{\rho r^2_0}{2}, \frac{r^2_0}{2}))^{\frac{1}{p}}||f^w||_{C^1}\leq C r_0^{\frac{4}{p}-1}.
\leqno(3.131)$$
Here, we use the  fact that $f^w$ is smooth and varies in a finite dimension set $E^D_f$ with
bounded $L^p_1$ norm. $\Box$.

   The same analysis will also implies that
\vskip 0.1in
\noindent
{\bf Lemma 3.31: }{\it $$||\frac{\partial}{\partial v}|_{v=v_0} f^{v,w}||_{
L^p_{1,v_0}}\leq C||w||_{L^p_1}\leqno(3.132)$$
for $f^w\in E^D_f$.}
\vskip 0.1in
We leave it to readers to fill out the detail. Let $F$ be the inverse of $exp_{f^{v_0}}$.
Then, 
$$||\frac{\partial}{\partial v}|_{v=v_0}w^v||\leq C(F)||\frac{
\partial}{\partial v}|_{v=v_0} f^{v,w}||_{L^p_{1,v}} \leq C(F)||w||_{L^p_1}.\leqno(3.133)$$

  Putting all the estimate together, we obtain
\vskip 0.1in
\noindent
{\bf Proposition 3.32: }{\it $$||\frac{d}{d\mu}|_{\mu=0} \xi(v, w_{\mu})||_{L^p_{1,v}}
\leq Cr^{\frac{4}{p}-1}.\leqno(3.134)$$
$$||\frac{\partial}{\partial v}|_{v=v_0} \xi(v, w)||_{L^p_{1,v}} \leq C r_0^{\frac{4}{p}-1}.
\leqno(3.135)$$}
\vskip 0.1in
{\bf Proof: } Recall that
$$\xi(v,w)=H(v,w,\xi(v,w))=-Q_{f^{v,w}}\S_e(f^{v,w})-Q_{f^{v,w}}N_{f^{v,w}}(\xi(v,w)).\leqno(3.136)$$
By Lemma 3.25-3.32, we have bound derivatives for all the term of $H$. Moreover, the derivative of error
term $\S_e(f^{v,w})$ is of the order $r^{\frac{4}{p}}$. Recall $\xi(v,w)$ is obtained by iterating $H$.
Hence, the derivative of $\xi(v, w)$ is bounded by $\delta$ in (3.91) when $r$ is
small. 
$$ \xi'(v,w)=Q'_{f^{v,w}}\S_e(f^{v, w})
-Q_{f^{v, w}}\S'_e(f^{v, w})-(Q'_{f^{v, w}}N_{f^{v,w}}+Q_{f^{v,w}}N'_{f^{v,w}})(\xi(v, w))
-Q_{f^{v, w}}N_{f^{v,w}}(\xi'(v,w)).\leqno(3.137)$$
By Lemma 3.25-3.32 and formula 3.133,
$$||\frac{d}{d\mu}|_{\mu=0} \xi(v,w_{\mu})||_{L^p_{1,v}}\leq C_1 r^{\frac{4}{p}-1}
+C_2 ||\xi(v,w)||_{L^p_{1,v}}+
C_3||\frac{d}{d\mu}|_{\mu=0} \xi(v,w_{\mu})||^2_{L^p_{1,v}}.\leqno(3.138)$$
$$||\frac{d}{d\mu}|_{\mu=0} \xi(v,w_{\mu})||_{L^p_{1,v}}\leq\frac{1}{1-\delta C_3}
(C_1 r^{\frac{4}{p}-1}+C_2 ||\xi(v,w)||_{L^p_{1,v}}).\leqno(3.139)$$
Using (3.93), we obtain the  inequality (3.134).
The proof of the second inequality (3.135) is completely same. Only difference is that the derivative of
$Q_{f^{w,v}}, N_{f^{v,w}}$ has a order $\frac{1}{r_0}$. However, we have $\S_e(f^{v,w}), \xi(v, w)$
in the formula, where both have order $r_0^{\frac{4}{p}}$. Hence, we obtain the order
$r_0^{\frac{4}{p}-1}$. $\Box$

  Let $u$ be a map over $\Sigma_v$ and $\xi\in \Omega^0(f^*T_FV)$. We define 
$u_v=(u^1_v,  u^2_v)$ and $\xi_v=(\xi^1_v, \xi^2_v)$ as in (3.36), (3.73). Now, we want to embed
$\M_{\S_e}\cap \U_f$ into $\U_{f,D}\times \C_f$ by the map
$$exp_u \xi\rightarrow (exp_{u_v} \xi_v,v)\leqno(3.140)$$ 
for $u$ over $\Sigma_v$. 
Consider the composition of (3.140) with $\alpha_{v,w}$.
$$\alpha(v,w)_v: E^D_f\times B_{\delta_f}(\C_f)\rightarrow L^p_1(\U_{f,D})\times \C_f.\leqno(3.141)$$
\vskip 0.1in
\noindent
{\bf Proposition 3.33: }{\it $\alpha(v,w)_v$ is $C^1$-smooth.}
\vskip 0.1in
{\bf Proof: } Our proof is motivated by the following observation. Suppose that $f$ is
a continuous function over $\R$ such that $f(0)=0$ and $f$ is $C^1$ for $x\neq 0$.
If $|f'(x)|\leq Cx^{\alpha}$ for $\alpha>0$, by mean value theorem $f'(0)=0$ and $f'$ is
continuous at $x=0$.

 We first prove 
$$(f^w,v)\rightarrow f^{v,w}_v \leqno(3.142)$$
is a $C^1$-map. $f^{v,w}_v=f^w$ outside the gluing region. By symmetry, it is enough to consider
$D_p(\frac{2r^2}{\rho})$. Over $D_p(2r^2)$,
$$\begin{array}{lll}
f^{v,w}_v&=&\tilde{\beta}_r(f^w_1(s,t)+f^w_2(s,t)-\frac{1}{2\pi r^2}
\int_{S^1}(f^{w}_1(s,r^2)+f^w_2(\theta+s, r^2)))\\
&&+\frac{1}{2\pi r^2}\int_{S^1}(f^{w}_1(s,r^2)+f^w_2(\theta+s, r^2))
\end{array}. \leqno(3.143)$$
$$\begin{array}{lll}
\frac{d}{d\mu}|_{\mu=0} (f^{v, w_{\mu}}_v-f^{w_{\mu}})&=&
\tilde{\beta}_r(\tilde{w}_1(s,t)+\tilde{w}_2(s,t)-\frac{1}{2\pi r^2}
\int_{S^1}(\tilde{w}_1(s,r^2)+\tilde{w}_2(\theta+s, r^2)))\\
&&+\frac{1}{2\pi r^2}\int_{S^1}(\tilde{w}_1(s,r^2)+
\tilde{w}_2(\theta+s, r^2))
\end{array}.\leqno(3.144)$$
Note that
$$|\frac{1}{2\pi r^2}\int_{S^1}\tilde{w}(s,r^2)- \tilde{w}(s,0)|
\leq C r^2 ||\tilde{w}||_{C^1}. \leqno(3.145).$$
By inserting the term $\tilde{w}(s,0)$ in the formula (3.144), 
$$||\frac{d}{d\mu}|_{\mu=0}( f^{v, w_{\mu}}-f^{w_{\mu}})||_{L^p(D_p(2r^2)}\leq C 
vol(D_p(2r^2))^{\frac{1}{p}} ||\tilde{w}||_{C^1}\leq C r^{\frac{4}{p}} 
||\tilde{w}||_{L^p_1}.\leqno(3.146)$$
Here, we use the fact that $\tilde{w}$ varies in a finite dimensional space. 
$$\begin{array}{lll}
&&||\bigtriangledown \frac{d}{d\mu}|_{\mu=0}( f^{v, w_{\mu}}-f^{w_{\mu}})||_{L^p(D_p(2r^2)}\\
&\leq & ||\bigtriangledown \tilde{\beta}_r
(\tilde{w}_1(s,t)+\tilde{w}_2(s,t)-\frac{1}{2\pi r^2}
\int_{S^1}(\tilde{w}_1(s,r^2)+\tilde{w}_2(\theta+s, r^2))||_{L^p}\\
&&+
||\tilde{\beta}_r\bigtriangledown (\tilde{w}_1+\tilde{w}_2)||_{L^p(D_p(2r^2)}\\
&\leq & C (vol(D_p(2r^2)))^{\frac{1}{p}} ||\tilde{w}||_{C^1}\\
&\leq &C r^{\frac{4}{p}}||\tilde{w}||_{L^p_1}
\end{array}. \leqno(3.147)$$
Over $N_p(\frac{\rho r^2}{2}, \frac{r^2}{2})$,
$$f^{v,w}_v=f^w_2(s,t)+\bar{\beta}_r(f^w_1(s+\theta, t)-f^w_1(0)).\leqno(3.147.1)$$
$$\frac{d}{d\mu}|_{\mu=0}(f^{v, w_{\mu}}_v-f^{w_{\mu}})=\bar{\beta}_r(\tilde{w}_2(s,t)-
\tilde{w}_2(0)).\leqno(3.147.2)$$
The same argument shows that
$$||\frac{d}{d\mu}|_{\mu=0}(f^{v, w_{\mu}}_v-f^{w_{\mu}})||_{L^p_1(N_p(\frac{\rho r^2}{2}, 
\frac{r^2}{2}))}\leq
Cr^{\frac{4}{p}}||\tilde{w}||_{L^p_1}.\leqno(3.147.3)$$
Using previous argument and Lemma 3.26, we can also show that
$$||\frac{\partial}{\partial v}|_{v=v_0} (f^{v,w}_v-f^w)||_{L^p_1}\leq Cr_0^{\frac{4}{p}-1}.
\leqno(3.148)$$
Therefore,
$$||(f^{v,w}_v)'-(f^w)'||\leq Cr^{\frac{4}{p}-1}.\leqno(3.149).$$
 $f^{v,w}_v$  is $C^1$ for $v\neq 0$. At $v=0$, the estimate (3.149) implies 
$$(f^{v,w}_v)'=(f^w)'  \mbox{ at } v=0.\leqno(3.150)$$
Moreover, $(f^{v,w}_v)'$ is continuous. The same argument together with Proposition 3.32 shows that
$$||(\xi(v,w)_v)'||\leq C r^{\frac{4}{p}-1}.\leqno(3.151)$$
Hence, $\xi(v,w)_v$ is a $C^1$-map and has derivative zero at $v=0$.
In general,
$$(exp_{f^{v,w}_v} \xi(v,w)_v)'=D_1exp_{f^{v,w}_v}\xi(v,w)_v (f^{v,w}_v)'
+D_2exp_{f^{v,w}_v}\xi(v,w)_v (\xi(v,w)_v)',\leqno(3.152)$$
where $D_1, D_2$ are the partial derivatives of $exp$-function.
$$\begin{array}{lll}
&&||(exp_{f^{v,w}_v} \xi(v,w)_v)'-(f^w)'||\\
&\leq &||D_1exp_{f^{v,w}_v}\xi(v,w)_v 
(f^{v,w}_v)'-(D_1exp_{f^{v,w}_v} 0)(f^{v,w}_v)'||\\
&&+||(f^{v,w}_v)'-(f^w)'||
+||D_2exp_{f^{v,w}_v}\xi(v,w)_v (\xi(v,w)_v)'||\\
&\leq&C||\xi(v,w)_v||_{L^{\infty}}||(f^{v,w}_v)'||+Cr^{\frac{4}{p}-1}\\
& \leq &C||\xi(v,w)_v||_{L^p_1}||(f^w)'||+Cr^{\frac{4}{p}-1}\\
&\leq &C r^{\frac{4}{p}-1}
\end{array}. \leqno(3.153)$$
Notes that $\alpha(v,w)_v$ is identity on $\C_f$-factor. Hence,  
we prove the proposition. Moreover, the derivative of $\alpha(v,w)_v$ is identity at $v=0$,
since $\frac{\partial}{\partial w} f^w=id$.
$\Box$
\vskip 0.1in
\noindent
{\bf Theorem 3.34: }{\it With the coordinate system given by $(E^D_f\times B_{\delta_f}(\C_f), 
\alpha_f(v,w))$, $\M_{\S_e}\cap \U_{\S_e}$ is a $C^1$-oriented V-manifold.}
\vskip 0.1in
{\bf Proof: } Recall the definition 2.1. Suppose that $\alpha_{\bar{f}}(E^{\bar{D}}_{\bar{f}}
\times B_{\delta_{\bar{f}}}(\C_{\bar{f}}))\subset \alpha_f(E^D_f\times B_{\delta_f}(\C_f))$.
Then, $stb_{\bar{f}}\subset stb_f$ and we can assume that $E^{\bar{D}}_{\bar{f}}\subset \U_{\bar{f}}
\subset \U_f$.
It is clear that $\bar{D}$ is either a higher strata than $D$ or $D$. Let's consider
the case that $\bar{D}$ is a higher strata. The proof for the second case is the same.
To be more precise, let's
consider the case that $D$ has three components $\Sigma_1\wedge \Sigma_2\wedge \Sigma_3$ and
$\bar{D}$ has two components $\Sigma_1\wedge \Sigma_2\#_{v_2} \Sigma_3$ for $v_2\neq 0$. The general case is similar and we 
leave it to readers. Suppose that the gluing parameters are $(v_1, v_2)\in \C_1\times \C_2$.
To construct Banach manifold $L^p_1(\U_{\bar{f}})$, we need a trivialization
of $\bigcup_{v_2}\Sigma_2\#_{v_2} \Sigma_3$. As we discuss in the beginning of this section, we can
choose any trivialization. Here, we choose the one given by $\Phi_{v_2}$ Lemma 3.26. 
Clearly, $\alpha_f(v,w)$ maps
an open subset of $E^D_f\times B_{\delta_f}(\C_2)$ onto $E^{\bar{D}}_{\bar{f}}$ as a 
diffeomorphism. Now, we embed
$\M_{\S_e}\cap \U_{\bar{f}}$ into $\B_{\bar{D}}$ by (3.140). By Proposition 3.33, both
$$\alpha_f(v,w)_{v_1}, \alpha_{\bar{f}}(v_1, w)_{v_1}\leqno(3.154)$$
are injective $C^1$-map. Hence, we can view the image of $\M_{\S_e}\cap \U_{\bar{f}}$
as a $C^1$-submanifold of $\B_{\bar{D}}\times \C_{\bar{f}}$ and both $\alpha_f(v,w)_{v_1}, \alpha_{\bar{f}}(v_1, w)_{v_1}$
as $C^1$-diffeomorphisms to this submanifold. Hence,
$$(\alpha_f(v,w))^{-1}\alpha_{\bar{f}}(v_1,w)=(\alpha_f(v,w)_{v_1})^{-1}
\alpha_{\bar{f}}(v_1, w)_{v_1}.\leqno(3.155)$$
is a $C^1$-diffeomorphism. 

 Next, we consider the orientation. First of all, it was proved in \cite{RT1} (Theorem 6.1)
that both $\alpha_f(v,w)$ and $\alpha_{\bar{f}}(v_1, w)$ are orientation
preserving diffeomorphism when $v_1\neq 0, v_2\neq 0$. 
Therefore, it is enough to consider the case $v_1=0$. By our argument in Proposition 3.33 (3.150,
3.151), 
$$(\alpha_f(v,w)_{v_1})'|_{v_1=0}=(\alpha_f(v_2, w))'|_{v_1=0}\times id_{C_1}, \
(\alpha_{\bar{f}}(v_1,w))'|_{v_1=0}=id. \leqno(3.157)$$
Moreover, $\alpha_f(v_2,w)$ is an orientation preserving diffeomorphism. Hence, the transition map
is an orientation preserving diffeomorphism. We finish the proof. $\Box$.

\section{GW-invariants of a family of symplectic manifolds}
  In this section, we shall give a detail construction of GW-invariants
for a family of symplectic manifolds. Furthermore, we will prove 
composition law and $k$-reduction formula. Let's recall the construction
in the introduction.

  Let 
$$p: Y\rightarrow M\leqno(4.1)$$ 
be an oriented fiber bundle such that the fiber $X$ and the base
$M$ are smooth, compact, oriented manifolds. Then, $Y$ is also a smooth, 
compact,
oriented manifold. Let $\omega$ be a closed 2-form on $Y$ such that $\omega$
restricts to a symplectic form over each fiber. Hence, we can view $Y$ as
a family of symplectic manifolds.  A $\omega$-tamed almost complex
structure $J$ is an automorphism of the vertical tangent bundle $V(Y)$ such that 
$J^2=-Id$
and $\omega(w, Jw)>0$ for any vertical tangent vector $w\neq 0$. Suppose $A\in H_2(
V, \Z)\subset H_2(Y, \Z)$. Let $\M_{g,k}$ be the moduli space of 
genus g Riemann surfaces with $k$-marked points such that $2g+k>2$ and 
$\overline{\M}_{g,k}$ be its Deligne-Mumford compactification. We shall use
$$f: \Sigma\stackrel{F}{\rightarrow} Y$$
to indicate that the $im(f)$ is contained in a fiber. Consider its compactification-
the moduli space of stable holomorphic maps $\overline{\M}_A(Y,g,k,J)$.

 Using the machinery of section 2 and 3,  we can define a virtual neighborhood 
invariant $\mu_{\S}$. Here, we have to specify the cohomology class $\alpha$ 
in the definition of virtual neighborhood invariant $\mu_{\S}$. Recall that we 
have two natural maps
$$\Xi_{g,k}: \overline{\B}_A(Y,g,k)\rightarrow Y^k\leqno(4.2)$$
defined by evaluating $f$ at marked points
and 
$$\chi_{g,k}: \overline{\B}_A(Y,g,k)\rightarrow \overline{\M}_{g,k}\leqno(4.3)$$
defined by forgetting the map and contracting the unstable components of 
the domain. 
Notes that $\overline{\M}_{g,k}$ is a V-manifold. Suppose $K\in H^*(\overline{
\M}_{g,k}, \R)$ and $\alpha_i\in H^*(V, \R)$ are represented by differential forms.
\vskip 0.1in
\noindent
{\bf Definition 4.1: }{\it We define
$$\Psi^Y_{(A,g,k)}(K; \alpha_1, \cdots, \alpha_k)=\mu_{\S}(\chi^*_{g,k}(K)
\wedge\Xi_{g,k}^*(\prod_i \alpha_i)).\leqno(4.4)$$}
\vskip 0.1in
\noindent
{\bf Theorem 4.2 }{\it (i).$\Psi^Y_{(A,g,k)}(K; \alpha_1, \cdots, \alpha_k)$ is
well-defined, multi-linear and skew symmetry.
\vskip 0.1in
\noindent
(ii). $\Psi^Y_{(A,g,k)}(K; \alpha_1, \cdots, \alpha_k)$ is independent of the
choice of forms $K, \alpha_i$ representing the cohomology classes $[K], [
\alpha_i]$,  and the choice of virtual neighborhoods.
\vskip 0.1in
\noindent
(iii). $\Psi^Y_{(A,g,k)}(K; \alpha_1, \cdots, \alpha_k)$ is independent of $J$
and is a symplectic deformation invariant.
\vskip 0.1in
\noindent
(iv). When $Y=V$ is semi-positive and some multiple of $[K]$ is represented by
an immersed V-submanifold, $\Psi^Y_{(A,g,k)}(K; \alpha_1, \cdots, 
\alpha_k)$ agrees with the definition of \cite{RT2}.}
\vskip 0.1in
{\bf Proof: } (i) follows from the definition and we omit it. (ii) follows from
 Proposition 2.7.

To prove (iii), suppose that $\omega_t$ is a family of symplectic structures
and $J_t$ is a family of almost complex structures such that $J_t$ is tamed
with $\omega_t$. Then, we can construct a weakly smooth Banach  
cobordism
$(\B_{(t)}, \F_{(t)}, \S_{(t)})$ of 
$$\overline{\M}_{A}(Y, g,k, J_{(t)})=\cup_{t\in [0,1]}\overline{\M}_A(Y,g,k,
J_t)\times \{t\}. \leqno(4.5)$$
Then, (iii) follows from Proposition 2.8 and section 3.

  To prove (iv), recall the construction of \cite{RT2}. To avoid the confusion,
we will use $\Phi$ to denote the invariant defined in \cite{RT2}. The 
construction of \cite{RT1} starts from an inhomogeneous Cauchy-Riemann equation.
It was  known that $\overline{\M}_{g,k}$ does not admit a universal family, which
causes a problem to define inhomogeneous term. To overcome this difficulty, Tian
 and the author choose a finite cover 
$$p_{\mu}: \overline{\M}^{\mu}_{g,k}\rightarrow \overline{\M}_{g,k}.\leqno(4.6)$$
such that $\overline{\M}^{\mu}_{g,k}$ admits a universal family. One can use
the universal family of $\overline{\M}^{\mu}_{g,k}$ to define an inhomogeneous
term $\nu$ and inhomogeneous Cauchy-Riemann
equation $\bar{\partial}_J f=\nu$. Any $f$ satisfying this equation is called {\em a
$(J,\nu)$-map}.  Choose a generic $(J, \nu)$ such that
the moduli space $\M^{\mu}_A(\mu, g,k,J, \nu)$ of $(J, \nu)$-map is smooth and the certain 
contraction  $\overline{\M}^r_A(\mu, g,k,J,\nu)$ of $\overline{\M}_A(\mu, g,k,J,\nu)$ is of codimension 2 
boundary. Define
$$\Xi^{\mu, \nu}_{g,k}: \overline{\M}_A(\mu, g,k,J,\nu)\rightarrow X^k$$
and
$$\chi^{\mu, \nu}_{g,k}: \overline{\M}_A(\mu,g,k,J,\nu)\rightarrow \overline{\M
}^{\mu}_{g,k}$$
similarly. Then, we can choose Poincare duals (as pseudo-submanifolds)
$K^*, \alpha^*$ of $K, \alpha_i$ such that $K^*, \alpha^*$ did not meet
the image of $\overline{\M}_A(\mu,g,k,J,\nu)-\M_A(\mu,g,k,J,\nu)$ under
the map $\chi^{\mu,\nu}_{g,k}\times \Xi^{\mu,\nu}_{g,k}$  and intersects
transversely to the  restriction  of  $\chi^{\mu,\nu}_{g,k}\times \Xi^{\mu,
\nu}_{g,k}$
to $\M_A(\mu, g,k,J,\nu)$. Once this is done,
$\Phi^X_{(A,g,k, \mu)}$ is defined as the number of the points of
$(\chi^{\mu,\nu}_{g,k}\times \Xi^{\mu,\nu}_{g,k})^{-1}(K^*\times \prod_i \alpha^*_i)$, counted
by the orientation. Then, we define 
$$\Phi^V_{(A,g,k)}(K; \alpha_1, \cdots, \alpha_k)=\frac{1}{\lambda^{\mu}_{g,k}}
\Phi^V_{(A,g,k, \mu)}(p_{\mu}^*(K); \alpha_1, \cdots, \alpha_k),$$
 where $\lambda^{\mu}_{g,k}$ is the order of
cover map $p_{\mu}$ (4.6).

   The proof of (iv) is divided into 3-steps. First we observe that we can replace $\overline{\M}_{g,k}$ by
$\overline{\M}^{\mu}_{g,k}$ in our construction. Let $\pi_{\mu}: \overline{\B}^{\mu}_{g,k}$
be the projection and $(\E_{g,k}, s_{g,k})$ be the stablization terms for $\overline{\M}_{g,k}$.
Then, we can choose $(\pi^*_{\mu}\E_{g,k}, \pi^*_{\mu}s_{g,k})$ to be the stablization term
of $\overline{\M}^{\mu}_{g,k}$. Suppose that the 
resulting finite dimensional virtual neighborhoods are $(U, E, S), (U^{\mu}, E^{\mu}, S^{\mu})$ and invariant are $\Psi^Y_{(A,g,k)}, \Psi^Y_{(A,g,k,\mu)}$, respectively. 
Then, we have a commutative diagram
$$\begin{array}{ccc}
U^{\mu}&\rightarrow &E^{\mu} \\
\downarrow&         & \downarrow \\
U     &\rightarrow  & E
\end{array}\leqno(4.7)$$
and
$$\begin{array}{ccc}
U^{\mu} &  \rightarrow & V^k\times \overline{\M}^{\mu}_{g,k}\\
\downarrow &           &\downarrow \\
U      & \rightarrow  & V^k\times \overline{\M}_{g,k}
\end{array}.\leqno(4.8)$$
Let $\lambda$ be the order of the cover $p_U: U^{\mu}\rightarrow U$ and
$\lambda'$ be the order of the cover $p_G: E^{\mu}\rightarrow E$.
One can check that 
$$\lambda=\lambda' \lambda^{\mu}_{g,k}.\leqno(4.9)$$
Let $\Theta$ be a Thom-form supported in a neighborhood of zero section of $E$. Then,
$$\begin{array}{lll}
\Psi^V_{(A,g,k)}(K; \alpha_1, \cdots, \alpha_k)&=&\int_{U}\chi_{g,k}^*(K)\wedge\Xi_{g,k}^*(\prod_i 
   \alpha_i)\wedge  S^*(\Theta)\\
 &=&\frac{1}{\lambda}\int_{U^{\mu}}(\chi^{\mu}_{g,k})^*(p_{\mu}^*(K))\wedge 
(\Xi^{\mu}_{g,k})^*(\prod_i \alpha_i)\wedge (p_U S)^*(\Theta)\\
 &=&\frac{1}{\lambda^{\mu}_{g,k}}
\int_{U^{\mu}}(\chi^{\mu}_{g,
k})^*(p^*_{\mu}(K))\wedge (\Xi^{\mu}_{g,k})^*(\prod_i \alpha_i)\wedge (S^{\mu})^*(\frac{1}{\lambda'}p_G^*(\Theta))\\
 &=&\frac{1}{\lambda^{\mu}_{g,k}}\Psi^V_{(A,g,k,\mu)}(p_{\mu}^*(K); \alpha_1,
                \cdots, \alpha_k)
\end{array}\leqno(4.10)$$
where $\frac{1}{\lambda'} p_G^*(\Theta)$ is a Thom form of $E^{\mu}$.
Therefore, it is enough to show that 
$$\Psi^V_{(A,g,k,\mu)}=\Phi^V_{(A,g,k,\mu)}. $$

The second step is to deform Cauchy-Riemann equation $\bar{\partial}_Jf=0$ to 
inhomogeneous equation $\bar{\partial}_Jf=\nu$. Consider a family of equations
$\bar{\partial}_J f=t\nu$. We can repeat the argument of (ii) to show
that $\Psi^Y_{(A,g,k,\mu)}$ is independent of $t$. 

Let $(\B^{\mu,\nu}_{g,k}, \F^{\mu, \nu}_{g,k}, \S^{\mu, \nu}_{g,k})$ be VNA smooth compact V-triple  of 
$\overline{\M}^{\mu}_A(g,k,J,\nu)$ and define $\Xi^{\mu, \nu}_{g,k}, \chi^{\mu,
\nu}_{g,k}$ similarly. For the same reason, the virtual neighborhood 
construction applies. 
The third step is to construct a particular finite dimensional virtual  
neighborhood $(U^{\mu}_{\nu}, E^{\mu}_{\nu}, S^{\mu}_{\nu})$ such that
the restriction 
$$\chi^{\mu,\nu}_{g,k}\times \Xi^{\mu,\nu}_{g,k}: U^{\mu}_{\nu}\rightarrow X^k\times
\overline{\M}^{\mu}_{g,k}\leqno(4.11)$$
is transverse to $K^*\times \prod_i \alpha^*_i$.

 First of all,
since we work over $\R$, we can assume that each $\alpha^*$ is represented by
a bordism class, and hence an immersed submanifold by ordinary transversality.
By the linearity (i), we can assume that  $K^*$ is represented by an immersed V-submanifold.
Hence, $K^*\times \prod_i \alpha^*$ is represented by an immersed -submanifold
(still denoted by $K^*\times \prod_i \alpha^*$). We first assume that
$K^*\times \prod_i \alpha^*$ is an embedded V-submanifold. Recall that
$K^*\times \prod_i\alpha^*$ does not meet the image of 
$\overline{\M}_A(\mu,g,k,J,\nu)-\M_A(\mu,g,k,J,\nu)$ and intersects 
transversely to the image $\M_A(\mu, g,k,J,\nu)$. Therefore,
$$(\chi^{\mu,\nu}_{g,k}\times \Xi^{\mu,\nu}_{g,k})^{-1}(K^*\times \prod_i\alpha^*)\cap
\overline{\M}^{\mu}_A(g,k,J,\nu)\leqno(4.12)$$ 
is a collection of the smooth points of $\M_A(g,k,J,\nu)$. It implies that $L_x$
is surjective at $x\in (\chi^{\mu,\nu}_{g,k}\times \Xi^{\mu,\nu}_{g,k})^{-1}(
K^*\times 
\prod_i\alpha^*)\cap \overline{\M}^{\mu}_A(g,k,J,\nu)$ and 
$$\delta(\chi^{\mu,\nu}_{g,k}\times \Xi^{\mu}_{g,k}): Ker L_A\rightarrow X^k\times
\overline{\M}^{\mu}_{g,k}\leqno(4.13)$$
is surjective onto the normal bundle of $K^*\times \prod_i \alpha_i$.
Hence, the same is true over an open neighborhood $\U'$ of $(\chi^{\mu,\nu}_{g,k}\times \Xi^{\mu,\nu}_{
g,k})^{-1}(K^*\times \prod_i\alpha^*)\cap
\overline{\M}_A(\mu,g,k,J,\nu)$.
We cover $\overline{\M}_A(\mu,g,k,J,\nu)$ by $\U'$ and $\U''$ such that
$$\bar{\U}''\cap (\chi^{\mu,\nu}_{g,k}\times \Xi^{\mu}_{g,k})^{-1}(K^*\times \prod_i
\alpha^*_i)=\emptyset. \leqno(4.14)$$
Then, we construct $(\E^{\mu}_{\nu}, s^{\mu}_{\nu})$ such that $s^{\mu}_{
\nu}=0$ over $\U'-\bar{\U}''$. 
Suppose that $(U^{\mu}_{\nu}, E^{\mu}_{\nu}, S^{\mu}_{\nu})$ is the
finite dimensional virtual neighborhood constructed by $(\E^{\mu}_{\nu},
s^{\mu}_{\nu})$. It is easy to check that 
$$(\chi^{\mu,\nu}_{g,k}\times \Xi^{\mu,\nu}_{g,k})^{-1}(K^*\times \prod_i\alpha^*_i)\cap
U^{\mu}_{\nu}\subset \U'-\bar{\U}''.\leqno(4.15)$$
On the other hand,
$$s^{\mu}_{\nu}=0 \mbox{ over } \U'-\bar{\U}''.\leqno(4.16)$$
It implies that
$$(\chi^{\mu,\nu}_{g,k}\times \Xi^{\mu,\nu}_{g,k})^{-1}(K^*\times \prod_i\alpha^*)\cap
U^{\mu}_{\nu}=E^{\mu}_{\nu}|_{((\chi^{\mu,\nu}_{g,k}\times \Xi^{\mu,\nu}_{g,k})^{-1}(K^*\times \prod_i
\alpha^*_i)\cap \overline{\M}_A(\mu,g,k,J,\nu))}.\leqno(4.17)$$
It is easy to observe that the restriction of $\chi^{\mu,\nu}_{g,k}\times 
\Xi^{\mu,\nu}_{g,k}$ to $U^{\mu}_{\nu}$ is transverse to $K^*\times \prod_i\alpha^*_i$.

Since $K^*\times \prod_i\alpha^*_i$ is Poincare dual to 
$K\times \prod_i\alpha_i$, $(\chi^{\mu,\nu}_{g,k}\times \Xi^{\mu,\nu}_{g,k}
)^{-1}(K^*\times \prod_i\alpha^*)$ is Poincare dual to $(\chi^{\mu,\nu}_{g,k}
\times \Xi^{\mu,\nu}_{g,k})^*(K\times \prod_i \alpha_i)$. Therefore,
$$\begin{array}{lll}
\Psi^V_{(A,g,k,\mu)}(K; \alpha_1, \cdots, \alpha_k)&=&\int_{U^{\mu}_{\nu}}(
\Xi^{\mu,\nu}_{g,k}\times \chi^{\mu,\nu}_{g,k})^*(K\times\prod_i
\alpha_i)\wedge (S^{\mu}_{\nu})^*(\Theta)\\
&=&\int_{(\Xi^{\mu,\nu}_{g,k}\times \chi^{\mu,\nu}_{g,k})^{-1}(K^*\times \prod_i\alpha^*
)\cap U^{\mu}_{\nu}}(S^{\mu}_{\nu})^*(\Theta)\\
&=&\Phi^V_{(A,g,k,\mu)}(K; \alpha_1, \cdots, \alpha_k)
\end{array}.$$
When $K^*\times \prod_i \alpha^*_i$ is an immersed V-submanifold, there is a
V-manifold $N$ and a smooth map
$$H:N\rightarrow X^k\times \overline{\M}^{\mu}_{g,k}$$
whose image is $K^*\times \prod_i \alpha^*_i$. Then, we replace $\chi^{\mu,
\nu}_{g,k}\times \Xi^{\mu,\nu}_{g,k}$ by $\chi^{\mu,\nu}_{g,k}\times
\Xi^{\mu,\nu}_{g,k}\times N$ and $K^*\times\prod_i \alpha^*_i$ by the
diagonal of $(X^k\times \overline{\M}^{\mu}_{g,k})^2$ in the previous argument.
It will implies (iv). $\Box$

It is well-known that the projection map $p: Y\rightarrow X$ defines a modular
structure on $H^*(Y, \R)$ by $H^*(M, \R)$, defined by
$$\alpha \cdot \beta=p^*(\alpha)\wedge \beta\leqno(4.18)$$
where $\alpha\in H^*(M, \R)$ and $\beta\in H^*(Y, \R)$. GW-invariant we defined
behave nicely over this modulo structure, which is the basis of the modulo
structure of equivariant quantum cohomology (Theorem I).
\vskip 0.1in
\noindent
{\bf Proposition 4.3: }{\it Suppose that $\alpha_i\in H^*(Y, \R), \alpha\in H^*(M,
\R)$. Then
$$\Psi^Y_{(A,g,k)}(K; \alpha_1, \cdots, \alpha\cdot \alpha_i,\cdots, \alpha_j,
\cdots, \alpha_k)$$
$$=\Psi^Y_{(A,g,k)}(K; \alpha_1, \cdots,  \alpha_i,\cdots, \alpha\cdot\alpha_j,
\cdots, \alpha_k). \leqno(4.19)$$}
\vskip 0.1in
{\bf Proof: } By the definition,
$$\Psi^Y_{(A,g,k)}(K; \alpha_1, \cdots,  \alpha_i,\cdots, \alpha_j,
\cdots, \alpha_k)$$
$$=\int_{U}\chi^*_{g,k}(K)\wedge \Xi^*_{g,k}(\prod_i \alpha_i)\wedge S^*(
\Theta).$$
Let 
$$p: Y^k\rightarrow V^k$$
and $\Delta$ be the diagonal of $V^k$.
A crucial observation is that 
$$\Xi^*_{g,k}: \overline{\B}_A(Y,g,k)\rightarrow Y^k$$
is factored through
$$\B_{g,k}\stackrel{\Xi'_{g,k}}{\rightarrow} p^{-1}(\Delta) \stackrel{i_{p^{-1}
(\Delta)}}{\rightarrow} Y^k.\leqno(4.20)$$
Furthermore, for any $i$
$$\begin{array}{lll}
&&i^*_{p^{-1}(\Delta)}(\alpha_1\times \cdots \times\alpha\cdot\alpha_i\times \cdots \alpha_k)\\
&=&p^*(i^*_{\Delta}(1\times\cdots\times \alpha^{(i)}\times\cdots\times 1 ))
  \wedge i^*_{p^{-1}(\Delta)}(\alpha_1\times \cdots \times\alpha_i\times \cdots
 \alpha_k)
\end{array}\leqno(4.21)$$
where we use $\alpha^{(i)}$ to indicate that $\alpha$ is at the $i$-th 
component. However,
$$i^*_{\Delta}(1\times\cdots\times \alpha^{(i)}\times\cdots\times 1 )=\alpha=
i^*_{\Delta}(1\times\cdots\times \alpha^{(j)}\times\cdots\times 1 ).\leqno(4.22)$$
Hence, 
$$\Xi^*_{g,k}(\alpha_1\times \cdots \times\alpha\cdot\alpha_i\times \cdots 
\alpha_k)=\Xi^*_{g,k}(\alpha_1\times \cdots \times\alpha\cdot\alpha_j\times 
\cdots \alpha_k).\leqno(4.23)$$
Then,
$$\Psi^Y_{(A,g,k)}(K; \alpha_1, \cdots, \alpha\cdot \alpha_i,\cdots, \alpha_j,
\cdots, \alpha_k)$$
$$=\Psi^Y_{(A,g,k)}(K; \alpha_1, \cdots,  \alpha_i,\cdots, \alpha\cdot\alpha_j,
\cdots, \alpha_k).$$
$\Box$

As we mentioned in the introduction, there is a natural map 
$$\pi: \overline{\M}_{g,k}\rightarrow \overline{\M}_{g, k-1} \leqno(4.24)$$
by forgetting the last
marked point and contracting the unstable rational component. 
 One should be aware that
there are two exceptional cases $(g,k)=(0,3), (1,1)$ where $\pi$ is not well
defined. $\pi$ is not a fiber bundle, but a Lefschetz fibration. However, the
integration over the fiber still holds for $\pi$. In another words, we have a
map
$$\pi_*: H^*(\overline{\M}_{g,k}, \R)\rightarrow H^{*-2}(\overline{\M}_{g,k-1},
\R).\leqno(4.25)$$
 For a stable $J$-map $f\in \overline{\M}_A(Y,g,k,J)$, let's also forget the
last marked point $x_k$. If the resulting map is unstable, the unstable 
component is either a constant or non-constant map. If it is a constant map, we
simply contract this component. If it is non-constant map, we divided it by the
larger automorphism group. Then, we obtain a stable $J$-map in $\overline{\M}_A(Y,
g,k-1,J)$. Furthermore, we have a commutative diagram
$$\begin{array}{ccc}
\chi_{g,k}: \overline{\M}_A(Y,g,k,J)&\rightarrow &\overline{\M}_{g,k}\\
           \downarrow \pi          &             &\downarrow \pi\\
\chi_{g,k-1}: \overline{\M}_A(Y,g,k-1, J)&\rightarrow &\overline{\M}_{g,k-1}
\end{array}\leqno(4.26)$$
Associated with $\pi$,
we have two {\em k-reduction formulas} for $\Psi^Y_{(A, g, k)}$.
\vskip 0.1in
\noindent
{\bf Proposition 4.4. }{\it Suppose that $(g,k)\neq (0,3),(1,1)$. 
\vskip 0.1in
\noindent
(1) For any $\alpha _1, \cdots , \alpha _{k-1}$ in $H^*(Y, \R)$, 
we have} 
$$\Psi ^Y_{(A,g,k)}(K; \alpha _1, \cdots,\alpha _{k-1}, 1)~=~
\Psi ^Y_{(A,g,k-1)}(\pi_*(K); \alpha _1, \cdots,\alpha _{k-1})
\leqno(4.27)$$
\vskip 0.1in
\noindent
(2)  Let $\alpha _k$ be in $H^2(Y, \R)$, then
$$\Psi ^Y_{(A,g,k)}(\pi^*(K); \alpha _1, 
\cdots,\alpha _{k-1}, \alpha _k)~=~\alpha_k (A)
\Psi ^Y_{(A,g,k-1)}(K; \alpha _1, \cdots,\alpha _{k-1})
\leqno (4.28)$$
where $ \alpha^* _k$ is the Poincare dual of $\alpha _k$.
\vskip 0.1in
\noindent
{\bf Proof:} Let $(\overline{\B}_A(Y,g,k), \overline{\F}_A(Y,g,k), \S^A_{g,k})$
 be the VNA smooth Banach compact V-triple of
$\overline{\M}_A(Y,g,k,J)$. Following from our construction of last section, we
 have commutative diagram
$$\begin{array}{ccc}
\chi_{g,k}: \overline{\B}_A(Y,g,k)&\rightarrow &\overline{\M}_{g,k}\\
           \downarrow \pi          &             &\downarrow \pi\\
\chi_{g,k-1}:\overline{\B}_A(Y,g,k)&\rightarrow &\overline{\M}_{g,k-1}
\end{array}\leqno(4.29)$$
 Furthermore, 
$\overline{\F}_A(Y,g,k)=\pi^* \overline{\F}_A(Y,g,k-1)$. Using the virtual 
neighborhood technique, we 
construct $(\E, s)$ and a finite dimensional virtual neighborhood 
$(U_{g,k-1}, E_{g,k-1}, S_{g,k-1})$ of $\overline{\M}_A(Y,g,k-1,J)$. We observe that
the same $(\E, s)$ also works in the construction of finite dimensional
virtual neighborhood of $\overline{\M}_A(Y,g,k,J)$. Let $(U_{g,k}, E_{g,k}, 
S_{g, k})$ be the virtual neighborhood. Then, $E_{g,k}$ is the pull  back of $E_{g,k-1}$ by
$\pi: \overline{\B}_{g,k} \rightarrow \overline{\B}_{g,k-1}$. There is a projection
$$\pi: U_{g,k}\rightarrow U_{g,k-1}. \leqno(4.30)$$
Then,
$$S_{g,k}=S_{g,k-1}\circ \pi.\leqno(4.31)$$
Hence,
$$S^*_{g,k}(\Theta)=\pi^* S^*_{g,k-1}(\Theta).\leqno(4.32)$$
Moreover,
$$\Xi^*_{g,k}(\prod^{k-1}_1 \alpha_i\wedge 1)=(\Xi_{g,k-1}\pi)^*(\prod^{k-1}_1
\alpha_i)=\pi^*\Xi^*_{g,k-1}(\prod^{k-1}_1\alpha_i).\leqno(4.33)$$
Furthermore,
$$\pi_*\chi_{g,k}^*(K)=\chi^*_{g,k-1}(\pi_*(K)).\leqno(4.34)$$
So,
$$\begin{array}{lll}
\Psi^Y_{(A,g,k)}(K; \alpha_1, \cdots,\alpha_{k-1}, 1)&=&\int_{U_{g,k}} \chi_{g,
 k}^*(K)\wedge  \Xi^*_{g,k}(\prod^{k-1}_1 \alpha_i \wedge 1)\wedge S^*_{g,k}(
  \Theta)\\
&=&\int_{U_{g,k-1}}\pi_*( \chi_{g, k}^*(K)\wedge \pi^*(\Xi^*_{g,k}(\prod^{k-1}_1\alpha_i
)\wedge S^*_{g,k-1}(  \Theta)))\\
&=&\int_{U_{g,k-1}}\chi^*_{g,k-1}(\pi_*K)\wedge \Xi^*_{g, k-1}(\prod^{k-1}_1\alpha_i)
   \wedge S_*(\Theta)\\
&=&\Psi^Y_{(A,g,k)}(\pi_*(K); \alpha_1, \cdots, \alpha_{k-1})
\end{array}.\leqno(4.35)$$
On the other hand, for $\alpha_k\in H^2(Y, \R)$,
$$\Xi^*_{g,k}(\prod^{k-1}_1\alpha_i\wedge \alpha_k)=\pi^*\Xi^*_{g,k-1}(\prod^{
k-1}_1\alpha_i)\wedge e^*_k (\alpha_k),\leqno(4.36)$$
where 
$$e_i: \overline{\B}_A(Y,g,k)\rightarrow  Y \leqno(4.37)$$
is the evaluation map at the marked point $x_k$.
One can check that
$$\pi_*(e^*_k(\alpha_k))=\alpha_k(A).\leqno(4.38)$$
Therefore,
$$\begin{array}{lll}
\Psi^Y_{(A,g,k)}(\pi^*(K); \alpha_1, \cdots, \alpha_{k-1}, \alpha_k)&=&
\int_{U_{g,k}} \chi_{g,
 k}^*(\pi^*(K))\wedge  \Xi^*_{g,k}(\prod^{k-1}_1 \alpha_i \wedge \alpha_k)\wedge S^*_{g,k}(
  \Theta)\\
&=&\int_{U_{g,k-1}}\pi_*( \chi_{g,
 k}^*(\pi^*(K))\wedge  \Xi^*_{g,k}(\prod^{k}_1\alpha_i\wedge \alpha_k)\wedge S^*_{g,k}(
  \Theta))\\
&=&\int_{U_{g,k-1}}\chi^*_{g,k-1}(K)\wedge \Xi^*_{g,k-1}(\prod^{k-1}_1\alpha_i)
\wedge S^*_{g,k-1}(\Theta)\wedge \pi_*(e^*_k(\alpha_k))\\
&=&\alpha_k(A)\Psi^Y_{(A,g,k-1)}(K; \alpha_1, \cdots, \alpha_{k-1})
\end{array}.\leqno(4.39)$$
$\Box$

Let $\overline \U_{g,k}$ be the universal  curve
over $\overline \M_{g,k}$. Then each marked point $x_i$ gives rise to
a section, still denoted by $x_i$, of the fibration
$\overline \U_{g,k} \mapsto 
\overline \M_{g,k}$. If ${\cal K}_{\U |\M}$ denotes the cotangent
bundle to fibers of this fibration, we put ${\cal L}_{(i)} = x_i^* (
{\cal K}_{\U |\M})$. Following Witten, we put
$$\langle \tau _{d_1,\alpha _1}\tau_{d_2,\alpha _2}\cdots \tau _{d_k,\alpha _k}
\rangle _g (q) ~=~\sum _{A \in H_2(X,\Z)} \Psi^X_{(A,g,k)}(K_{d_1,\cdots,d_k};
\{\alpha _i\}) \, q^A \leqno (4.40)$$
where $\alpha _i \in H_*(V,\Q)$ and $[K_{d_1,\cdots,d_k}]=c_1({\cal L}_{(1)})^{d_1} \cup \cdots \cup
c_1({\cal L}_{(k)})^{d_k}$ and $q$ is an element of Novikov ring.
Symbolically, $\tau _{d,\alpha}$'s denote ``quantum field theory operators''.
For simplicity, we only consider the cohomology classes
of even degree. Choose a 
basis $\{\beta_a\}_{1\le a\le N}$ of
$H^{*,\rm even}(V, \Z)$ modulo torsion. We introduce formal variables $t_r^a$,
where $r= 0, 1, 2, \cdots$ and $1\le a \le N$. Witten's generating
function (cf. [W2]) is now simply defined to be
$$F^X(t^a_r ; q) = 
\langle e^{\sum _{r,a} t^a_r \tau _{r, \beta _a}} \rangle (q)\lambda^{2g-2}
=\sum _{n_{r,a}} \prod _{r,a} {(t^a_r)^{n_{r,a}} \over {n_{r,a}}!}
\left \langle \prod _{r,a} \tau _{r,\beta _a}^{n_{r,a}} 
\right \rangle (q)\lambda^{2g-2}
\leqno (4.41)$$
where $n_{r,a}$ are arbitrary collections of nonnegative integers,
almost all zero, labeled by $r, a$. The summation in (4.40) is over all
values of the genus $g$ and all homotopy classes $A$ of $J$-maps.
Sometimes, we write $F_g^X$ to be the part of $F^X$ involving only
GW-invariants of genus $g$. Using the argument of Lemma 6.1 (\cite{RT2}), 
Proposition 4.4 implies that the generating function satisfies several equation. 
\vskip 0.1in
\noindent
{\bf Corollary 4.5. }{\it Let $X$ be a symplectic manifold. $F^X$ satisfies the generalized string
equation
$${\partial F^X\over \partial t^1_0} = {1\over 2} \eta _{ab} t^a_0 t^b_0
+ \sum \limits_{i=0}^\infty \sum \limits _{a} t^a_{i+1} 
{\partial F^X\over \partial t^a_i}.\leqno(4.42)$$
$F^X_g$ satisfies the dilaton equation
$$\frac{\partial F^X_g}{\partial t^1_1}=(2g-2+\sum_{i=1}^{\infty}\sum_{a}t^a_i
\frac{\partial }{\partial t^a_i})F^X_g+\frac{\chi(X)}{24}\delta_{g,1}, 
\leqno(4.43)
$$
where $\chi(X)$ is the Euler characteristic of $X$.}
\vskip 0.1in

Next, we prove the composition law. Recall the construction in the introduction.
Assume $g=g_1+g_2$ and $k=k_1+k_2$ with $2g_i + k_i \ge 3$.
Fix  a decomposition $S=S_1\cup S_2$ of $\{1,\cdots , k\}$ with
$|S_i|= k_i$. Recall that 
$\theta _S: \overline \M_{g_1,k_1+1}\times \overline \M_{g_2,k_2+1}
\mapsto \overline \M_{g,k}$, which assigns to marked
curves $(\Sigma _i; x_1^i,\cdots ,x_{k_1+1}^i)$ ($i=1,2$), their 
union $\Sigma _1\cup \Sigma _2$ with $x^1_{k_1+1}$ identified to
$x^2_1$ and remaining points renumbered by $\{1,\cdots,k\}$ according to $S$.
Clearly, $im(\theta_S)$ is a V-submanifold of $\overline{\M}_{g,k}$, where the
Poincare duality holds. Recall the transfer map 
\vskip 0.1in
\noindent
{\bf Definition 4.6: }{\it Suppose that $X, Y$ are two topological space such that
Poincare duality holds over $\R$. Let $f: X\rightarrow Y$. Then, the transfer
map
$$f_{!}: H^*(X, \R)\rightarrow H^*(Y, \R)\leqno(4.44)$$
is defined by $f_{!}(K)=PD(f_*(PD(K)))$.}
\vskip 0.1in
We have another natural map defined in the introduction 
$\mu : \overline \M_{g-1, k+2}
\mapsto \overline \M_{g,k}$ by gluing together the last two marked 
points. Clearly, $im(\mu)$ is also a V-submanifold of $\overline{\M}_{g,k}$.

Choose a homogeneous 
basis $\{\beta _b\}_{1\le b\le L}$ of $H^*(Y,\R)$. Let $(\eta _{ab})$ be its 
intersection matrix. Note that
$\eta _{ab} = \beta _a \cdot \beta _b =0$ if the dimensions of
$\beta _a$ and $\beta _b$ are not complementary to each other.
Put $(\eta ^{ab})$ to be the inverse of $(\eta _{ab})$. Let $\delta\subset Y
\times Y$ be the diagonal. Then, its Poincare dual 
$$\delta^*=\sum_{a,b} \eta^{ab} \beta_a\otimes \beta_b.\leqno(4.45)$$
Now we can state the composition law, which consists of two formulas.
\vskip 0.1in
\noindent
{\bf Theorem 4.7: }{\it Let $K_i \in H_*(\overline \M_{g_i,k_i+1}, \R)$ $(i=1,
2)$ and $K_0 \in H_*(\overline \M_{g-1,k +2}, \R)$. For any $\alpha _1,
\cdots,\alpha _k$ in $H^*(Y,\R)$.
Then we have
(1).$$\begin{array}{rl}
&\Psi ^Y_{(A,g,k)}((\theta _{S})_{!}(K_1\times K_2])\{\alpha _i\})\\
=(-1)^{deg(K_2)\sum^{k_1}_{i=1} deg (\alpha_i)} ~& \sum \limits _{A=A_1+A_2} \sum \limits_{a,b}
\Psi ^Y_{(A_1,g_1,k_1+1)}(K_1;\{\alpha _{i}\}_{i\le k}, \beta _a)
\eta ^{ab}
\Psi ^Y_{(A_2,g_2,k_2+1)}(K_2;\beta _b,
\{\alpha _{j}\}_{j> k}) \\
\end{array}
\leqno (4.46) 
$$
(2).$$
\Psi ^Y_{(A,g,k)}(\mu_{!}(K_0);\alpha _1,\cdots, \alpha _k)
=\sum _{a,b} \Psi ^Y_{(A,g-1,k+2)}( K_0;\alpha _1,\cdots, \alpha _k,
\beta _a,\beta _b) \eta ^{ab}\leqno (4.47)
$$
}
\vskip 0.1in
{\bf Proof: } The proof of the theorem is divided into two steps.
First of all,
$$\chi_{g,k}: \overline{\B}_A(Y,g,k)\rightarrow \overline{\M}_{g,k} \leqno(4.48)$$
is a submersion.  $\B_{im(\theta_S)}=\chi^{-1}_{g,k}(Im(\theta_S))$ is a union of
some lower strata of $\overline{\B}_A(Y,g,k)$. Moreover, it is also weakly smooth.
Consider weakly smooth Banach compact-V triple $(\B_{im(\theta_S)},\F_{im(\theta_S)}, S_{im(\theta_S)})$. We can
use it to define invariant $\Psi_{(A, \theta_S)}$. The first step is to show
that
$$\Psi^Y_{(A,g,k)}(i_{!}(K); \alpha_1, \cdots, \alpha_k)=\Psi_{(A, \theta_S)}(
K; \alpha_1, \cdots, \alpha_k),\leqno(4.49)$$
Let $(im(\theta_S))^*$ be the Poincare dual of $im(\theta_S)$. $(im(\theta_S))^
*$ can be chosen to be supported in a tubular neighborhood of $im(\theta_S)$,
which can be identified with a neighborhood of zero section of normal bundle.
For any $K\in H^*(im(\theta_S), \R)$, we can pull it back to the total space
of normal bundle (denoted by $K_{\overline{\M}_{g,k}}$). Then, $K_{\overline{
\M}_{g,k}}$ is defined over a tubular neighborhood of $im(\theta_S)$. Since
$(im(\theta_S))^*$ is supported in the tubular neighborhood, 
$$(im(\theta_S))^*\wedge K_{\overline{\M}_{g,k}} \leqno(4.50)$$
is a closed differential form defined over $\overline{\M}_{g,k}$. In fact,
$$i_{!}(K)=(im(\theta_S))^*\wedge K_{\overline{\M}_{g,k}}.\leqno(4.51)$$
First we construct that $(\E, s)$ for $(\B_{im(\theta_S)},
\F_{im(\theta_S)}, S_{im(\theta_S)})$.
Suppose that the  virtual neighborhood is
$(U_{im(\theta_S)}, E_{im{\theta_S}}, S_{im(\theta_S)})$. We first extend $s$ over a neighborhood in
$\overline{\B}_A(Y,g,k)$. Then, we construct $s'$ supported away from $im(\theta_S)$. Suppose
that the stabilization term is $(\E\oplus \E', s+s')$ such that 
$$L_x+s+s'+\delta(\chi_{g,k}): T_x\B_{g,k}\oplus \E\oplus \E' 
\rightarrow \F_x\times T_{\chi_{g,k}(x)}\overline{\M}_{g,k}$$
is surjective over $\U$ in the construction of (4.14-4.16).  Suppose that the
resulting finite dimensional virtual neighborhood is $(U_{g,k}, E\oplus E', S_{g,k})$. Then, 
$$\chi_{g,k}: U_{g,k}\rightarrow \overline{\M}_{g,k}\leqno(4.52)$$
is a submersion and
$$\chi^{-1}_{g,k}(im(\theta_S))=E'_{U_{im(\theta_S)}}\subset 
U_{g,k}\leqno(4.53)$$
is a V-submanifold. Then, $\chi^*_{g,k}((im(\theta_S))^*)$ is Poincare dual to
$E'_{U_{im(\theta_S)}}$. Choose Thom forms $\Theta_1, \Theta_2$ of $E,E'$
Therefore,
$$\begin{array}{lll}
\Psi^Y_{(A,g,k)}(i_{!}(K); \alpha_1, \cdots, \alpha_k)&=&\int_{U_{g,k}}(im(\theta_S))^*)
\wedge \chi^*_{g,k}(K_{\overline{\M}_{g,k}} \wedge \Xi^*_{g,k}(\prod_i
\alpha_i)\wedge S^*_{g,k}(\Theta_1\wedge\Theta_2)\\
&=&\int_{U_{im(\theta_S)}\times \R^{m'}/G'}\chi^*_{
g,k}(K_{\overline{\M}_{g,k}})\wedge \Xi^*_{g,k}(\prod_i
\alpha_i)\wedge S^*_{g,k}(\Theta_1\wedge\Theta_2)\\
&=&\int_{U_{im(\theta_S)}}\chi^*_{g,k}(K)\wedge\Xi^*_{g,k}(\prod_i
\alpha_i)\wedge S^*_{g,k}(\Theta_1)\\
&=&\Psi^Y_{(A, \theta_S)}(K; \alpha_1, \cdots, \alpha_k)
\end{array}.\leqno(4.55)$$

The second step is to show that $\Psi^Y_{(A, \theta_S)}$ can be expressed
by the formula (1).
By the construction in the last section, we have a submersion
$$e^{A_1}_{k_1+1}\times e^{A_2}_{k_2+1}: \overline{\B}_A(Y, g_1, k_1+1)\times 
\overline{\B}_A(Y,g_2, k_2+1)\rightarrow Y\times Y\leqno(4.56)$$ 
such that 
$$\bigcup_{A_1+A_2=A}(e^{A_1}_{k_1+1}\times e^{A_2}_{k_2+1})^{-1}(\Delta)=
\B_{Im(\theta_S)}\leqno(4.57)$$
where $\delta\subset Y\times Y$ is the diagonal. By Gromov-compactness theorem,
there are only finite many such pairs $(A_1, A_2)$ we need. Notes that
$$(e^{A_1}_{k_1+1}\times e^{A_2}_{k_2+1})^{-1}(\Delta)\cap (e^{A'_1}_{k_1+1}
\times e^{A'_2}_{k_2+1})^{-1}(\Delta)\leqno(4.58)$$
may be nonempty for some $(A_1, A_2)\neq (A'_1, A'_2)$. But it is in lower strata
 of $\B_{im(\theta_S)}$ of codimension at 
least two. Furthermore, by the construction of the section 3
$$\overline{\F}_A(Y,g,k)|_{\B_{im(\theta_S)}}=\overline{\F}_A(Y,g_1, k_1+1)
\times \overline{\F}_A(Y, g_2, k_2+1)|_{\B_{
im(\theta_S)}}.\leqno(4.59)$$
We want to construct  {\em a system of stabilization terms compatible with the stratification}. 
The idea is to start from the bottom strata and construct inductively the stabilization term
supported away from lower strata. The same construction is crucial in the construction of Floer homology.
We choose to wait until the last section to give the detail (called a system of stablization terms
compatible with the corner structure in the last section). 
Let $s_1, s_2$ be the stablization terms
for 
$$(\overline{\B}_{A_1}(Y, g_1, k_1+1), \overline{\F}_{A_1}(g_1, K_1+1), 
\S^{A_1}_{g_1, k_1+1}), (\overline{\B}_{A_2}(Y,g_1, k_1+1), \overline{\F}_{A_2}
(Y,g_1, K_1+1), \S^{A_2}_{g_1, k_1+1}).$$
  Suppose that
the resulting 
 virtual neighborhoods are
 $$(U^{A_1}_{g_1,k_1+1}, E, S^{A_1}_{g_1, k_1+1}), (U^{A_2}_{g_1,k_1+1}, E', S^{A_2
}_{g_1, k_1+1}).$$
 By (4.56) and adding sections if necessary, we can assume that
$$e^{A_1}_{k_1+1}\times e^{A_2}_{k_2+1}: U^{A_1}_{g_1, k_1+1}\times U^{A_2}_{
g_2, k_2+1}\rightarrow Y\times Y\leqno(4.60)$$
is a submersion.  Let
$$U_{A_1, A_2}=(e^{A_1}_{k_1+1}\times e^{A_2}_{k_2+1})^{-1}(\Delta)\subset 
U^{A_1}_{g_1, k_1+1}\times U^{A_2}_{
g_2, k_2+1}.\leqno(4.61)$$
One consequence of our system of stabilization compatible with the stratification is  
$$U_{A_1, A_2}\cap U_{A'_1, A'_2}$$
is a V-submanifold of codimension at least two for both $U_{A_1, A_2}$ and
$U_{A'_1, A'_2}$ if $(A_1, A_2)\neq (A'_1, A'_2)$. Then, 
$$(\bigcup_{A_1+A_2=A} U_{A_1,A_2}, E\oplus E', S_{A_1}\times
S_{A_2})\leqno(4.62)$$
is a finite dimensional virtual neighborhood of $(\B_{im(\theta_S)}, \F_{im(
\theta_S)}, \S_{im(\theta_S)})$. Moreover, we can choose stabilization term such that
both $E$ and $E'$ are of even rank. Let $\delta^*$ be the Poincare dual of $\delta$. 
Then, $(e^{A_1}_{g_1, k_1+1}
\times e^{A_2}_{g_2, k_2+1})^*(\delta^*)$ is Poincare dual to $U_{A_1, A_2}$.
 Therefore,
$$\begin{array}{lll}
&&\Psi^Y_{(A, \theta_S)}(K_1\times K_2; \{\alpha_i\})\\
&=&\int_{\cup_{A_1+A_2=A} 
U_{A_1,A_2}} \chi^*_{g,k}(K_1\times K_2)\wedge \Xi^*_{g,k}(\prod_i \alpha_i)
\wedge S^*_{A_1}(\Theta_1)\wedge S^*_{A_2}(\Theta_2)\\
&=&\sum_{A_1+A_2=A}\int_{U_{A_1, A_2}} \chi^*_{g,k}(K_1\times K_2)\wedge \Xi^*_{g,k}(\prod_i \alpha_i)
\wedge S^*_{A_1}(\Theta_1)\wedge S^*_{A_2}(\Theta_2)\\
&=&\sum_{A_1+A_2=A} \int_{U^{A_1}_{g_1, k_1+1}\times U^{A_2}_{g_2, k_2+1}}
(e^{A_1}_{g_1, k_1+1}\times e^{A_2}_{g_2, k_2+1})^*(\delta^*)\wedge\chi^*_{g,k}(K_1\times K_2) \\
&&\wedge \Xi^*_{g,k}(\prod_i \alpha_i)
\wedge S^*_{A_1}(\Theta_1)\wedge S^*_{A_2}(\Theta_2)\\
&=&\sum_{A_1+A_2=A}\sum_{a,b}\eta^{ab}\int_{U^{A_1}_{g_1, k_1+1}\times U^{A_2}_{g_2, k_2+1}}
 (e^{A_1}_{g_1, k_1+1})^*{\beta_a}\wedge (e^{A_2}_{g_2, k_2+1})^*(\beta_b)\\
&&\wedge \chi^*_{g,k}(K_1\times K_2)\wedge \Xi^*_{g,k}(\prod_i \alpha_i)
\wedge S^*_{A_1}(\Theta_1)\wedge S^*_{A_2}(\Theta_2)\\
&=&(-1)^{deg(K_2)\sum^{k_1}{i=1} deg(\alpha_i)}\sum_{A_1+A_2=A}\\
&&\sum_{a,b}\eta^{ab}
\int_{U^{A_1}_{g_1, k_1+1}} \chi^*_{g_1, k_1+1}(K_1)\wedge \Xi^*_{g_1, k_1}(
\prod^{k_1}_i \alpha_i)(e^{A_1}_{g_1, k_1+1})^*{\beta_a}\wedge S^*_{A_1}(\Theta_1)\\
&&\int_{U^{A_2}_{g_2, k_2+1}}(
\chi^*_{g_2, k_2+1}(K_2)e^{A_2}_{g_2, k_2+1})^*{\beta_b}\wedge \Xi^*_{
g_2, k_2}(\prod_{j>k_1} \alpha_j)\wedge S^*_{A_2}(\Theta_2)\\
&=&(-1)^{deg(K_2)\sum^{k_1}{i=1} deg(\alpha_i)}\sum_{A_1+A_2=A}\\
&&\sum_{a,b}\eta^{ab}\Psi^Y_{(A,g_1,k_1+1)}(K_1; \{\alpha_i\}_
{i\leq k_1}, \beta_a)\Psi^Y_{(A_2, g_2, k_2+1)}(K_2; \{\alpha_j\}_{j>k_1}, 
\beta_b).
\end{array}\leqno(4.63)$$

The Proof of (2) is similar. We leave it to  readers.
\vskip 0.1in
\noindent
{\bf Corollary 4.8: }{\it Quantum multiplication is associative and hence there
is a quantum ring structure over any symplectic manifolds.}
\vskip 0.1in
{\bf Proof: } The proof is well-known (see \cite{RT1}). We omit it. $\Box$

Here, we give another application to higher dimensional algebraic geometry.
Recall that a Kahler manifold $W$ is called uniruled if $W$ is covered by 
rational curves. As we mentioned in the beginning, Kollar showed that if $W$
is a 3-fold, the uniruledness is a symplectic property \cite{K1}. Combined 
Kollar's argument with our construction, we generalize this result to general
symplectic manifolds.
\vskip 0.1in
\noindent
{\bf Proposition 4.9: }{If a smooth Kahler manifold $W$ is symplectic deformation 
equivalent to a uniruled manifold, $W$ is uniruled.}
\vskip 0.1in
First we need following 
\vskip 0.1in
{\bf Lemma 4.10: }{\it Suppose that $N\subset Y$ is a submanifold such that
for any $x\in \M_N=(\overline{\M}_A(Y,g,k,J)\cap e^{-1}_1(N))$
$$Coker L_x=0 \mbox{ and } \delta(e_1): L_x\rightarrow T_{e_1(x)}Y\leqno(4.64)$$
is surjective onto the normal bundle of $N$. Then, $\M_N$ is a smooth 
V-manifold of dimensional $ind-Cod(N)$ and
$$\Psi^Y_{(A,g,k+1)}(K; N^*, \alpha_1, \cdots, \alpha_k)=(-1)^{deg(K)deg(N^*)}
\int_{\M_N}\chi^*_{g,k+1}(K)\wedge \Xi^*_{g,k}(\prod_i\alpha_i).\leqno(4.65)$$}
\vskip 0.1in
{\bf Proof: } Since $e_1: \overline{\B}_{g,k+1}\rightarrow Y$ is a submersion, 
 we can construct $(\E, s)$ such that $s=0$ over a neighborhood
of $\M_N$ and 
$$e_1|_U: U\rightarrow Y$$
is transverse to $N$, where $(U, E, S)$ is the   virtual
neighborhood constructed by $(\E, s)$. Therefore, 
$$(e_1|_U)^{-1}(N)=E_{\M_N}\leqno(4.66)$$
is a smooth V-submanifold of $U$. Thus, $e^*_1(N^*)$ is Poincare dual to
$E_{\M_N}$. 
$$\begin{array}{lll}
\Psi^Y_{(A,g,k+1)}(K; N^*, \alpha_1, \cdots, \alpha_k)&=&\int_U \chi^*_{g,k+1}
(K)\wedge e^*_{1}(N^*)\wedge\Xi^*_{g,k}(\prod_i\alpha_i)\wedge S^*(\Theta)\\
&=&(-1)^{deg(K)deg(N^*)}\int_{E_{\M_N}}\chi^*_{g,k+1}
(K)\wedge\Xi^*_{g,k}(\prod_i\alpha_i)\wedge S^*(\Theta)\\
&=&(-1)^{deg(K)deg(N^*)}\int_{\M_N}\chi^*_{g,k+1}(K)\wedge\Xi^*_{g,k}(\prod_i\alpha_i)
\end{array}.\leqno(4.67)$$
$\Box$
\vskip 0.1in
{\bf Proof of Proposition 4.9:} If $\Psi^Y_{(A,0,k+1)}(K; pt, \cdots)\neq 0$, then
$W$ is covered by rational curves. Otherwise, there is a point $x_0$ where there
is no rational curve passing through $x_0$.  
$$\M_N=\overline{\M}_A(Y,0,k,J)\cap e^{-1}_1(N)=\emptyset\leqno(4.68)$$
for any $A, k$. The condition of Lemma 4.10 is obviously satisfied. By 
Lemma 4.10, 
$$\Psi^Y_{(A,0,k+1)}(K; pt, \cdots)=0$$
 and this is a contradiction.

 Since GW-invariant $\Psi^Y_{(A,0,k+1)}(K; pt, \cdots)$ is a symplectic
deformation invariant property, it is enough to show that if $W$ is uniruled,
$\Psi^Y_{(A,0,k+1)}(K; pt, \alpha_1, \cdots, \alpha_k)\neq 0$ for some
$K$, $\alpha_1, \cdots, \alpha_k$. Assuming Lemma 4.10, Kollar showed
some $\Psi^W_{(A,0,3)}(pt; pt, \alpha, \beta)$ is not zero for some $A$ and $\alpha, \beta$.
 His argument uses Mori's machinery. Here we give a more
elementary  argument to show that 
$$\Psi^Y_{(A,0,k+1)}(pt; pt, \alpha_1, \cdots, \alpha_k)\neq 0\leqno(4.69)$$
for some $A$ and some $\alpha_i$ with $k>>0$. Then, using the composition law we proved, we can
derive Kollar's calculation.

First, we repeat some of Kollar's argument. By \cite{K}, for a generic point
$x_0$, $\M_A(W, 0,k,J)\cap e^{-1}_1(x_0)$ satisfies the condition of Lemma 4.10 for any
$A$. Next choose $A_0$ such that
$$H(A_0)=min_A\{H(A); \M_A(W,0,k,J)\cap e^{-1}_1(x_0)\neq \emptyset\}.
\leqno(4.70)$$
where $H$ is an ample line bundle. Then, one can check that 
$$(\overline{\M}_A(W,0,k,J)-\M_A(W,0,k,J))\cap e^{-1}_1(x_0)=\emptyset.$$
Furthermore, $\M_{x_0}=\M_A(W, 0,k,J)\cap e^{-1}_1(x_0)$ is  a compact, smooth,
complex manifold. In particular, it carries a fundamental class. 

Next, we show that
$$\Xi_{0,k}: \M_{x_0}\rightarrow W^k \leqno(4.71)$$
is an immersion for large $k>>0$. For any $f\in \M_{x_0}$, 
$$T_f\M_{x_0}=\{v\in H^0(f^*TV); v(x_0)=0\}.\leqno(4.72)$$
Since $v_f\in H^0(f^*TV)$ is holomorphic, there are finite many points $x_2,
\cdots, x_{k+1}$ such that if for any $v_f$ with $v_f(x_i)=0$ for every $i$, 
$v_f=0$ . One can check
that
$$\delta (\Xi_{0,k})_f(v)=(v(x_2), \cdots, v(x_k)).\leqno(4.73)$$
Therefore, $\delta(\Xi_{0,k})$ is injective.

Since $\Xi_{0,k}$  is an immersion, $\Xi_{0,k}(\M_{x_0})\subset W^k$ is a compact
complex subvariety of the same dimension. Hence, it carries a nonzero homology
class $[\Xi_{0,k}(\M_{x_0})]$. Furthermore, $(\Xi_{0,k})_*([\M_{x_0}])=\lambda
[\Xi_{0,k}(\M_{x_0})]$ for some $\lambda>0$. By Poincare duality, there are 
$\alpha_1, \cdots, \alpha_k$ such that 
$$\prod_i \alpha_i([\Xi_{0,k}(\M_{x_0})])\neq 0. \leqno(4.74)$$
By Lemma 4.10,
$$\begin{array}{lll}
\Psi^W_{(A,g,k+1)}(pt; pt, \alpha_1, \cdots, \alpha_k)&=&
         \int_{\M_{x_0}}\Xi^*_{0,k}(\prod_i \alpha_i)\\
&=&(\prod_i \alpha_i)(\Xi_*([\M_{x_0}]))\neq 0
\end{array}\leqno(4.75)$$
$\Box$

\section{Equivariant GW-invariants and Equivariant quantum cohomology}
We will study the equivariant GW-invariants and the equivariant quantum cohomology in
detail in this section. The equivariant theory is an important topic. It has
been studied by several authors \cite{AB}, \cite{GK}. As we mentioned in the 
\cite{R4}, equivariant theory is the one that usual Donaldson method has
trouble to deal with, where there are topological obstructions to choose a 
``generic'' parameter. But our virtual neighborhood method is particularly
suitable to study equivariant theory.  In our case, one can attempt to choose
an equivariant almost complex structure and apply the equivariant virtual neighborhood 
technique.
However, a technically simpler approach is to view the equivariant GW-invariants
 as a limit of GW-invariants for the families of symplectic manifolds. This
approach was advocated by \cite{GK}, where they formulated some conjectural
properties for the equivariant GW-invariants and the equivariant quantum cohomology.
First work to give a rigorous foundation of the equivariant GW-invariants and the
equivariant quantum
cohomology was given by Lu \cite{Lu} for monotonic symplectic manifolds, where
he used the  method of \cite{RT1}, \cite{RT2}. Here, we use the invariants
we established in last section to establish the equivariant GW-invariants and the
equivariant quantum cohomology for general symplectic manifolds.

Let's recall the construction of the introduction. 
Suppose that $G$ acts on $(X, \omega)$ as symplectomorphisms. Let $BG$ be the 
classifying space of $G$ and $EG\rightarrow BG$ be the 
universal $G$-bundle. Suppose that 
$$BG_1\subset BG_2\cdots\subset BG_m \subset BG \leqno(5.1)$$
such that $BG_i$ is a smooth oriented compact manifold and $BG=\cup_i BG_i$. Let
$$EG_1\subset EG_2\cdots\subset EG_m \subset BG\leqno(5.2)$$
be the corresponding universal bundles. We can also form the approximation of
homotopy quotient $X_G=X\times EG/G$ by $X^i_G=X\times EG_i/G$. Since $\omega$
is invariant under $G$, its pull-back on $V\times EG_i$ descends to $V^i_G$.
So, we have a family of symplectic manifold $P_i: X^i\rightarrow BG_i$. 
Applying our previous construction, we obtain GW-invariant $\Psi^{X^i_G}_{(A, 
g,k)}$. We define equivariant GW-invariant 
$$\Psi^G_{(A,g,k)}\in Hom(H^*(\overline{\M}_{g,k}, \R)\otimes(H^*(V_G, \R))^{
\otimes k}, H^*(BG, \R)) \leqno(5.3)$$
as follow:

 For any $D\in H_*(BG, \Z)$, $D\in H_*(BG_i, \Z)$ for some $i$. For any $K\in
H^*(\overline{\M}_{g,k}, \R)$, $\pi^*(K)\in H^*(\overline{\M}_{g,k+1}, \R)$.
Let $i_{X^i_G}: X^i_G\rightarrow  X_G$. 
\vskip 0.1in
\noindent
{\bf Definition 5.1: }{\it For $\alpha_i\in H^*_G(X,\R)$, we define 
$$\Psi^G_{(A,g,k)}(K, \alpha_1, \cdots, \alpha_k)(D)=\Psi^{X^i_G}_{(A, g,k+1)}
(\pi^*(K); i^*_{X^i_G}(\alpha_1), \cdots, i^*_{X^i_G}(\alpha_k), P^*_i(D^*_{BG_i})),\leqno(5.4)$$
where $D^*_{BG_i}$ is the Poincare dual of $D$ with respect to $BG_i$.}
\vskip 0.1in
\noindent
{\bf Theorem 5.2: }{\it (i). $\Psi^G_{(A, g,k)}$ is independent of the choice of
$BG_i$.
\vskip 0.1in
\noindent
(ii). If $\omega_t$ is a family of $G$-invariant symplectic forms, $\Psi^G_{(A,
g,k)}$ is independent of $\omega_t$.}
\vskip 0.1in
{\bf Proof: } The proof is similar to the third step of the proof of Proposition 4.2(iv).
 Choose a $G$-invariant tamed almost complex structure $J$ on $X$.
It induces a tamed almost complex structure (still denoted by $J$) over every
$X^i_G$. Clearly, there is a natural inclusion map
$$\overline{\M}_A(X^i_G, g,k, J)\subset \overline{\M}_A(X^j_G, g, k,J) \mbox{
 for } i\leq j. \leqno(5.5)$$
Suppose that $(\B_i, \F_i, \S_i)$ is the configuration space of 
$\overline{\M}_A(X^i_G,g,k,J)$. Then, there is a natural inclusion.
$$(\B_i, \F_i, \S_i)\subset (\B_j, \F_j, \S_j) \mbox{ for } i\leq j.\leqno(5.6)$$
 We first construct $(\E_i, s_i)$ for $(\B_i, \F_i, 
\S_i)$. Suppose that the resulting finite dimensional virtual neighborhood is
$(U_i, E_i, S_i)$.  Then, we extend $s_i$ over $\B_j$. Since
$L_A+s_i$ is surjective over $\U_i\subset\B_i$. We can construct 
$(\E_j, s_j)$
such that $s_j=0$ over $\U_i$ and $L_A+s_i+s_j$ is surjective over $\U_j$.
 Suppose that the resulting
finite dimensional virtual neighborhood is $(U_j, E_i\oplus E_j,
S_j)$. Then, 
$$U_j\cap (\E_j)_{\B_i}=(E_j)_{U_i}\subset U_j$$
 is a V-submanifold. Let
$$e^j_{k+1}: \B_j\rightarrow X^j_B\leqno(5.7)$$
be the evaluation map at $x_{k+1}$. Then,  we can choose $s_i, s_j$ such that
  the restriction of
$e^j_{k+1}$ to $U_j$ is a submersion. Furthermore, since $(e^j_{k+1})^{-1}(
X^i_G)=\B_i$, 
$$(e^j_{k+1})^{-1}(X^i_G)\cap U_j= (E_j)_{U_i}.\leqno(5.8)$$
Notes that 
$$S_j\circ i=S_i, \leqno(5.9)$$
where $i:(E_j)_{U_i} \rightarrow U_j$ is the inclusion. Choose Thom forms
$\Theta_i, \Theta_j$  of $E_i,E_j$.
Let's use $I_{ij}$ to denote the inclusion $\B_i\subset \B_j$, $BG_i\subset 
BG_j$ and $X^i_G\subset X^j_G$ and define $\Xi^i_{g,k}, \chi^i_{g,k}$ similarly.
Then
$$\Xi^j_{g,k}\circ I_{ij}=I_{ij}\Xi^i_{g,k}, \mbox{ and } \chi^j_{g,k}\circ 
I_{ij}=\chi^i_{g,k}.\leqno(5.11)$$
Furthermore, 
$$D^*_{BG_j}=(I_{ij})_{!} D^*_{BG_i}.\leqno(5.12)$$
Let $(BG_i)^*_j$ be the Poincare dual of $BG_i$ in $BG_j$. Choose $(BG_i)^*_j$ 
supported in a tubular neighborhood of $BG_i$. By Lemma 2.10,
$$D^*_{BG_j}=(D^*_{BG_i})_{BG_j}\wedge (BG_i)^*_j.\leqno(5.13)$$
Furthermore, $P^*_j((BG_i)^*_j)$ is Poincare dual to $X^i_G$ in $X^j_G$. Hence,
$(e^j_{k+1})^*P^*_j((BG_i)^*_j)$ is Poincare dual to $(E_j)_{U_i}$.
$$\begin{array}{lll}
&&\Psi^{X^j_G}_{(A,g,k+1)}(\pi^*(K), i^*_{X^j_G}(\alpha_1),\cdots, i^*_{X^j_G}(\alpha_{
k}), P^*_j(D^*_{BG_j}))\\
&=&\int_{U_j} \chi^j_{g,k+1}(\pi^*(K))\wedge \Xi^i_{g,k}(\prod_m i^*_{X^j_G}(
\alpha_m))\wedge (e^j_{k+1})^*P^*_j(D^*_{BG_j})\wedge S_j^*(\Theta_i\times 
\Theta_j)\\
&=&\int_{U_j}\chi^j_{g,k+1}(\pi^*(K))\wedge \Xi^j_{g,k}(\prod_m i^*_{X^j_G}(
\alpha_m))\wedge (e^j_{k+1})^*P^*_j((D^*_{BG_i})_{BG_j})\wedge (e^j_{k+1})^*P^*_j(
(BG_i)^*_j)\wedge S_j^*(\Theta_i\times \Theta_j)\\
&=&\int_{(E_j)_{U_i}}\chi^i_{g,k+1}(\pi^*(K))\wedge \Xi^i_{g,k}(\prod_m i^*_{X^i_G}(
\alpha_m))\wedge (e^i_{k+1})^*P^*_i(D^*_{BG_i})\wedge S_j^*(\Theta_i\times \Theta_j)\\
&=&\int_{U_i}\chi^i_{g,k+1}(\pi^*(K))\wedge \Xi^i_{g,k}(\prod_m i^*_{X^i_G}(
\alpha_m))\wedge (e^i_{k+1})^*P^*_i(D^*_{BG_i})\wedge S_i^*(\Theta_i)\\
&=&\Psi^{X^i_G}_{(A,g,k+1)}(\pi^*(K); i^*_{X^i_G}(\alpha_1),\cdots, i^*_{X^i_G}
(\alpha_{k}), P^*_i(D^*_{BG_i}))
\end{array}.\leqno(5.14)$$
(ii) follows from the same property of $\Psi^{X^i_G}$. $\Box$

As we discussed in the introduction,  for any  equivariant cohomology class 
$\alpha\in H^*_G(X)$, we can evaluate over the fundamental class of $X$
$$\alpha [X]\in H^*(BG).\leqno(5.15)$$
Furthermore, there
is a modulo structure by $H^*_G(pt)=H^*(BG)$, defined by using the projection map
$$X_G\rightarrow BG.\leqno(5.16)$$
The equivariant quantum multiplication is a new multiplication structure
over $H^*_G(X, \Lambda_{\omega})=H^*(X_G, \Lambda_{\omega})$ as follows. We first
define a total 3-point function
$$\Psi^G_{(X,\omega)}(\alpha_1, \alpha_2, \alpha_3)=\sum_A \Psi^G_{(A,0,3)}(
pt; \alpha_1, \alpha_2, \alpha_3)q^A.\leqno(5.17)$$
Then, we define an equivariant quantum multiplication by
$$(\alpha\times_{QG}\beta)\cup \gamma [X]=\Psi^G_{(X,\omega)}(\alpha_1, 
\alpha_2, \alpha_3).\leqno(5.18)$$
\vskip 0.1in
\noindent
{\bf Theorem I: }{\it (i) The equivariant quantum multiplication is
 skew-symmetry.
\vskip 0.1in
\noindent
(ii) The equivariant quantum multiplication is commutative
with the modulo structure of $H^*(BG)$.
\vskip 0.1in
\noindent
(iii) The equivariant quantum multiplication is associative.

Hence, there is a equivariant quantum ring structure for any $G$ and any symplectic
manifold $V$}
\vskip 0.1in
{\bf Proof: } (i) follows from the definition.
By the proposition 5.2, for any $\alpha\in H^*(BG, \R)$,
$$\begin{array}{lll}
&&\Psi^{X^i_G}_{(A,g,k+1)}(\pi^*(K); i_{X^i_G}^*(\alpha_1 ),\cdots, P^*_i(i_{
BG_i})^*(\alpha)\wedge i_{X^i_G}^*(\alpha_j), \cdots,i_{X^i_G}(\alpha_k), 
P^*_i(D^*_{BG_i}))\\
&=&\Psi^{X^i_G}_{(A,g,k+1)}(\pi^*(K);  i_{X^i_G}^*(\alpha ), 
i_{X^i_G}^*(\alpha_2), \cdots,i_{X^i_G}(\alpha_k), P^*_i(i_{
BG_i})^*(\alpha)\wedge P^*_i(D^*_{BG_i}))\\
&=&\Psi^{X^i_G}_{(A,g,k+1)}(\pi^*(K);  i_{X^i_G}^*(\alpha ), 
i_{X^i_G}^*(\alpha_2), \cdots,i_{X^i_G}(\alpha_k), P^*_i(i_{
BG_i})^*(\alpha\wedge D^*_{BG_i}))\\
&=&\Psi^{X^i_G}_{(A,g,k+1)}(\pi^*(K);  i_{X^i_G}^*(\alpha ), 
i_{X^i_G}^*(\alpha_2), \cdots,i_{X^i_G}(\alpha_k), P^*_i(i_{
BG_i})^*((\alpha(D)^*_{BG_i}))\\
&=&\Psi^G_{(A,g,k)}(K, \alpha_1, \alpha_2, \cdots, \alpha_k)(\alpha(D))
\end{array}. \leqno(5.19)$$
Then, (ii) follows from the definition.

The proof of (iii) is the same as the case of the ordinary quantum cohomology.
We omit it. $\Box$

\section{Floer homology and Arnold conjecture}
  In this section, we will extend our construction of previous sections to
Floer homology to remove the semi-positive condition. Floer homology was first
introduced by Floer in an attempt to solve Arnold conjecture \cite{F}. The 
original Floer
homology  was only defined for monotonic symplectic manifolds. Floer 
solved Arnold conjecture in the same category. 
The Floer homology for semi-positive symplectic
manifolds was defined by Hofer and Salamon \cite{HS}.  Arnold conjecture
for semi-positive symplectic manifolds were solved by 
\cite{HS} and \cite{O}.  Roughly
speaking, there are two difficulties to solve Arnold conjecture for general
symplectic manifolds,i.e., (i) to extend Floer homology to general 
symplectic manifolds and (ii) to show that Floer homology is the same as 
ordinary homology. For the second problem, the traditional method is to deform
a Hamiltonian function to a small Morse function and calculate its Floer homology
directly. This approach involved some delicate analysis about the contribution
of trajectories
which are not gradient flow lines of a Morse function. It has only been carried
out for semi-positive symplectic manifolds \cite{O}. But the author and Tian
showed \cite{RT3} that this part of difficulties can be avoided by introducing
a Bott-type Floer homology, where we can deform a Hamiltonian function to zero.
The difficulty to extend Floer homology for  a general symplectic manifold is
 the same as the difficulty to
extend GW-invariant to a general symplectic manifold. Once we establish the
GW-invariant for general symplectic manifolds, it is probably not surprising to
experts that the same technique can work for Floer homology. Since many of the
construction here is similar to that of last several sections,  we shall be 
sketch in this section.

 Let's recall the set-up of \cite{HS}. Let $(X,\omega)$ be a closed symplectic manifold. Given any function
 $H$ on $X\times S^1$, we can associate a vector field $X_H$ 
as follow:
$$ \omega (X_H(z,t), v) = v(H)(z,t),~~~~~{\rm for~any}~v~\in ~T_zV
\leqno (6.1)$$
We call $H$ a periodic Hamiltonian and $X_H$ a Hamiltonian vector field
associated to $H$. Let $\phi_t(H)$ be the integral flow of the Hamiltonian
vector field $X_H$. Then $\phi _1(H)$ is a Hamiltonian symplectomorphism.
\vskip 0.1in
\noindent
{\bf Definition 6.1.} {\it  We 
call a periodic Hamiltonian $H$ to be non-degenerate 
if and only if the fixed-point set $F(\phi_1(H))$ of $\phi_1(H)$
is  non-degenerate. }
\vskip 0.1in
\noindent
Let $\L(X)$ be the space of contractible maps (sometimes called
contractible loops)
from $S^1$ 
into $X$ and $\tilde{\L}(X)$ be the universal cover of $\L(X)$, 
namely, $\tilde{\L}(X)$ is as follows:
$$\tilde{\L}(X)=\{(x,u)|x\in \L(X), u:D^2\rightarrow X \mbox{ such 
that } x=u|_{\partial D^2}\}/\sim,
\leqno(6.2)$$
where the equivalence relation $\sim$ is the homotopic equivalence of $x$.
The covering group of $\tilde {\L}$ over $\L$
is $\pi_2(V)$. We can
define a symplectic action functional on $\tilde{\L}(X)$, 
$$a_H((x,u))=\int_{D^2 }u^*\omega+\int_{S^1} H(t, x(t)) dt
\leqno(6.3)$$
It follows from the closeness of $\omega $
that $a_H$ descends to the quotient space 
by $\sim$. The Euler-Lagrange equation of $a_H$ is
$$\dot{u}-X_H(t, u(t))=0 \leqno(6.4)$$
Let $\R(H)$ be the critical point set of $a_H$, i.e., the set of smooth 
contractible loops satisfying the Euler-Lagrange equation. 
The image $\bar{\R}(H)$ of $\R(H)$ in $\L(V)$ one-to-one corresponds to the fixed points of
$\phi_1(H)$ and hence is a finite set. Since $\phi_1(H)$ is non-degenerate,
it implies that $\R(H)$ is the set of non-degenerate critical points of $a(H)$.
But $\R(H)$ may have infinitely many points,
which are generated by the covering transformation group $\pi_2(V)$.

Given $(x, u)\in \R(H)$, choose a symplectic trivialization 
$$\Phi(t): \R^{2n}\rightarrow T_{x(t)}V$$
of $u^*TV$ which extends over the disc $D$. Linearizing the Hamiltonian
differential equation  along $x(t)$, we obtain a path of symplectic
matrices
$$A(t)=\Phi(t)^{-1}d\phi_t(x(0))\Phi(0)\in Sp(2n, \R).$$
Here the symplectomorphism $\phi_t: X\rightarrow X$ denotes the time-$t$-map
of the Hamiltonian flow
$$\dot{\phi_t}=\nabla H_t(\phi_t).$$
Then, $A(0)=Id$ and $A(1)$ is conjugate to $d\phi_1(x(0))$. Non-degeneracy
means that $1$ is not an eigenvalue of $A(1)$.  Then, we can assign a 
Conley Zehnder index for $A(t)$.
We can decomposed $\R(H)$ as
$$\R(H)=\cup_i \R_i(H),$$
where $\R_i(H)$ consists of critical points in $\R(H)$ with
the Conley-Zehnder index $i$. 

 To define  Floer homology, we first construct
a chain complex and a boundary map $(C_*(X, H), \delta)$.
The chain complex 
$$C^*(X, H)=\otimes_i C_i(X,H).\leqno(6.5)$$
where $C_i(X, H)$ is a $\R$-vector space consisting of $\sum_{\mu(\tilde{x})=i}
\xi(\tilde{x})\tilde{x}$ where the coefficients $\xi(\tilde{x})$ satisfy
the finiteness condition that 
$$\{\tilde{x}; \xi(\tilde{x})\neq 0, a_H(\tilde{x})>c\}$$
is a finite set for any $c\in \R$. We recall that the Novikov ring $\Lambda _
\omega$ is defined as the set of formal sum $\sum_{A\in \pi_2(X)}\lambda_A e^A$
such that  for each $c>0$, the number of nonzero $\lambda_A$ with $\omega(A)
\leq c$ is finite. For each $(x, u_x)\in \R(H)$, we define
$$e^A(x, u_x)=(x, u_x\#A),$$
where $\#$ is the connected sum operation in the interior of disc $u_x$. It is
easy to check that 
$$\mu(e^A(x, u_x))=2C_1(A)+\mu(x, u_x).\leqno(6.6)$$
It induces an action of Novikov ring $\Lambda_{\omega}$  on $C_*(V,H)$.

 Next we consider the boundary map, where we have to study the moduli space of
trajectories. Let $J(x)$ be a  compatible almost complex structure of $\omega$.
We can consider the perturbed gradient flow equation of $a_H$:
$$\F(u(s,t))=\frac{\partial u}{\partial s}+ J(u)\frac{\partial u}{\partial t}+ 
\bigtriangledown H(t, u)=0,$$
where we use $s$ to denote the time variable and $t$
to denote the circle variable. At this point, we ignore the homotopic class of
disc, which we will discuss later.  Let
$$\tilde{\M}=\{ u: S^1\times \R\rightarrow \R\,|\,
 \F(u)=0, E(u)=\int_{S^1\times \R}( 
|\frac{\partial u}{\partial s}|^2+|J(u)\frac{\partial u}{\partial t}+ 
\bigtriangledown
 H(t, u)|^2)dsdt<\infty\}.$$
Because $a(H)$ has only non-degenerate critical points, the following lemma is 
well-known.
\vskip 0.1in
\noindent
{\bf Lemma 6.3.} {\it
For every $u\in \tilde{\M}$, $u_s(t)=u(s, t)$ converges to $x_{\pm}
(t)\in \bar{\R}(H)$ when $s\rightarrow \pm \infty$. If $H$ is 
non-degenerate, $u_s$ 
converges exponentially to its limit, i.e., $|u_s-u_{\pm \infty}|<C e^{-\delta 
|s|}$ for $s\geq |T|.$}
\vskip 0.1in
By this lemma, we can divide $\tilde{\M}$ into
$$\tilde{\M}=\bigcup_{x^-, x^+\in \bar{\R}} \M(x^-, x^+; H, J),$$
where 
$$\tilde{\M}(x^-,x^+; H, J)=\{u\in \tilde{\M}; \lim_{s\rightarrow -\infty}u_s=x^-,
 \lim_{s \rightarrow \infty}u_s=x^+\}.$$
Clearly, $\R ^1$ acts on $\tilde{\M}(x^-,x^+; H, J)$ as translations in time. 
Let 
$$\M(x^-,x^+; H, J)=\tilde{\M}(x^-,x^+; H, J)/\R ^1.\leqno(6.7)$$ 
$\M(x^-, x^+; H, J)$
consists of the different components of different dimensions. For each $(x^-, u^-),
(x^+, u^+)\in \R(H)$, let $\M((x^-, u^-),(x^+, u^+); H,J)$ be the components of
$\M(x^-,x^+; H,J)$ satisfying that 
$$(x^+, u^-\#u)\cong (x^+, u^+)$$
 for any
$u\in\M((x^-, u^-),(x^+, u^+); H,J)$. Then, the virtual dimension of 
$\M((x^-, u^-),(x^+, u^+); H,J)$ is $\mu(x^+, u^+)-\mu(x^-, u^-)-1$.

  Next, we need a stable compactification of $\M(x^-, x^+; H, J)$.
\vskip 0.1in
\noindent
{\bf Definition 6.4: }{\it A stable trajectory (or symplectic gradient flow line)
$u$ between $x^-, x^+$ consists of 
trajectories $u_0\in \M(x^-, x_1;H,J), u_1\in \M(x_1, x_2; H, J) \cdots, u_k\in
\M(x_k, x^+)$ and finite many genus zero stable $J$-maps $f_, \cdots, f_m$ with
one marked point such that the marked point is attached to the interior of some
$u_i$. Furthermore, if $u_i$ is a constant trajectory, there is at least one stable map
attaching to it (compare with ghost bubble). We call two stable trajectories to be
equivalent if they are different by an automorphism of the domain. For each stable 
map $f$, we define $E(f)=\omega(
A)$ and denote the sum
of the energy from each component by $E(u)$. If we drop the perturbed 
Cauchy Riemann equation from the definition of trajectory and Cauchy Riemann
equation from the definition of genus zero stable maps, we simply call it a flow line. }
\vskip 0.1in
Suppose that 
$\overline{\M}((x^-,u^-), (x^+, u^+); H,J)$ is the set of the equivalence classes of
stable trajectories $u$ between
$x^-, x^+$ such that $E(u)=a(x^+)- a(x^-)$ and $(x^+, u^-\# u)\cong (x^+, u^+)$
. Let $\overline{\B}((x^-,u^-), (x^+, u^+))$ be the space of corresponding flow
lines. A slight modification of \cite{PW} shows that
\vskip 0.1in
{\bf Theorem 6.5:(\cite{PW})}{\it $\overline{\M}((x^-,u^-), (x^+, u^+); H,J)$ is 
compact.}
\vskip 0.1in
We will leave the proof to readers.

 The configuration space is $\overline{\B}_{\delta}((x^-,u^-), (x^+, 
u^+))$-the space of flow lines converging exponentially to the periodic orbits
$(x^-,u^-), (x^+, u^+)$. Next, we construct a virtual neighborhood using the
construction of  section 3. Since the construction is similar, we shall
outline the difference and leave to readers to fill out the detail. The unstable 
component is either a unstable bubble or a unstable trajectory 
$u\in \B_{\delta}((x^-,u^-), (x^+,u^+))$ where $\B_{\delta}((x^-,u^-), (x^+,u^+))$
is the space of $C^{\infty}$-map
from $S^1\times (-\infty, \infty)$ converging expentially to the periodic orbits.
When $u$ is a unstable trajectory, $u$ is a non-constant trajectory and has
no intersection point in the interior.  Theirfore, $\R$ acts freely on 
$Map_{\delta}((x^-,u^-), (x^+,u^+))$
We want to show that 
$$\B_{\delta}((x^-,u^-), (x^+,u^+))=Map_{\delta}((x^-,u^-), (x^+,u^+))/\R \leqno(6.8)$$
is a Hausdorff Frechet manifold.  Using the same method of Lemma 3.4,
we can show that 
$$\B_{\delta}((x^-,u^-), (x^+,u^+))$$
 is Hausdorff. For any $u\in
\B_{\delta}((x^-,u^-), (x^+,u^+))$, we can construct a slice 
$$W_u=\{u^w; w\in \Omega^0(u^*TV), w_s \mbox{ converges expentially to zero  and } 
||w||_{L^p_1}<\epsilon, ||w||_{C^1(D_{\delta_0}(e))}<\epsilon, w\perp \frac{\partial u}{\partial s}(e)\},\leqno(6.9)$$
where $\frac{\partial u}{\partial s}$ is injective at $e$.
Let $u\in \overline{\B}_{\delta}((x^-, u^-), (x^+, u^+))$ be a stable trajectory.
Recall that for closed case, the gluing parameter for each nodal point is $\C$.
For the trajectory, it satisfies the perturbed Cauchy Riemann equation. In particular,
the Hamiltonian perturbation term depends on the circle parameter. Therefore, the rotation
along circle is not a automorphism of the equation. The gluing parameter is only a real
number in $\R^+$. If we have more than two components of broken trajectories. The gluing
parameter is a small ball of 
$$I_k=\{(v_1, \cdots, v_k); v_i\in \R \& v_i\geq 0\},\leqno(6.10)$$
where $k+1$ is the number of broken trajectories of $u$. We call $u$  {\em a corner point}.

\vskip 0.1in
\noindent
{\bf Remark: }{\it A minor midification of Siebert's construction (Appendix) is needed
in this case. For the trajectory component, $H^0, H^1$ should be understood as the space
of sections which are exponentially decay at infinity. Recall that the vanishing theorem
 of $H^1$ was proved by certain Weitzenbock formula, which still holds in this case.}
\vskip 0.1in

The obstruction bundle $\overline{\F}_{\delta}((x^-, u^-), (x^+, u^+))$ 
can be constructed similarly. Sometimes, we shall drop $u^-, u^+$ from the notation without
any confusion.

For the corner point,  a special care is need to construct stabilizing
term $s_{x^-, x^+}$. The idea is to construct a stabilized term first in a neighborhood of bottom
strata. Then, we process to the next strata until we reach to the top. Furthermore,
we need to construct stabilization terms for all the moduli spaces of stable trajectories
at the same time. We can do it by the induction on the energy. Since there is a minimal
energy  for all the stable  trajectories, the set of the possible values of the energy of
stable trajectories are discrete. We can first construct a stabilization term for the stable
trajectories of the smallest energy and then proceed to next energy level. By the compactness
theorem, there are only finite many topological  type of stable trajectories below any 
energy level.
To simplify the notation,
let's assume that the maximal number of broken trajectories for the element of 
$\overline{\B}_{\delta}(x^-, x^+)$ is $3$ and there are three energy levels. We leave 
to  readers to fill out
the detail for general case. Suppose that $u=(u_1, u_2, u_3)$, where $u_i$ is a trajectory
connecting $x^{i-1}$ to $x^i$ attached by some genus zero stable maps.
Moreover, $x^0=x^-, x^1,x^2, x^3=x^+$. Since $u_i$ is not a corner point,
we can construct $s_{u_i}$ in the same way as  section 3. Here, we  require 
the value of $s_{u_i}$ to be compactly supported away form the gluing region. Note that
over 
$$\overline{\B}_{\delta}(x^-, x^1)\times \overline{\B}_{\delta}(x^1,
 x^2)\times \overline{\B}_{\delta}(x^2, x^+),$$ 
the obstruction bundle $\overline{\F}_{\delta}(x^-, x^+)$ is naturally decomposed as
$$\overline{\F}_{\delta}(x^-, x^1)\times \overline{\F}_{\delta}(x^1,
 x^2)\times \overline{\F}_{\delta}(x^2, x^+).\leqno(6.11)$$ 
Then, $s_{u_1}
\times s_{u_2}\times s_{u_3}$ is a section on 
$$\overline{\B}_{\delta}(x^-, x^1)\times \overline{\B}_{\delta}(x^1,
 x^2)\times \overline{\B}_{\delta}(x^2, x^+)$$ 
supported in a neighborhood of $u$. Since its value is  supported away from the gluing
region, it extends naturally over a neighborhood of $u$ in $\overline{\B}_{\delta}(
x^-, x^+)$. We multiple it by a cut-off function as we did in the section 3. 
Then, we can treat $s_{u_1}\times s_{u_2}\times s_{u_3}$ as a section supported in a 
neighborhood of $u$ in $\overline{\B}_{\delta}(
x^-, x^+)$. By the assumption, 
$$\overline{\M}(x^-, x^1)\times \overline{\M}(x^1, x^2)\times \overline{\M}(x^2, x^+)$$
is compact. We construct finite many such sections such that the linearlization of the extend equation
$$\S_e=\bar{\partial}_J+\bigtriangledown H+\sum s_{u_i}$$
is surjective over the bottom strata. 
Let 
$$s_3=\sum_i s_{u_i}$$
to indicate that it is supported in neighborhood of third strata.
Next, let's consider the next strata 
$$\overline{\M}(x^-, x^1)\times \overline{\M}(x^1, x^+) \cup
\overline{\M}(x^-, x^2)\times \overline{\M}(x^2, x^+).$$
Two components are not disjoint from each other. Then have a common boundary in the bottom
strata.  By our construction, the restriction of $s_3$ over next strata is naturally 
decomposed as
$$s^3_{(x^-, x_1)}\times s^3_{(x_1, x^+)}, s^3_{(x^-, x_2)}\times s^3_{(x_2, x^+)}.$$
Then, we construct a section of the form
$$s^2_{(x^-, x_1)}\times s^2_{(x_1, x^+)}, s^2_{(x^-, x_2)}\times s^2_{(x_2, x^+)}$$
supported away from the bottom strata. Then, we extend it over a neighborhood of the second
strata in $\overline{\B}_{\delta}(x^-, x^+)$. Over the top strata,
we construct a section supported away from the lower strata.  In general, the stabilization
term $s_{x^-, x^+}$ is the summation of $s_i$, where $s_i$ is supported in a neighborhood of $i$-th
strata and away from the lower strata. Suppose that the corresponding vector spaces
are 
$$\E^{m_{x^-, x^+}}=\prod_i \E_i.\leqno(6.12)$$
We shall choose 
$$\Theta_{x^-, x^+}=\prod_i\Theta_i, \leqno(6.13)$$
where $\Theta_i$ is  a Thom form  supported in a neighborhood of zero section of $E_i$.
with integral $1$. We call such $(s_{x^-, x^+}, \Theta_{x^-, x^+})$ {\em compatible with
the corner structure} and the set of $(s_{x^-, x^+}, \Theta_{x^-, x^+})$ for all $x^-, x^+$
{\em a system of 
 stabilization terms compatible with
the corner structure}.  Suppose that $(s_{x^-, x^+}, \Theta_{x^-, x^+})$ is compatible with
the corner structure. It has following nice property.  (i)  $s_{x^-, x^+}=s^t+s_l$, where $s^t$ is 
supported away from lower strata and $s_l$ is supported in a neighborhood of
strata. (ii) the restriction of $s_l$ to any boundary component preserves the product structure. 
Namely, we view 
$$\partial \overline{\B}_{\delta}(x^-, x^+)=\bigcup_x \overline{\B}_{\delta}(x^-, x)\times 
\overline{\B}_{
\delta}(x, x^+).\leqno(6.14)$$
The restriction of $s_l$ is of the form
$$\bigcup_x s_{x^-, x}\times s_{x, x^+}\times \{0\}.\leqno(6.15)$$

Let $(U_{x^-, x^+}, \E^{x^-, x^+}, S_{x^-, x^+})$ be the virtual 
neighborhood. Then, $U_{x^-, x^+}$ is
a finite dimensional V-manifold with the corner. 
$$\partial U_{x^-,x^+}=\bigcup_{x}E^{ot}_{U_{x^-,x}\times U_{x,x^+}},$$
where $U_{x^-,x}, U_{x, x^+}$ are the virtual neighborhoods constructed by $s_{x^-, x}, 
s_{x, x^+}$ and $E^{ot}$ is the product of other $E_i$ factors. 

When $\mu(x^+)=\mu(x^-)+1$, $dim U_{x^-, x^+}=deg \Theta_{x^-, x^+}$. We define
$$<(x^+,u^+), (x^-,u^-)>=\int_{U_{x^-,x^+}} S^*_{x^-,x^+}\Theta_{x^-, x^+},$$
where $(s_{x^-, x^+}, \Theta_{x^-, x^+})$ is compatible with the corner structure.
When $\mu(x^+)<\mu(x^-)+1$, $dim U_{x^-, x^+}<deg \Theta_{x^-, x^+}$, we define
$$<(x^+,u^+), (x^-,u^-)>=\int_{U_{x^-,x^+}} S^*_{x^-,x^+}\Theta_{x^-, x^+}=0,\leqno(6.17)$$
For any $x\in C_k(X,H)$, we  define a boundary operator as
$$\delta x=\sum_{y\in C_{k-1}}<x,y>y.\leqno(6.18)$$

Novikov ring naturally  acts on $C_*(V, H)$ by $e^A(x,u)=(x, u\#A)$ for
$A\in \pi_2(X)$. 
Furthermore, it is commutative with the boundary operator. Next, we show that
\vskip 0.1in
\noindent
{\bf Proposition 6.6: }{\it $\delta^2=0$.}
\vskip 0.1in
{\bf Proof: } 
$$\delta^2 x=\sum_{z\in C_{k-2}}\sum_{y\in C_{k-1}}<x y><y,z>z.\leqno(6.19)$$
Let $<x,z>^2=\sum_{y\in  C_{k-1}}<x y><y,z>$. It is enough to show that
$$<x, z>^2=0. \leqno(6.20)$$
Consider $\M(x, z; H, J)$. Its stable compactification $\overline{\M}(x,z; 
H,J)$ consists of broken trajectories of the form $(u_0, u_1; f_1, \dots, f_m)$
for $u_0\in \overline{\M}(x, y; H, J), u_1\in \overline{\M}(y, z; H, J)$. Choose
compatible $(s_{x,z}, \Theta_{x, z})$.   The boundary 
components
$$\partial \overline{\B}_{\delta}(x,z)=\bigcup_{y} \overline{\B}_{x,y}\times 
\overline{\B}_{y, z},\leqno(6.21)$$
where $\overline{\B}_{x, y}, \overline{\B}_{y,z}$ are the configuration spaces of
$\overline{\M}(
x,y, H, J)$, $\overline{\M}(y,z; H, J)$, respectively. Furthermore, $\overline{
\F}_{x,z}$
is naturally decomposed,i.e.,
$$\overline{\F}_{x,z}|_{\overline{\B}_{x,y}\times \overline{\B}_{y, z}}=
\overline{\F}_{x,y}\times \overline{\F}_{y,z}.\leqno(6.22)$$
Suppose that the 
resulting virtual neighborhood by $s_{x,z}$ is $(U_{x,z}, E^{x,z}, S_{x,z})$. 
Then,
$$\partial{U_{x,z}}=\bigcup_{y}E^{ot}_{ U_{x,y}\times U_{y,z}}.\leqno(6.23)$$
Note that  $dim U_{x,z}=deg \Theta_{x,z}+1$.
$$\begin{array}{ccl}
0&=&\int_{U_{x,z}}S^*_{x,z}d(\Theta_{x,z})\\
&=&\int_{\partial U_{x, z}}S^*_{x,z}(\Theta_{x,z})\\
&=&\sum_{y}\int_{U_{x,y}\times U_{y,z}}(S_{x,y}\times S_{y,z})^*
  (\Theta_{x,y}\times \Theta_{y,z})\\
&=&\sum_y <x,y><y,z>\\
&=&\sum_{y\in C_{k-1}} <x,y><y,z>,
\end{array} \leqno(6.24)$$
where the last equality comes from (6.17).
We finish the proof. $\Box$
\vskip 0.1in
\noindent
{\bf Definition 6.7: }{\it We define Floer homology $HF_*(X,H)$ as the homology of
chain complex $(C_*(X, H), \delta)$}
\vskip 0.1in
Since the action of Novikov ring $\Lambda_{\omega}$ is commutative with the boundary operation
$\delta$, Novikov ring acts on $HF_*(X, H)$ and we can view 
$HF_*(X,H)$ as a $\Lambda_{\omega}$-module.
\vskip 0.1in
\noindent
{\bf Remark 6.8: }{\it The boundary operator $\delta$ may depend on the choice of compatible
$\Theta_{x^-, x^+}$. However, Floer homology is independent of such a choice.}
\vskip 0.1in
\noindent
{\bf Proposition 6.9: }{\it $HF_*(X, H)$ is independent of $(H,J)$ and the 
construction of the  virtual neighborhood and the choice of compatible $\Theta_{x^-, x^+}$.}
\vskip 0.1in
  The proof is  routine. We leave it to the readers.
\vskip 0.1in
\noindent
{\bf Theorem 6.10: }{\it $HF_*(X, H)=H_*(X, \Lambda_{\omega})$ as a $\Lambda_{
\omega}$-module.}
\vskip 0.1in
\noindent
{\bf Corollary 6.11: }{\it Arnold conjecture holds for any symplectic manifold.}
\vskip 0.1in
 The basic idea is to view  $HF_*(X, H)$ and $H_*(X, \Lambda_{\omega})$ as the 
special cases of the Bott-type Floer homology \cite{RT3}, where $H_*(X, \Lambda_{\omega})$
is Floer homology of zero Hamiltonian function. The isomorphism between them
is interpreted as the independence of Bott-type Floer homology from Hamiltonian
functions.  Instead of giving 
the general construction of Bott type Floer homology, we shall construct the
isomorphism between $HF_*(X, H)$ and $H_*(X, \Lambda_{\omega})$ directly. It
 consists of several lemmas. 

 Let $\Omega_i(X)$ be the space of the differential
forms of degree $i$. Let $C_m(V, \Lambda_{\omega})=\oplus_{i+j=m}\Omega^{2n-i}(
X)\otimes \Lambda^j_{\omega}$, where we define $deg(e^A)=2C_1(X)(A)$.
For $\alpha\in \Omega^{2n-i}(X)$, define $\delta(\alpha)=d\alpha\in \Omega^{2n-
(i-1)}$.
The boundary operator is defined by
$$\delta(\alpha\otimes \lambda)=\delta(\alpha)\otimes \lambda \in C_{m-1}(V, \Lambda_{
\omega}).\leqno(6.25)$$
Clearly, its homology
$$H(C_*(V, \Lambda_{\omega}), \delta)=H_*(V, \Lambda_{\omega}).\leqno(6.26)$$
Consider a family of Hamiltonian function $H_s$ 
such that $H_s=0$ for $s<-1$ and $H_s=H$ for $s<1$. Furthermore, we choose 
a family of  compatible almost complex structures $J(s,x)$ such that $J_s=J$ for $s<-1$
 is $H$-admissible. Moreover, $J_s=J_0$ for $s>1$. Consider the moduli space of the solutions of
equation
$$\F((J_s), (H_s))=\frac{\partial u}{\partial s}+J(t,s, u(t,s))\frac{
\partial u}{\partial t}-\bigtriangledown H$$
$S^1\times (-\infty, +\infty)$ is conformal equivalent to $\C-0$ by the map
$$e^z: S^1\times (-\infty, +\infty)\rightarrow \C.\leqno(6.27)$$ 
Hence, we 
can view $u$ as map from $\C-\{0\}$ to $V$ which is holomorphic near zero. By
removable singularity theorem, $u$ extends to a map over $\C$ with a marked 
point at zero. In another words, $\lim_{s\rightarrow -\infty}u_s=pt$.
Furthermore,
when the energy $E(u)<\infty$, $u(s)$ converges to a periodic orbit when 
$s\rightarrow \infty$ by Lemma 6.3. Let $\M(pt,x^+)$ be the space of  $u$ such
that $\lim_{s\rightarrow \infty} u_s=x^+$. $\M(pt, x^+)$ has many components
of different dimensions. We use $\M(pt, A; x^+, u^+)$ to denote the components
satisfying $u\#u^+=A$. Consider the stable compactification $\overline{\M}(pt,
 A; x^+, u^+)$ in the same fashion. The virtual dimension of $\M(pt, A; x^+, u^+)$ is
$\mu(x^+, u^+)-2C_1(V)(A)$. Choose the stabilization terms $(s_{pt, A, x^+}, 
\Theta_{pt, A, x^+})$ compatible with the corner structure.
 Its virtual neighborhood
$(U(A; x^+, u^+), E(A; x^+,u^+), S(A; x^+, u^+))$ is a smooth V-manifold with corner.
 Notice
$$\partial(\overline{\B}(A; x^+, u^+))= \bigcup_{(x, u)}\overline{\B}(pt,A; x,u)
\times \overline{\B}((x,u); (x^+,u^+)).\leqno(6.28)$$
By our construction, 
$$\partial(U(A; x^+, u^+))\cong \bigcup_{(x,u)}E^{ot}_{U(A; x, u)\times U((x,u); (x^+,u^+))}
.\leqno(6.29)$$
Moreover, 
$$S(A, x^+, u^+)|_{\partial(U(A; x^+, u^+))}=\bigcup_{(x,u)}S(A; x, u)\times 
S((x,u); (x^+,u^+)),\leqno(6.30)$$
Let $e_{-\infty}$ be the evaluation map at $-\infty$.
 We define a map
$$\psi: C_m(V, \Lambda_{\omega})\rightarrow C_m(V, H)$$
by
$$\psi(\alpha, A;x^+, u^+)=\sum_{i=\mu(x^+, u^+)-2C_1(V)(A)}<\alpha,A; x^+, 
\mu^+> (x^+, u^+),\leqno(6.31)$$
where 
$$<\alpha, A; x^+, \mu^+>=\int_{U(A; x^+, u^+)}e^*_{-\infty}\alpha\wedge S(A; 
x^+, u^+)^*\Theta(A; x^+, u^+).\leqno(6.32)$$
\vskip 0.1in
\noindent
{\bf Lemma 6.12:}{\it (i) $\delta\psi=\psi\delta$.
\vskip 0.1in
\noindent
(ii)$\psi$ is independent of the 
virtual neighborhood compatible with the corner structure.}
\vskip 0.1in
{\bf Proof of Lemma:} The proof of (ii) is routine. We omit it. 

To prove (i),
for $\alpha\in \Omega^{2n-(i+1)}(X)$,
$$<\delta \alpha, A; x^+, \mu^+>=\int_{\partial U(A; 
x^+, u^+)}e^*_{-\infty}\alpha\wedge S(A; x^+, u^+)^*\Theta(A; x^+, u^+)\leqno(6.33)$$
$$=\sum_{(x,u)}\int_{
U(A; x,u)}e^*_{-\infty}(\alpha)\wedge S(A; x, u)^*\Theta(A; x, u)\int_{U((x,u); (x^+, u^+))}
S((x,u); (x^+,u^+))^*\Theta(x,u); (x^+,u^+)).$$
However, 
$$\dim (U(A; x, u))-deg(\Theta(A; x,u))=\mu(x,u)-2C_1(V)(A)<deg(\alpha)$$
unless $\mu(x,u)=\mu(x^+,u^+)+1$.
Hence, 
$$\begin{array}{lll}
&&\int_{\partial U(A; 
x^+, u^+)}\beta\wedge S(A; x^+, u^+)^*\Theta(A; x^+, u^+)\\
&=&\sum_{\mu(x,u)\mu(
x^+, u^+)+1}\int_{
U(A; x,u)}\alpha\wedge S(A; x, u)^*\Theta(A; x, u)\\
&&\int_{U((x,u); (x^+, u^+))}
S((x,u); (x^+,u^+))^*\Theta(x,u); (x^+,u^+))\\
&=&\psi\delta(x^+,u^+).
\end{array}\leqno(6.34)$$
$\Box$

Therefore, $\psi$ induces a homomorphism on Floer homology.

Consider a family of Hamiltonian function $H_s$ 
such that $H_s=0$ for $s>1$ and $H_s=H$ for $s<-1$. Furthermore, we choose 
a family of  compatible almost complex structures $J(s,x)$ such that $J_s=J$ for $s<-1$.
 Moreover, $J_s=J_0$ for $s>1$. Consider the moduli space of the solutions of
equation
$$\F((J_s), (H_s))=\frac{\partial u}{\partial s}+J(t,s, u(t,s))\frac{\partial u}{
\partial t}-\bigtriangledown H$$
$S^1\times (-\infty, +\infty)$ is conformal equivalent to $\C-0$ by the map
$$e^{-z}: S^1\times (-\infty, +\infty)\rightarrow \C.\leqno(6.35)$$ 
Hence, we 
can view $u$ as map from $\C-\{0\}$ to $V$ which is holomorphic near zero. By
removable singularity theorem, $u$ extends to a map over $\C$ with a marked 
point at zero. In another words, $\lim_{s\rightarrow \infty}u_s=pt$.
Furthermore,
when the energy $E(u)<\infty$, $u(s)$ converges to a periodic orbit when 
$s\rightarrow -\infty$ by Lemma 6.3. Let $\M(pt,x^-)$ be the space of  $u$ such
that $\lim_{s\rightarrow -\infty} u_s=x^-$. $\M(pt,x^-)$ has many components
of different dimension. We use $\M(x^-, u^-; pt,A)$ to denote the components
satisfying $u^-\#u=A$. The virtual dimension of $\M(x^-, u^-)$ is
$\mu(x^-, u^-)-2C_1(V)(A)$. Consider the stable compactification $\overline{\M}(
 x^-, u^-; pt,A)$ and its configuration space $\overline{\B}_{\delta}(x^-, u^-; pt, A)$.
Choose the stabilization terms $(s_{x^-; pt}, \Theta_{x^-, pt})$ 
compatible with the corner structure.  Furthermore, by adding more sections, we
can assume that the evaluation map $e_{\infty}$ is a submersion. Then, we
define
$$\phi: C_m(V,H)\rightarrow C_m(V, \Lambda_{\omega})$$
by 
$$\phi(x^-,u^-)=\sum_{A}<x^-, u^-; A>e^A.\leqno(6.36)$$
where
$$<x^-, u^-; A>=(e_{\infty})_*S(x^-,u^-;A)^*\Theta(x^-,u^-;A)\in \Omega^{2n-i}
(X)\leqno(6.37)$$
for $i=\mu(x^-, u^-)-2C_1(X)(A)$.
\vskip 0.1in
\noindent
{\bf Lemma 6.13:}{\it (i)$\phi\delta=\delta \phi$. (ii)$\phi$ is independent
of the choice of the virtual neighborhood compatible with the corner structure.}
\vskip 0.1in
{\bf Proof:} The proof of (i) is routine and we omit it. To prove (i),
$$\begin{array}{lll}
&&d<x^-, u^-; A>\\
&=&(e_{\infty})_*dS(x^-,u^-;A)^*\Theta(x^-,u^-;A)=(e_{\infty}|_
{\partial U(x^-, u^-; A)})_*S(x^-,u^-;A)^*\Theta(x^-,u^-;A)\\
&=&\sum_{\mu(x,u)=\mu(x^-,u^-)-1}(e_{\infty})_*S(x,u;A)^*\Theta(
x,u;A)\int_{U((x^-, u^-); (x, u))}S((x^-, u^-); (x, u))^*\Theta((x^-, u^-); (x, u))\\
&=&\phi\delta(x^-, u^-).
\end{array} \leqno(6.38)$$
$\Box$
\vskip 0.1in
\noindent
{\bf Lemma 6.14:}{\it $\phi\psi=Id$ and $\psi\phi=Id$ as the homomorphisms on Floer
homology.}
\vskip 0.1in
{\bf Proof: } The proof is tedious and routine. We omit it.

\section{Appendix}
   This appendix is due to B. Seibert \cite{S1}. We use the notation of the section 2.
\vskip 0.1in
\noindent
{\bf Lemma A1: }{\it Any local V-bundle of $\overline{\B}_A(Y,g, k)$ is dominated by
a global $V$-bundle.}
\vskip 0.1in
{\bf Proof: } The construction of global $V$-bundle imitates the similar construction
in algebraic geometry. First of all, we can slightly deform $\omega$ such that $[\omega]$
is a rational class. By taking multiple, we can assume $[\omega]$ is an integral class.
Therefore, it is Poincare dual to a complex line bundle $L$. We choose a unitary connection
$\bigtriangledown$ on $L$. There is a line bundle associated with the domain of stable
maps called dualized tangent sheaf $\lambda$. The restriction of $\lambda_C$ on $C$ is
$\lambda_C(x_1, \dots, x_k)$-the sheaf of meromorphic 1-form with simple pole at the
intersection points $x_1, \dots, x_k$. $\lambda_C$ varied continuously the domain of $f$.
For any $f\in \overline{\B}_A(Y,g,k)$, $f^*L$ is a line bundle over $dom(f)$
with a unitary connection. It is well-known in differential geometry that $f^*L$ has
a holomophic structure compatible with the unitary connection. Note that $L$ doesn't 
have holomoprhic structure in general. Therefore, $f^*L\otimes \lambda_C$ is a holomorphic
line bundle. Moreover, if $D$ is not a ghost component, $\omega(D)>0$ since it is
represented by a $J$-map. Therefore, $C_1(f^*L)(D)>0$. For ghost component, $\lambda_C$
is positive. By taking the higher power of $f^*L\otimes \lambda_C$, we can assume that
$f^*L\otimes\lambda_C$ is very ample.  Hence, $H^1(f^*L\otimes \lambda_C)=0$. Therefore,
$\E_f=H^0(f^*L\otimes \lambda_C)$ is of constant rank. It is easy to prove that $\E=\cup_f
\E_f$ is
bundle in terms of topology defined in Definition 3.10.

   To show that $\E$ dominates any local $V$-bundle, we recall that the group ring of
any finite group will dominate (or map surjectively to) any of its irreducible 
representation. So it is enough to construct a copy of group ring from $\E_f$. However,
$stb_f$ acts effectively on $dom(f)$. We can pick up a point $x\in dom(f)$ in the smooth
part of $dom(f)$ such that $stb_f$ acts on $x$ effectively. Then, $stb_f(x)$ is of
cardinality $|stb_f|$. By choose higher power of $f^*L\otimes \lambda_C$, we can assume that
there is a section $v\in \E_f$ such that $v(x)=1, v(g(v))=0$ for $g\in stb_f, g\neq id$.
Then, $stb_f(v)$ generates a copy of the group ring of $stb_f$.

\end{document}